\documentclass[10pt,reqno]{amsart}
\usepackage{mathptmx}
\usepackage{helvet}
\usepackage{courier}
\usepackage{amssymb}
\usepackage{amsmath,amsthm}
\usepackage{enumerate}
\usepackage{natbib}
\usepackage{hyperref}
\newcommand{\heads}{\mathit{h}} \newcommand{\tails}{\mathit{t}}  \newcommand{\reals}{\mathbb{R}} \newcommand{\nats}{\mathbb{N}}  \newcommand{\rv}{X} \newcommand{\values}{\mathcal{X}} \newcommand{\pr}{P} \newcommand{\lpr}{{\underline{\pr}}} \newcommand{\upr}{{\overline{\pr}}} \newcommand{\apr}{Q} \newcommand{\alpr}{{\underline{\apr}}}  \newcommand{\nex}{E} \newcommand{\lnex}{{\underline{\nex}}} \newcommand{\unex}{{\overline{\nex}}} \newcommand{\samp}{S} \newcommand{\lsamp}{{\underline{\samp}}} \newcommand{\usamp}{{\overline{\samp}}} \newcommand{\banlim}{L}   \newcommand{\gambles}{\mathcal{L}} \newcommand{\linprevs}{\mathbb{P}} \newcommand{\domain}{\mathcal{K}} \newcommand{\invars}{\mathcal{I}} \newcommand{\atoms}{\mathcal{A}} \newcommand{\atom}[2]{\left[#1\right]_{#2}} \newcommand{\transfos}{\mathcal{T}} \newcommand{\permuts}{\mathcal{P}} \newcommand{\desirs}{\mathcal{D}} \newcommand{\edesirs}{\mathcal{E}} \newcommand{\rdesirs}{\mathcal{R}} \newcommand{\cone}{\mathcal{C}} \newcommand{\lattice}{\mathcal{S}} \newcommand{\partit}{\mathcal{B}}  \newcommand{\counts}{\mathcal{N}} \newcommand{\sample}[1]{\mathbf{#1}} \newcommand{\cnt}[1]{\mathbf{#1}} \newcommand{\cntf}{\mathbf{T}} \newcommand{\set}[2]{\left\{#1\colon#2\right\}}   \newcommand{\iden}{\operatorname{id}} \newcommand{\abs}[1]{\lvert#1\rvert} \newcommand{\norm}[1]{\lVert#1\rVert} \newcommand{\asa}{\Leftrightarrow} \newcommand{\solp}{\mathcal{M}} \newcommand{\dif}{\,\mathrm{d}} \newcommand{\indif}{\approx} \newcommand{\pref}{\succeq} \newcommand{\spref}{\succ} \newcommand{\incomp}{\parallel} \newcommand{\succeeds}{\geqslant} \newcommand{\citegen}[2][]{\citeauthor*{#2}'s \citeyearpar[#1]{#2}} 
\theoremstyle{plain} 
\newtheorem{theorem}{Theorem} 
\newtheorem{proposition}[theorem]{Proposition}
\newtheorem{corollary}[theorem]{Corollary} 
 
\newtheorem*{RPIR}{Revised Principle of Insufficient Reason (RPIR)}

\newtheorem*{SSP}{Strong Symmetry Principle (SSP)} 
\newtheorem*{SP}{Symmetry Principle (SP)} 
\newtheorem*{EP}{Embedding Principle (EP)} 
\theoremstyle{remark} 
\newtheorem{example}{Example} 
\newtheorem*{dice}{The dice example} 
\newtheorem*{dice-con}{The dice example (cont.)} 
\newtheorem*{shifts}{The example of shift transformations}
\newtheorem*{shifts-con}{The example of shift transformations (cont.)}

\theoremstyle{definition} 
\newtheorem{definition}{Definition}

\begin{document}
\title{Symmetry of models versus models of symmetry} 
\author{Gert de Cooman} 
\address{Ghent   University, SYSTeMS Research Group, Technologiepark -- Zwijnaarde 914, 9052 Zwijnaarde, Belgium} \email{gert.decooman@ugent.be} 
\author{Enrique Miranda} 
\address{Rey Juan Carlos University, Dep. of   Statistics and Operations Research, C-Tulip\'an, s/n 28933, M\'ostoles, Spain} 
\email{enrique.miranda@urjc.es}

\begin{abstract}
  A model for a subject's beliefs about a phenomenon may exhibit symmetry, in the sense that it is   invariant under certain transformations.  On the other hand, such a belief model may be intended to   represent that the subject believes or knows that the phenomenon under study exhibits symmetry.  We   defend the view that these are fundamentally different things, even though the difference cannot be   captured by Bayesian belief models. In fact, the failure to distinguish between both situations   leads to Laplace's so-called Principle of Insufficient Reason, which has been criticised   extensively in the literature.
\par
We show that there are belief models (imprecise probability models, coherent lower previsions) that generalise and include the Bayesian belief models, but where this fundamental difference can be captured. This leads to two notions of symmetry for such belief models: weak invariance (representing symmetry of beliefs) and strong invariance (modelling beliefs of symmetry).  We discuss various mathematical as well as more philosophical aspects of these notions. We also discuss a few examples to show the relevance of our findings both to probabilistic modelling and to statistical inference, and to the notion of exchangeability in particular.
\end{abstract}

\date{19 April 2006}
\maketitle

\section{Introduction}\label{sec:intro}
This paper deals with symmetry in relation to models of beliefs. Consider a model for a subject's beliefs about a certain phenomenon. Such a \emph{belief model} may be \emph{symmetrical}, in the sense that it is invariant under certain transformations. On the other hand, a belief model may try to capture that the subject believes that the phenomenon under study exhibits symmetry, and we then say that the belief model \emph{models symmetry}. We defend the view that there is an important conceptual difference between the two cases: symmetry of beliefs should not be confused with beliefs of symmetry.\footnote{This echoes \citegen[Section~9.5.6, p.~466]{walley1991} view that `symmetry of   evidence' is not the same thing as `evidence of symmetry'.}
\par
Does this view need defending at all? That there is a difference may strike you as obvious, and yet we shall argue that Bayesian belief models, which are certainly the most popular belief models in the literature, are unable to capture this difference.
\par
To make this clearer, consider a simple example. Suppose I will toss a coin, and you are ignorant about its relevant properties: it might be fair but on the other hand it might be heavily loaded, or it might even have two heads, or two tails (situation $A$). To you the outcomes of the toss that are practically possible are $\heads$ (for heads) and $\tails$ (for tails). Since you are ignorant about the properties of the coin, any model for your beliefs should not change if heads and tails are permuted, so the model that `faithfully' captures your beliefs about the outcome of the toss should be symmetrical too, i.e., invariant under this permutation of heads and tails.
\par
Suppose on the other hand that you know that the coin (and the tossing mechanism) I shall use is completely symmetrical (situation $B$). Your belief model about the outcome of the toss should capture this knowledge, i.e., it should model your beliefs about the symmetry of the coin.
\par
Our point is that belief models should be able to catch the important difference between your beliefs in the two situations. Bayesian belief models cannot do this. Indeed---the argument is well-known---the only symmetrical probability model, which is in other words invariant under permutations of heads and tails, assigns equal probability 1/2 to heads and tails. But this is automatically also the model that captures your beliefs that the coin is actually symmetrical, so heads and tails should be equally likely.
\par
The real reason why Bayesian belief models cannot capture the difference between symmetry of models and modelling symmetry, is that they do not allow for \emph{indecision}. Suppose that I ask you to express your preferences between two gambles, whose reward depends on the outcome of the toss. For first one, $a$, you will win one euro if the outcome is heads, and lose one if it is tails. The second one, $b$, gives the same rewards, but with heads and tails swapped.
\par
In situation $B$, because you believe the coin to be symmetrical, it does not matter to you which gamble you get, and you are \emph{indifferent} in your choice between the two.
\par
But in situation $A$, on the other hand, because you are completely ignorant about the coin, the available information gives you \emph{no reason to (strictly) prefer} $a$ over $b$ or $b$ over $a$. You are therefore \emph{undecided} about which of the two gambles to choose.
\par
Because decision based on Bayesian belief models leaves you no alternative but to either strictly prefer one action over the other, or to be indifferent between them, the symmetry of the model leaves you \emph{no choice but to act as if you were indifferent} between $a$ and $b$. We strongly believe that it is wrong to confuse indecision with indifference in this example (and elsewhere of course), but Bayesian belief models leave you no choice but to do so, unless you want to let go of the principle that if your evidence or your beliefs are symmetrical, your belief model should be symmetrical as well.  The problem with Laplace's Principle of Insufficient Reason is precisely this: if you use a Bayesian probability model then the symmetry present in ignorance forces you to treat indecision (or insufficient reason to decide) between $a$ and $b$ as if it were indifference.\footnote{This may seem a good explanation why \citet[p.~83]{keynes1921} renamed the   `Principle of Insufficient Reason' the `Principle of Indifference'. He (and others, see   \cite{zabell1989a}) also suggested that the principle should not be applied in a state of complete   ignorance, but only if there is good reason to justify the indifference (such as when there is   evidence of symmetry). By the way, Keynes was also among the first to consider what we shall call   imprecise probability models, as his comparative probability relations were not required to be   complete.} Or in other words, it forces you to treat symmetry of beliefs as if there were beliefs of symmetry.
\par
If on the other hand, we consider belief models that allow for indecision, we can sever the unholy link between indecision and indifference, because in a state of complete ignorance, we are then allowed to remain undecided about which of the two actions to choose: in the language of preference relations, they simply become \emph{incomparable}, and you need not be indifferent between them. As we shall see further on, similar arguments show that such belief models also allow us to distinguish between `symmetry of models' and `models of symmetry' in those more general situations where the symmetry involved is not necessarily that which goes along with complete ignorance.
\par
So, it appears that in order to better understand the interplay between modelling beliefs and issues of symmetry, which is the main aim of this paper, we shall need to work with a language, or indeed, with a type of belief models that, unlike the Bayesian ones, take indecision seriously. For this purpose, we shall use the language of the so-called \emph{imprecise probability models} \citep{walley1991}, and in particular coherent lower previsions, which have the same behavioural pedigree as the more common Bayesian belief models (\textit{in casu} coherent previsions, see \citet{finetti19745}), and which contain these models as a special case. We give a somewhat unusual introduction to such models in Section~\ref{sec:ip-models}.\footnote{For other brief and perhaps more   conventional introductions to the topic, we refer to   \Citet{walley1996,cooman2004b,cooman2003c,cooman2004a}. A much more detailed account of the   behavioural theory of imprecise probabilities can be found in \citet{walley1991}.} In Section~\ref{sec:transformations}, we provide the necessary mathematical background for discussing symmetry: we discuss monoids of transformations, and invariance under such monoids. After these introductory sections, we start addressing the issue of symmetry in relation to belief models in Section~\ref{sec:belief-model-invariance}. We introduce two notions of invariance for the imprecise probability models introduced in Section~\ref{sec:ip-models}: \emph{weak invariance}, which captures symmetry of belief models, and \emph{strong invariance}, which captures that a model represents the belief that there is symmetry. We study relevant mathematical properties of these invariance notions, and argue that the distinction between them is very relevant when dealing with symmetry in general, and in particular (Section~\ref{sec:complete-ignorance}) for modelling complete ignorance.  Further interesting properties of weak and strong invariance, related to inference, are the subject of Sections~\ref{sec:weakly-invariant-lpr} and~\ref{sec:strongly-invariant-lpr}, respectively. We show among other things that a weakly invariant coherent lower prevision can always be extended to a larger domain, in a way that is as conservative as possible. This implies that, for any given monoid of transformations, there always are weakly invariant coherent lower previsions.  This is not generally the case for strong invariance, however, and we give and discuss sufficient conditions such that for a given monoid of transformations, there would be strongly invariant coherent (lower) previsions.  We also give various expression for the smallest strongly invariant coherent lower prevision that dominates a given weakly invariant one (if it exists). In Section~\ref{sec:shift-invariance}, we turn to the important example of coherent (lower) previsions on the set of natural numbers, that are shift-invariant, and we use them to characterise the strongly invariant coherent (lower) previsions on a general space provided with a single transformation. Further examples are discussed in Section~\ref{sec:finite}, where we characterise weak and strong invariance with respect to finite groups of permutations. In particular, we discuss \citegen{walley1991} generalisation to lower previsions of \citegen{finetti1937} notion of exchangeability, and we use our characterisation of strong permutation invariance to prove a generalisation to lower previsions of de Finetti's representation results for finite sequences of exchangeable random variables. Conclusions are gathered in Section~\ref{sec:conclusions}.
\par
We want to make it clear at this point that this paper owes a significant intellectual debt to Peter Walley. First of all, we use his behavioural imprecise probability models \citep{walley1991} to try and clarify the distinction between symmetry of beliefs and beliefs of symmetry. Moreover, although we like to believe that much of what we do here is new, we are also aware that in many cases we take to their logical conclusion a number of ideas about symmetry that are clearly present in his work (mainly \citet[Sections~3.5, 9.4 and~9.5]{walley1991} and \cite{pericchi1991}), sometimes in embryonic form, and often more fully worked out.

\section{Imprecise probability models}\label{sec:ip-models}
Consider a very general situation in which uncertainty occurs: a subject is
uncertain about the value that a variable $\rv$ assumes in a set of possible
values $\values$. Because the subject is uncertain, we shall call $\rv$ an
\emph{uncertain}, or \emph{random}, variable.
\par
The central concept we shall use in order to model our subject's uncertainty
about $\rv$, is that of a \emph{gamble} (on $\rv$, or on $\values$), which is
a bounded real-valued function $f$ on $\values$. In other words, a gamble $f$
is a map from $\values$ to the set of real numbers $\reals$ such that
\begin{equation*}
  \sup f:=\sup\set{f(x)}{x\in\values}
  \text{ and }
  \inf f:=\inf\set{f(x)}{x\in\values}
\end{equation*}
are (finite) real numbers.  It is interpreted as the reward function for a
transaction which may yield a different (and possibly negative) reward $f(x)$,
measured in units (called \emph{utiles}) of a pre-determined linear
utility,\footnote{This utility can be regarded as amounts of money, as is the
  case for instance in \citet{finetti19745}. It is perhaps more realistic, in
  the sense that the linearity of the scale is better justified, to interpret
  it in terms of probability currency: we win or lose lottery tickets
  depending on the outcome of the gamble; see
  \citet[Section~2.2]{walley1991}.} for each of the different values $x$ that
the random variable $\rv$ may assume in $\values$.
\par
We denote the set of all gambles on $\rv$ by $\gambles(\values)$. For any two
gambles $f$ and $g$, we denote their point-wise sum by $f+g$, and we denote
the point-wise (scalar) multiplication of $f$ with a real number $\lambda$ by
$\lambda f$. $\gambles(\values)$ is a real linear space under these
operations. We shall always endow this space with the \emph{supremum norm},
i.e., $\norm{f}=\sup\abs{f}=\sup\set{\abs{f(x)}}{x\in\values}$, or
equivalently, with the topology of uniform convergence, which turns
$\gambles(\values)$ into a Banach space.
\par
An \emph{event} $A$ is a subset of $\values$. If $\rv\in A$ then we say that
the event \emph{occurs}, and if $\rv\not\in A$ then we say that $A$
\emph{doesn't occur}, or equivalently, that the \emph{complement(ary event)}
$A^c=\set{x\in\values}{x\not\in A}$ occurs. We shall identify an event with a
special $\{0,1\}$-valued gamble $I_A$, called its \emph{indicator}, and
defined by $I_A(x)=1$ if $x\in A$ and $I_A(x)=0$ elsewhere. We shall often
write $A$ for $I_A$, whenever there is no possibility of confusion.

\subsection{Coherent sets of really desirable gambles}
\label{sec:really-desirable}
Given the information that the subject has about $\rv$, she will be disposed
to accept certain gambles, and to reject others. The idea is that we model a
subject's beliefs about $\rv$ by looking at which gambles she accepts, and to
collect these into a \emph{set of really desirable gambles} $\rdesirs$.

\begin{dice}
  Assume that our subject is uncertain about the outcome $\rv$ of my tossing a
  die. In this case $\values=\values_6:=\{1,2,3,4,5,6\}$ is the set of
  possible values for $\rv$. If the subject is rational, she will accept the
  gamble which yields a positive reward whatever the value of $\rv$, because
  she is certain to improve her `fortune' by doing so. On the other hand, she
  will not accept a non-positive gamble that is negative somewhere, because by
  accepting such a gamble she can only lose utility (we then say she
  \emph{incurs a partial loss}).  She will not accept the gamble which makes
  her win one utile if the outcome $\rv$ is 1, and makes her lose five utiles
  otherwise, unless she knows for instance that the die is loaded very heavily
  in such a way that the outcome $1$ is almost certain to come up.
  \par
  Real desirability can also be interpreted in terms of the betting behaviour
  of our subject. Suppose she wants to bet on the occurrence of some event,
  such as my throwing $1$ (so that she receives 1 utile if the event happens
  and 0 utiles otherwise). If she thinks that the die is fair, she should be
  disposed to bet on this event at any rate $r$ strictly smaller than
  $\frac{1}{6}$. This means that the gamble $I_{\{1\}}-r$ representing this
  transaction (winning $1-r$ if the outcome of $\rv$ is $1$ and losing $r$
  otherwise) will be really desirable to her for $r<\frac{1}{6}$.
  $\blacklozenge$
\end{dice}

Now, accepting certain gambles has certain consequences, and has certain
implications for accepting other gambles, and if our subject is rational,
which we shall assume her to be, she should take these consequences and
implications into account. To give but one example, if our subject accepts a
certain gamble $f$ she should also accept any other gamble $g$ such that
$g\geq f$, i.e., such that $g$ \emph{point-wise dominates} $f$, because
accepting $g$ is certain to bring her a reward that is at least as high as
accepting $f$ does.
\par
Actually, this requirement is a consequence [combine (D2) with (D3)] of the
following four basic rationality axioms for real desirability, which we shall
assume any rational subject's set of really desirable gambles $\rdesirs$ to
satisfy:
\begin{enumerate}[(D1)]
\item if $f<0$ then $f\not\in\rdesirs$ [avoiding partial loss];
\item if $f\geq0$ then $f\in\rdesirs$ [accepting sure gains];
\item if $f\in\rdesirs$ and $g\in\rdesirs$ then $f+g\in\rdesirs$ [accepting
  combined gambles]
\item if $f\in\rdesirs$ and $\lambda>0$ then $\lambda f\in\rdesirs$ [scale
  invariance].
\end{enumerate}
where $f<g$ is shorthand for $f\leq g$ and $f\not=g$.\footnote{So, here and in
  what follows, we shall write `$f<0$' to mean `$f\leq0$ and not $f=0$', and
  `$f>0$' to mean `$f\geq0$ and not $f=0$'.} We call any subset $\rdesirs$ of
$\gambles(\values)$ that satisfies these axioms a \emph{coherent} set of
really desirable gambles.
\par
It is easy to see that these axioms reflect the behavioural rationality of our
subject: (D1) means that she should not be disposed to accept a gamble which
makes her lose utiles, no matter the outcome; (D2) means that she should
accept a gamble which never makes her lose utiles; on the other hand, if she
is disposed to accept two gambles $f$ and $g$, she should also accept the
combination of the two gambles, which leads to a reward $f+g$; this is an
immediate consequence of the linearity of the utility scale.  This justifies
(D3). And finally, if she is disposed to accept a gamble $f$, she should be
disposed to accept the scaled gamble $\lambda f$ for any $\lambda>0$, because
this just reflects a change in the linear utility scale. This is the idea
behind condition (D4).

\par
\citet{walley1991,walley2000} has a further coherence axiom that sets of
really desirable gambles should satisfy, which turns out to be quite important
for conditioning, namely
\begin{enumerate}
\item[(D5)] if $\partit$ is a partition of $\values$ and if $I_Bf\in\rdesirs$
  for all $B$ in $\partit$, then $f\in\rdesirs$ [full conglomerability].
\end{enumerate}
Since this axiom is automatically satisfied whenever $\values$ is finite [it
is then an immediate consequence of~(D3)], and since we shall not be concerned
with conditioning unless when $\values$ is finite (see
Section~\ref{sec:finite}), we shall ignore this additional axiom in the
present discussion.
\par
A coherent set of really desirable gambles is a convex cone
[axioms~(D3)--(D4)] that includes the `non-negative orthant'
$\cone_+:=\set{f\in\gambles(\values)}{f\geq0}$ [axiom~(D2)] and has no gamble
in common with the `negative orthant'
$\cone_-:=\set{f\in\gambles(\values)}{f<0}$ [axiom~(D1)].\footnote{This means
  that the zero gamble $0$ belongs to the set of really desirable gambles.
  This is more a mathematical convention than a behavioural requirement, since
  this gamble has no effect whatsoever in the amount of utiles of our subject.
  See more details in \citet{walley1991}.} If we have two coherent sets of
really desirable gambles $\rdesirs_1$ and $\rdesirs_2$, such that
$\rdesirs_1\subseteq\rdesirs_2$, then we say that $\rdesirs_1$ is less
committal, or more conservative, than $\rdesirs_2$, because a subject whose
set of really desirable gambles is $\rdesirs_2$ accepts at least all the
gambles in $\rdesirs_1$. The least-committal (most conservative, smallest)
coherent set of really desirable gambles is $\cone_+$. Within this theory, it
seems to be the appropriate model for \emph{complete ignorance}: if our
subject has no information at all about the value of $\rv$, she should be
disposed to accept only those gambles which cannot lead to a loss of utiles
(see also the discussion in Section~\ref{sec:complete-ignorance}).
\par
Now suppose that our subject has specified a set $\rdesirs$ of gambles that
she accepts. In an elicitation procedure, for instance, this would typically
be a finite set of gambles, so we cannot expect this set to be coherent. We
are then faced with the problem of enlarging this $\rdesirs$ to a coherent set
of really desirable gambles that is as small as possible: we want to find out
what are the (behavioural) consequences of the subject's accepting the gambles
in $\rdesirs$, taking into account \emph{only} the requirements of coherence.
This inference problem is (also formally) similar to the problem of inference
(logical closure) in classical propositional logic, where we want to find out
what are the consequences of accepting certain propositions.\footnote{See
  \citet{moral1995} and \Citet{cooman2001b,cooman2003a} for more details on
  this connection between natural extension and inference in classical
  propositional logic.}
\par
The smallest convex cone including $\cone_+$ and $\rdesirs$, or in other
words, the smallest subset of $\gambles(\values)$ that includes $\rdesirs$ and
satisfies (D2)--(D4), is given by
\begin{equation*}
  \edesirs^r_\rdesirs:=
  \set{g\in\gambles(\values)}
  {\text{$g\geq\sum_{k=1}^n\lambda_k f_k$ for some $n\geq0$,
      $\lambda_k\in\reals^+$ and $f_k\in\rdesirs$}},
\end{equation*}
where $\reals^+$ denotes the set of non-negative real numbers.  If this convex
cone $\edesirs^r_\rdesirs$ intersects $\cone_-$ then it is easy to see that
actually $\edesirs^r_\rdesirs=\gambles(\values)$, and then it is impossible to
extend $\rdesirs$ to a coherent set of really desirable gambles [because (D1)
cannot be satisfied]. Observe that $\edesirs^r_\rdesirs\cap\cone_-'=\emptyset$
if and only if
\begin{equation*}
  \text{there are no $n\geq0$, $\lambda_k\in\reals^+$ and $f_k\in\rdesirs$
    such that $\sum_{k=1}^n\lambda_k f_k<0$},
\end{equation*}
and we then say that the set $\rdesirs$ \emph{avoids partial loss}. Let us
interpret this condition. Assume that it doesn't hold (so we say that
$\rdesirs$ \emph{incurs partial loss}). Then there are really desirable
gambles $f_1$, \dots, $f_n$ and positive $\lambda_1$, \dots, $\lambda_n$ such
that $\sum_{k=1}^{n}\lambda_kf_k<0$.  But if our subject is disposed to accept
the gamble $f_k$ then by coherence [axioms(D2) and~(D4)] she should also be
disposed to accept the gamble $\lambda_kf_k$ for all $\lambda_k\geq0$.
Similarly, by coherence [axiom~(D3)] she should also be disposed to accept the
sum $\sum_{k=1}^{n}\lambda_kf_k$. Since this sum is non-positive, and strictly
negative in at least some elements of $\values$, we see that the subject can
be made subject to a partial loss, by suitably combining gambles which she
accepts. This is unreasonable.
\par
When the class $\rdesirs$ avoids partial loss, and only then, we are able to
extend $\rdesirs$ to a coherent set of really desirable gambles, and the
smallest such set is precisely $\edesirs^r_\rdesirs$, which is called the
\emph{natural extension} of $\rdesirs$ to a set of really desirable gambles.
This set reflects only the behavioural consequences of the assessments present
in $\rdesirs$: the acceptance of a gamble $f$ not in $\edesirs^r_\rdesirs$
(or, equivalently, a set of really desirable gambles strictly including
$\edesirs^r_\rdesirs$) is not implied by the information present in
$\rdesirs$, and therefore represents stronger implications that those of
coherence alone.

\subsection{Coherent sets of almost-desirable gambles}
\label{sec:almost-desirable}
Coherent sets of really desirable gambles constitute a very general and
powerful class of models for a subject's beliefs (see
\citet[Appendix~F]{walley1991} and \citet{walley2000} for more details and
discussion). We could already discuss symmetry aspects for such coherent sets
of really desirable gambles, but we shall instead concentrate on a slightly
less general and powerful type of belief models, namely coherent lower and
upper previsions. Our main reason for doing so is that this will allow us to
make a more direct comparison to the more familiar Bayesian belief models, and
in particular to \citegen{finetti19745} coherent previsions, or fair prices.
\par
Consider a gamble $f$. Then our subject's \emph{lower prevision}, or supremum
acceptable buying price, $\lpr(f)$ for $f$ is defined as the largest real
number $s$ such that she accepts the gamble $f-t$ for any price $t<s$, or in
other words accepts to buy $f$ for any such price $t$. Similarly, her
\emph{upper prevision}, or infimum acceptable selling price, $\upr(f)$ for the
gamble $f$ is the smallest real number $s$ such that she accepts the gamble
$t-f$ for any price $t>s$, or in other words accepts to sell $f$ for any such
price $t$.
\par
For an event $A$, the lower prevision $\lpr(I_A)$ of its indicator is also
called the \emph{lower probability} of $A$, and denoted by $\lpr(A)$. It can
be interpreted as the supremum rate for betting on the event $A$. Similarly,
$\upr(I_A)$ is called the \emph{upper probability} of $A$, and also denoted by
$\upr(A)$.
\par
Since selling a gamble $f$ for price $s$ is the same thing as buying $-f$ for
price $-s$, we have the following \emph{conjugacy} relationship between an
upper and a lower prevision:
\begin{equation*}
  \upr(f)=-\lpr(-f).
\end{equation*}
This implies that from a given lower prevision $\lpr$, we can always construct
the conjugate upper prevision $\upr$, so they are mathematically equivalent
belief models. In what follows, we shall mainly concentrate on lower
previsions.
\par
Now assume that our subject has a coherent set of really desirable gambles
$\rdesirs$, then it is clear from the definition of lower and upper prevision
that we can use $\rdesirs$ to define a lower prevision
\begin{equation}\label{eq:lpr}
  \lpr_\rdesirs(f)=\sup\set{s\in\reals}{f-s\in\rdesirs}\tag{D-LPR}
\end{equation}
and an upper prevision
\begin{equation*}
  \upr_\rdesirs(f)=\inf\set{s\in\reals}{s-f\in\rdesirs}
\end{equation*}
for every gamble $f$ on $\values$. So, given $\rdesirs$ we can construct two
real-valued functionals, $\lpr_\rdesirs$ and $\upr_\rdesirs$, whose
interpretation is that of a supremum acceptable buying price, and an infimum
acceptable selling price, respectively, and whose domain is
$\gambles(\values)$. We shall call these functionals \emph{lower} and
\emph{upper previsions}.
\par
We call a \emph{coherent lower prevision} on $\gambles(\values)$ any
real-valued functional on $\gambles(\values)$ satisfying the following three
axioms:
\begin{enumerate}[(P1)]
\item $\lpr(f)\geq\inf f$ [accepting sure gains]; \item
  $\lpr(f+g)\geq\lpr(f)+\lpr(g)$ [super-additivity]; \item $\lpr(\lambda
  f)=\lambda\lpr(f)$ [non-negative homogeneity].
\end{enumerate}
for all gambles $f$ and $g$ on $\values$, and all non-negative real $\lambda$.
\par
It follows from the coherence axioms (D1)--(D4) for $\rdesirs$ that the lower
prevision $\lpr_\rdesirs$ that corresponds to a coherent set of really
desirable gambles $\rdesirs$ is coherent.\footnote{To prove (P1), use (D2);
  for (P2) use (D3); and for (P3) use (D4) for $\lambda>0$ and (D1) and (D2)
  for $\lambda=0$.}
\par
So we see that with a coherent set of really desirable gambles $\rdesirs$, we
can define a coherent lower prevision on $\gambles(\values)$,
using~\eqref{eq:lpr}. We shall see further on that, conversely, given a
coherent lower prevision $\lpr$ on $\gambles(\values)$, we can always find a
coherent set of really desirable gambles $\rdesirs$ such that $\lpr$ and
$\rdesirs$ are related through~\eqref{eq:lpr}. But unfortunately, the
relationship between the two types of belief models is many-to-one: there are
usually many coherent sets of really desirable gambles that lead to the same
coherent lower prevision. This is why we said before that coherent sets of
really desirable gambles are a more general and powerful belief model than
coherent lower previsions.  The ultimate reason for this is the following:
suppose that a subject specifies her supremum buying price $\lpr(f)$ for a
gamble $f$. This implies that she accepts all the gambles $f-\lpr(f)+\delta$,
where $\delta>0$.  But the specification of $\lpr(f)$ says nothing about the
gamble $f-\lpr(f)$ (where $\delta=0$) itself: she might accept it, but then
again she might not.  And precisely because specifying a coherent lower
prevision says nothing about this border behaviour, it leads to a belief model
that is less powerful than coherent sets of really desirable gambles, where
this border behaviour would be determined.

\begin{dice-con}
  Let us go back to the die example. Consider, for any $x$ in
  $\values_6=\{1,\dots,6\}$, the event $\{x\}$ that the outcome $\rv$ of
  rolling the die is $x$. If, for some real number $r$, our subject accepts
  the gamble $I_{\{x\}}-r$, she is willing to pay $r$ utiles in return for the
  uncertain reward $I_{\{1\}}$, or in other words to bet \emph{on} the event
  $\{1\}$ at \emph{rate} $r$. So her lower probability $\lpr(\{x\})$ for
  $\{x\}$, or equivalently, her lower prevision $\lpr(I_{\{x\}})$ for
  $I_{\{x\}}$, is the supremum rate at which she is willing to bet on $\{x\}$.
  This means that she accepts the gamble $I_{\{x\}}-s$ for any
  $s<\lpr(\{x\})$. But it doesn't imply that she actually accepts the gamble
  $I_{\{x\}}-\lpr(\{x\})$: this gamble is only claimed to be almost-desirable,
  as we shall see further on.
  \par
  If she is completely ignorant about the properties of the die, her evidence
  about the die is symmetrical, i.e., doesn't change when the possible
  outcomes are permuted. A belief model that `faithfully' captures the
  available evidence should therefore be symmetrical with respect to such
  permutations as well, so we infer that in particular $\lpr(\{1\})$, \dots,
  $\lpr(\{6\})$ are all equal to some number $p$. Coherence [use (P1) and
  (P2)] then requires that $0\leq p\leq\frac{1}{6}$. Any such $p$ leads to a
  symmetrical lower probability defined on the singletons, and therefore
  reflects `symmetry of beliefs'. As we have indicated above, the model
  corresponding to $p=0$ is the one that reflects complete ignorance. We shall
  see further on (see Sections~\ref{sec:strong-invariance}
  and~\ref{sec:finite}) that the choice $p=\frac{1}{6}$ leads to the only
  model that captures the belief that the die is fair, i.e., that reflects
  `beliefs of symmetry'.  $\blacklozenge$
\end{dice-con}
\par
In order to better understand the relationship between coherent lower
previsions and coherent sets of really desirable gambles, we need to
introduce, besides \emph{real} desirability, an new and weaker notion, called
\emph{almost-desirability}, which will also play an important part in our
discussion of symmetry further on. This notion is inspired by the ideas in the
discussion above: we say that a gamble $f$ is \emph{almost-desirable} to a
subject, or that she \emph{almost-accepts} $f$, whenever she accepts
$f+\delta$, or in other words $f+\delta$ is really desirable to her, for any
strictly positive amount of utility $\delta>0$. By stating that $f$ is
almost-desirable to her, nothing is specified about whether the subject
accepts $f$ itself: she might, but then again she also might not. If we
generically denote by $\desirs$ a set of gambles that are almost-desirable to
our subject, we see that the set $\desirs_\rdesirs$ of almost-desirable
gambles that corresponds to a coherent set $\rdesirs$ of really desirable
gambles, is given by
\begin{equation}\label{eq:ad}
  \desirs_\rdesirs
  =\set{f\in\gambles(\values)}{(\forall\delta>0)f+\delta\in\rdesirs}
  =\bigcap_{\delta>0}[\rdesirs-\delta]\tag{D-M} 
\end{equation}
so $\desirs_\rdesirs$ is the closure (in the topology of uniform convergence
on $\gambles(\values)$) of the convex cone $\rdesirs$.
\par
We call any set of gambles $\desirs$ that satisfies the following five axioms
a \emph{coherent set of almost-desirable gambles}:
\begin{enumerate}[(M1)]
\item if $\sup f<0$ then $f\not\in\desirs$ [avoiding sure loss];
\item if $\inf f\geq0$ then $f\in\desirs$ [accepting sure gains];
\item if $f\in\desirs$ and $g\in\desirs$ then $f+g\in\desirs$ [accepting
  combined gambles];
\item if $f\in\desirs$ and $\lambda>0$ then $\lambda f\in\desirs$ [scale
  invariance];
\item if $f+\delta\in\desirs$ for all $\delta>0$ then $f \in \desirs$
  [closure].
\end{enumerate}
It is a closed and convex cone in $\gambles(\values)$ that includes the
non-negative orthant $\cone_+$ and does not intersect with the set
$\cone_{-}'=\set{f\in\gambles(\values)}{\sup f<0}\subset\cone_-$.  It is easy
to see that the set of almost-desirable gambles $\desirs_\rdesirs$ that
corresponds to a coherent set of really desirable gambles $\rdesirs$ is
actually also coherent.\footnote{To prove (M1), use (D1) with
  $\delta=-\frac{\sup f}{2}$; to prove (M2), use (D2); to prove (M3), use
  (D3); to prove (M4), use (D4); and to prove (M5), use
  $\epsilon=\frac{\delta}{2}$ and the definition of $\desirs_\rdesirs$ to
  prove that $f+\delta\in\rdesirs$ for all $\delta>0$.}
\par
It should at this point come as no surprise that coherent lower previsions and
coherent sets of almost-desirable gambles are actually equivalent belief
models.  Indeed, consider a coherent set of almost-desirable gambles
$\desirs$, i.e., $\desirs$ satisfies (M1)--(M5). Then the real-valued
functional $\lpr_\desirs$ defined on $\gambles(\values)$ by\footnote{The
  supremum in~Eq.~\eqref{eq:lpr} now becomes a maximum, simply because the set
  $\desirs$ is closed.}
\begin{equation}\label{eq:desirs-lpr}
  \lpr_\desirs(f):=\max\set{s\in\reals}{f-s\in\desirs}\tag{M-LPR}
\end{equation}
satisfies (P1)--(P3) and therefore is a coherent lower prevision on
$\gambles(\values)$.\footnote{(P1) follows from (M2), (P2) from (M3) and (P3)
  is a consequence of (M4).}
\par
Conversely, if we consider a coherent lower prevision $\lpr$ on
$\gambles(\values)$, i.e., $\lpr$ satisfies (P1)--(P3), then the set of
gambles
\begin{equation}\label{eq:lpr-desirs}
  \desirs_\lpr:=\set{f\in\gambles(\values)}{\lpr(f)\geq0}\tag{LPR-M}
\end{equation}
satisfies~(M1)--(M5) and is therefore a coherent set of almost-desirable
gambles.\footnote{First, conditions (P1) and (P2) imply that $\lpr$ is
  monotone. Now, (P2) and (P3) imply that $0=\lpr(0)\geq \lpr(f)+\lpr(-f)\geq
  \lpr(f)+\inf (-f)$, whence $\lpr(f)\leq\sup f$. From these two facts we
  deduce (M1).  (M2) is a consequence of (P1), (M3) of (P2) and (M4) of (P3).
  Finally, the monotonicity of $\lpr$ implies that $\lpr(\mu)=\mu$ for any
  constant value $\mu$, and from this we deduce that
  $\lpr(f+\delta)=\lpr(f)+\delta$ for any $\delta>0$. This implies (M5).}
Moreover, the relationships~\eqref{eq:desirs-lpr} and~\eqref{eq:lpr-desirs}
are bijective (one-to-one and onto), and they are each other's
inverses.\footnote{To see that they are each other inverses, it suffices to
  use that a coherent lower prevision satisfies $\lpr(f-s)=\lpr(f)-s$ for any
  gamble $f$ and any real number $s$, and, conversely, that $f\in\desirs_\lpr$
  if and only if $\lpr(f)\geq0$; this implies also that both transformations
  are bijective.}
\par
Finally, consider a coherent lower prevision $\lpr$ on $\gambles(\values)$,
and define the following set of gambles
\begin{equation*}
  \desirs_\lpr^+:=\set{f\in\gambles(\values)}
  {\text{$\lpr(f)>0$ or $f>0$}}.
\end{equation*}
Then $\desirs_\lpr^+\cup\{0\}$ is a coherent set of really desirable gambles,
i.e., it satisfies (D1)--(D4).\footnote{For (D1), use that a coherent lower
  prevision $\lpr$ satisfies $\lpr(f)\leq \sup f$ for any gamble $f$; for
  (D2), that $f\geq 0$ satisfies either $f>0$ or $f=0$; for (D3), use (P2) and
  the monotonicity of the coherent $\lpr$, and for (D4) use (P3).} Moreover,
any coherent set of really desirable gambles $\rdesirs$ that satisfies
\begin{equation*}
  \desirs_\lpr^+\cup\{0\}\subseteq\rdesirs\subseteq\desirs_\lpr,
\end{equation*}
i.e., the union of whose (relative) topological interior with $\cone^+$ is
$\desirs_\lpr^+\cup\{0\}$ and whose topological closure is $\desirs_\lpr$, has
$\lpr$ as its associated lower prevision, through~\eqref{eq:lpr}.  This
confirms what we claimed before: coherent lower previsions, or equivalently,
coherent sets of almost-desirable gambles, are less powerful belief models
than coherent sets of really desirable gambles. If a subject specifies a
coherent lower prevision $\lpr$, then she actually states that all gambles in
the union $\desirs_\lpr^+\cup\{0\}$ of $\cone_+$ with the relative topological
interior of $\desirs_\lpr$ are really desirable, but she doesn't specify
whether the gambles in the topological boundary
$\desirs_\lpr\setminus\desirs_\lpr^+$ of $\desirs_\lpr$ are: we only know that
they are almost-desirable to her.

\subsection{Natural extension for coherent lower previsions}
\label{sec:natex-lpr}
There is one important problem that we skipped over in the discussion above,
namely that of inference. Suppose a subject specifies a set $\desirs$ of
gambles that are almost-desirable to her. In an elicitation procedure, for
instance, this would typically be a finite set of gambles, so we cannot expect
this set to be coherent. We are then, as before for really desirable gambles,
faced with the problem of enlarging this $\desirs$ into a coherent set of
almost-desirable gambles that is as small as possible: we want to find out
what are the (behavioural) consequences of the subject's almost-accepting the
gambles in $\desirs$, taking into account \emph{only} the requirements of
coherence.
\par
The smallest closed convex cone including $\cone_+$ and $\desirs$, or in other
words, the smallest subset of $\gambles(\values)$ that includes $\desirs$ and
satisfies (M2)--(M5), is given by
\begin{equation}\label{eq:natex-desirs}
  \edesirs^m_\desirs:=
  \set{g\in\gambles(\values)}
  {(\forall\delta>0)(\exists n\geq0,\lambda_k\in\reals^+,f_k\in\desirs)
    g\geq\sum_{k=1}^n\lambda_k f_k-\delta}.\tag{M-NE}
\end{equation}
This is the topological closure of the set $\edesirs^r_\desirs$.  If this
convex cone $\edesirs^m_\desirs$ intersects
$\cone_-'=\set{f\in\gambles(\values)}{\sup f<0}$ then it is easy to see that
actually $\edesirs^m_\desirs=\gambles(\values)$, and then it is impossible to
extend $\desirs$ to a coherent set of almost-desirable gambles [because (M1)
cannot be satisfied].  Observe that $\edesirs^m_\desirs\cap\cone_-'=\emptyset$
if and only if\footnote{Actually, this condition is equivalent to the one
  where we always choose $\lambda_k=1$.}
\begin{equation}\label{eq:avoids-sure-loss-desirs}
  \sup\left[\sum_{k=1}^n\lambda_kf_k\right]\geq0
  \text{ for some $n\geq0$, $\lambda_k\in\reals^+$ and $f_k\in\desirs$},
  \tag{M-ASL}
\end{equation}
and we then say that the set $\desirs$ of almost-desirable gambles
\emph{avoids sure loss}. In that case, and only then, we are able to extend
$\desirs$ to a coherent set of almost-desirable gambles, and the smallest such
set is precisely $\edesirs^m_\desirs$, which is called the \emph{natural
  extension} of $\desirs$ to a set of almost-desirable gambles.
\par
What does natural extension mean for the equivalent model of coherent lower
previsions? Suppose our subject specifies a supremum acceptable buying price,
or lower prevision, $\lpr(f)$ for each gamble $f$ in some set of gambles
$\domain\subseteq\gambles(\values)$.\footnote{This set of gambles $\domain$
  need not have any predefined structure; in particular, it does not have to
  be a linear space.} We can then interpret $\lpr$ as a real-valued map on
$\domain$, and we call $\lpr$ a \emph{lower prevision on $\domain$}, and say
that $\domain$ is the \emph{domain} of $\lpr$.
\par
To study the problem of natural extension for this lower prevision, we shall
use what we already know about natural extension in the context of
almost-desirable gambles. Recall that specifying $\lpr$ on $\domain$ is
tantamount to stating that the gambles in the set
$\desirs:=\set{f-\lpr(f)}{f\in\domain}$ are almost-desirable. We now look at
the natural extension of this $\desirs$.
Using~\eqref{eq:avoids-sure-loss-desirs}, we know that such a natural
extension exists if and only if\footnote{Here too, this condition is
  equivalent to the one where we always choose $\lambda_k=1$.}
\begin{equation}\label{eq:avoids-sure-loss-lpr}
  \sup\left[\sum_{k=1}^n\lambda_k\left[f_k-\lpr(f_k)\right]\right]\geq0
  \text{ for all $n\geq0$, $\lambda_k\in\reals^+$ and $f_k\in\domain$},
  \tag{LPR-ASL}
\end{equation}
and we then say that the lower prevision $\lpr$ on $\domain$ \emph{avoids sure
  loss}. In this case, the natural extension $\edesirs^m_\desirs$ is the
smallest coherent set of almost-desirable gambles that includes $\desirs$, and
consequently the coherent lower prevision $\lpr_{\edesirs^m_\desirs}$
associated with $\edesirs^m_\desirs$ through
\begin{equation*}
  \lpr_{\edesirs^m_\desirs}(g)
  :=\max\set{s}{g-s\in\edesirs^m_\desirs}
\end{equation*}
is the point-wise smallest coherent lower prevision on $\gambles(\values)$
that dominates $\lpr$ on $\domain$. We call this coherent lower prevision the
\emph{natural extension} of $\lpr$ and we denote it by $\lnex_\lpr$. We deduce
from~\eqref{eq:natex-desirs} that for all gambles $g$ on $\values$:
\begin{equation}\label{eq:natex-lpr}
  \lnex_\lpr(g)
  =\sup_{\substack{\lambda_k\geq0,g_k\in\desirs\\k=1\dots,n,n\geq0}}
  \inf\left[g-\sum_{k=1}^n\lambda_kg_k\right]
  =\sup_{\substack{\lambda_k\geq0,f_k\in\domain\\k=1\dots,n,n\geq0}}
  \inf\left[g-\sum_{k=1}^n\lambda_k\left[f_k-\lpr(f_k)\right]\right].
  \tag{LPR-NE}
\end{equation}
If $\lpr$ \emph{incurs sure loss}, i.e., \eqref{eq:avoids-sure-loss-lpr} is
not satisfied, then $\edesirs^m_\desirs=\gambles(\values)$ and consequently
$\lnex_\lpr$ assumes the value $+\infty$ in every gamble.
\par
We shall call the lower prevision $\lpr$ on $\domain$ \emph{coherent},
whenever it can be extended to a coherent lower prevision on
$\gambles(\values)$, or in other words, whenever it coincides with its natural
extension $\lnex_\lpr$ on every gamble in its domain $\domain$. Taking into
account~\eqref{eq:natex-lpr}, we see that this happens exactly when
\begin{equation}\label{eq:coherence-lpr}
  \sup\left[\sum_{k=1}^n\lambda_k\left[f_k-\lpr(f_k)\right]
    -\lambda_0\left[f_0-\lpr(f_0)\right]\right]\geq0
  \text{ for all $n\geq0$, $\lambda_k\in\reals^+$ and $f_k\in\domain$},
  \tag{LPR-COH}
\end{equation}
This coherence condition implies that $\lpr$ avoids sure loss.
\par
Let us see if, for lower previsions, we can give a more immediate behavioural
interpretation for avoiding sure loss, coherence, and natural extension. This
should allow us to develop more intuition, as the approach we have followed so
far, which motivates these notions through the coherence axioms for real and
almost-desirable gambles, is admittedly quite abstract. We begin with avoiding
sure loss. Suppose that condition~\eqref{eq:avoids-sure-loss-lpr} is not
satisfied. Then there are $n\geq0$, $\lambda_1$, \dots, $\lambda_n$ in
$\reals^+$ and $f_1$, \dots $f_n$ in $\domain$ such that
$\sup\left[\sum_{k=1}^n\lambda_k\left[f_k-\lpr(f_k)\right]\right]<0$, which
implies that there is some $\delta>0$ for which
\begin{equation*}
  \sum_{k=1}^n\lambda_k\left[f_k-\lpr(f_k)+\delta\right]\leq-\delta.
\end{equation*}
Now, by the definition of $\lpr(f_k)$, our subject accepts each of the gambles
$f_k-\lpr(f_k)+\delta$, so she should also accept the combined gamble
$\sum_{k=1}^n\lambda_k[f_k-\lpr(f_k)+\delta]$ [use axioms~(D3) and~(D4) for
real desirability]. But this gamble leads to a sure loss of at least $\delta$.
In other words, if condition~\eqref{eq:avoids-sure-loss-lpr} doesn't hold,
there are gambles which the subject accepts and which, if properly combined,
make her subject to a sure loss.
\par
Next, assume that condition~\eqref{eq:coherence-lpr} fails to hold. Then there
are $n\geq0$, $\lambda_0$, \dots, $\lambda_n$ in $\reals^+$ and $f_0$, \dots
$f_n$ in $\domain$ such that
$\sup[\sum_{k=1}^n\lambda_k[f_k-\lpr(f_k)]-\lambda_0[f_0-\lpr(f_0)]]<0$.
Assume that $\lambda_0>0$, as we have already considered the case
$\lambda_0=0$ in our discussion of avoiding sure loss. Then there is some
$\delta>0$ such that
\begin{equation*}
  \sum_{k=1}^n\frac{\lambda_k}{\lambda_0}\left[f_k-\lpr(f_k)+\delta\right]
  \leq f_0-(\lpr(f_0)+\delta).
\end{equation*}
As before, the gamble on the left-hand side is a gamble that our subject
accepts. But then she should also accept the gamble $f_0-(\lpr(f_0)+\delta)$
since it point-wise dominates a gamble she accepts [use (D2) and (D3)]. This
implies that she should be willing to pay a price $\lpr(f_0)+\delta$ for
$f_0$, which is strictly higher than the supremum price $\lpr(f_0)$ she has
specified for it. Coherence avoids this kind of inconsistency.
\par
Finally, we turn to natural extension. Consider a gamble $g$ on $\values$,
then~\eqref{eq:natex-lpr} tells us that $\lnex_\lpr(g)$ is the supremum $s$
such that there are $n\geq0$, $\lambda_1$, \dots, $\lambda_n$ in $\reals^+$
and $f_1$, \dots $f_n$ in $\domain$ for which
\begin{equation*}
  g-s\geq\sum_{k=1}^n\lambda_k\left[f_k-\lpr(f_k)\right]
\end{equation*}
Now the expression on the right-hand side is almost-desirable, because it is a
non-negative linear combination of almost-desirable gambles [apply the
axioms~(M3) and~(M4)]. So $g-s$ should be almost-desirable as well [apply the
axioms~(M2) and~(M3)], and therefore our subject should be willing to buy $g$
for any price $t<s$. So we deduce that $\lnex_\lpr(g)$ is the supremum price
for $g$ that the subject can be forced to pay for the gamble $g$, by suitably
combining transactions that she is committed to accept by her specifying the
lower prevision $\lpr$ on $\domain$. In other words, $\lnex_\lpr(g)$ is the
lower prevision for $g$ that is implied by the assessments in $\lpr$ and
coherence \emph{alone}.

\subsection{Coherent previsions: the Bayesian belief models}
\label{sec:linear-previsions}
When a lower prevision $\lpr$ on $\domain$ is \emph{self-conjugate}, that is,
when $\lpr(f)=\upr(f)$ for any gamble $f$ in $\domain$, it is called a
\emph{prevision}. The common value $\pr(f)$ is then called the
\emph{prevision} of $f$; it is a \emph{fair price} for the gamble $f$ in the
sense of \citet{finetti19745}. Formally, a real-valued function $\pr$ on a
class of gambles $\domain$ is called a \emph{linear}, or \emph{coherent,
  prevision} whenever
\begin{equation}\label{eq:coherence-pr}
  \sup\left[\sum_{k=1}^n \left[f_k-\pr(f_k)\right]
    - \sum_{j=1}^{m} [g_j-\pr(g_j)]\right]\geq0
  \text{ for all $n,m\geq0$ and $f_k,g_j\in\domain$},
  \tag{PR-COH}
\end{equation}
A linear prevision is coherent, both as a lower and as an upper prevision.
Moreover, if its domain is the class of all gambles, $\gambles(\values)$, then
condition~\eqref{eq:coherence-pr} simplifies to
\begin{enumerate}[(PR1)]
\item $\pr(f+g)=\pr(f)+\pr(g)$ for any $f$ and $g$ in $\gambles(\values)$
  [linearity].
\item $\pr(f)\geq\inf f$ for any $f$ in $\gambles(\values)$ [accepting sure
  gains].
\end{enumerate}
Linear previsions are the familiar Bayesian belief models: any linear
prevision on all gambles is indeed a coherent prevision in the sense of
\citet{finetti19745}; and a prevision defined on an arbitrary set of gambles
is coherent exactly when it is the restriction of some coherent prevision on
all gambles. The restriction to (indicators of) events of a coherent prevision
on all gambles is a finitely additive probability. We shall denote by
$\linprevs(\values)$ the set of all coherent previsions on
$\gambles(\values)$.
\par
There is an interesting relationship between coherent previsions and coherent
lower previsions. Let $\lpr$ be a lower prevision with domain $\domain$, and
let us denote by
\begin{equation*}
  \solp(\lpr)
  :=\set{\pr\in\linprevs(\values)}
  {(\forall f \in \domain)\pr(f)\geq\lpr(f)}
\end{equation*}
the set of all coherent previsions on $\gambles(\values)$ that \emph{dominate}
$\lpr$ on its domain. Then it can be checked\footnote{See
  \cite[Sections~3.3--3.4]{walley1991} for proofs for these statements.} that
$\lpr$ avoids sure loss if and only if $\solp(\lpr)$ is non-empty, that is, if
and only if there is some coherent prevision on $\gambles(\values)$ that
dominates $\lpr$ on $\domain$, and $\lpr$ is coherent if and only if it is the
\emph{lower envelope} of $\solp(\lpr)$, meaning that for all $\lpr$ in
$\domain$,
\begin{equation*}
  \lpr(f)=\min\set{\pr(f)}{\pr\in\solp(\lpr)}.
\end{equation*}
Also, any lower envelope of a set of coherent previsions is a coherent lower
prevision. Moreover, the natural extension $\lnex_\lpr$ of $\lpr$ to all
gambles can be calculated using the set $\solp(\lpr)$ of coherent previsions:
for any gamble $f$ on $\values$, we have
\begin{equation*}
  \lnex_\lpr(f)=\min\set{\pr(f)}{\pr\in\solp(\lpr)}.
\end{equation*}
This means that from a \emph{mathematical} point of view, a coherent lower
prevision $\lpr$ and its set of dominating coherent lower previsions
$\solp(\lpr)$, are equivalent belief models.  It can be checked that this set
is convex and closed in the weak* topology.\footnote{The weak* topology on the
  set of all continuous linear functionals on $\gambles(\values)$ is the
  topology of point-wise convergence. For more details, see
  \citet[Appendix~D]{walley1991}.}  Moreover, there is a bijective
relationship between weak*-closed convex sets of coherent previsions and
coherent lower previsions (their lower envelopes).  This fact can (but need
not) be used to give coherent lower previsions a \emph{Bayesian sensitivity
  analysis interpretation}, besides the direct behavioural interpretation
given in Section~\ref{sec:almost-desirable}: we might assume the existence of
a precise but unknown coherent prevision $\pr$ expressing a subject's
behavioural dispositions, and we might model the information about $\pr$ by
means of a weak*-closed convex set of coherent previsions $\solp$ (the set of
possible candidates).  Then, this set is \emph{mathematically} equivalent to
its lower envelope $\lpr$, which is a coherent lower prevision. We shall come
back to the difference between the direct behavioural and the Bayesian
sensitivity analysis interpretation of a lower prevision in
Section~\ref{sec:strong-invariance}, when we discuss the interplay between
these interpretations and the notion of symmetry.
\par
Taking into account the bijective relationship that exists between coherent
lower previsions and sets of almost-desirable gambles, we may also establish a
bijective relationship between sets of coherent previsions and sets of
almost-desirable gambles: given a weak*-closed convex set $\solp$ of coherent
previsions on $\gambles(\values)$, the class
\begin{equation*}
  \desirs_\solp
  :=\set{f\in\gambles(\values)}{(\forall\pr\in\solp)\pr(f)\geq0}
\end{equation*}
is a coherent set of almost-desirable gambles, that is, it satisfies the
coherence conditions~(M1)--(M5). Conversely, given a coherent set of
almost-desirable gambles $\desirs$, the corresponding set of coherent
previsions
\begin{equation*}
  \solp(\desirs)
  :=\set{\pr\in\linprevs(\gambles)}{(\forall f\in\desirs)\pr(f)\geq0}
\end{equation*}
is a weak*-closed convex set of coherent previsions.
\par
Hence, there are at least three mathematically equivalent representations for
the behavioural dispositions of our subject: coherent sets of almost-desirable
gambles, coherent lower previsions, and weak*-closed convex sets of coherent
previsions. The bijective relationships between them are summarised in
Table~\ref{tab:equivalent-models}.

\begin{table}[htb]
  \centering
  \begin{tabular}{cccc}
    $\swarrow$ & $\desirs$ & $\lpr(\cdot)$ & $\solp$\\\\
    $\desirs$ & & $\set{f}{\lpr(f)\geq0}$
    & $\set{f}{(\forall\pr\in\solp)\pr(f)\geq0}$ \\\\
    $\lpr(\cdot)$ & $\max\set{s}{\cdot-s\in\desirs}$ &
    & $\min\set{\pr(\cdot)}{\pr\in\solp}$ \\\\
    $\solp$ & $\set{\pr}{(\forall f\in\desirs)\pr(f)\geq0}$
    & $\set{\pr}{(\forall f)\pr(f)\geq\lpr(f)}$ &\\\\
  \end{tabular}
  \caption{Bijective relationships between the equivalent belief models:
    coherent sets of almost-desirable gambles $\desirs$, coherent lower
    previsions $\lpr$ on $\gambles(\values)$, and weak*-closed convex sets
    $\solp$ of coherent previsions on $\gambles(\values)$}
  \label{tab:equivalent-models}
\end{table}

We now briefly discuss a number of belief models that constitute particular
instances of coherent lower previsions.  First, we consider $n$-monotone lower
previsions, where $n\geq1$. A lower prevision $\lpr$ is called
\emph{$n$-monotone}\footnote{See \Citet{cooman2005e,cooman2005d,cooman2005b}
  for a detailed discussion of $n$- and complete monotonicity for lower
  previsions.} when the following inequality holds for all $p\in\mathbb{N}$,
$p\leq n$, and all $f$, $f_1$, \dots, $f_p$ in $\gambles(\values)$:
\begin{equation*}
  \sum_{I\subseteq\{1,\dots,p\}}(-1)^{\abs{I}}
  \lpr\left(f\wedge\bigwedge_{i\in I}f_i\right)\geq0,
\end{equation*}
where, here and further on, $\abs{I}$ denotes the number of elements in a
finite set $I$. A similar definition can be given if the domain of $\lpr$ is
only a \emph{lattice of gambles}, i.e., a set of gambles closed under
point-wise minimum $\wedge$ and point-wise maximum $\vee$. Such $n$-monotone
lower previsions are particular instances of exact functionals
\citep{maass2003}, i.e., they are scalar multiples of some coherent lower
prevision. In particular, an $n$-monotone lower probability defined on a
lattice of events $\lattice$ that contains $\emptyset$ and $\values$ is
coherent if and only if $\lpr(\emptyset)=0$ and $\lpr(\values)=1$.
\par
A \emph{completely monotone} lower prevision is simply one that is
$n$-monotone for any natural number $n\geq1$. When it is defined on indicators
of events, it is called a completely monotone lower probability. When
$\values$ is finite, this leads to \emph{belief functions} in the terminology
of \citet{shafer1976}.
\par
Two particular cases of belief functions and their conjugate upper
probabilities are \emph{probability charges}, or finitely additive
probabilities defined on a field of events \citep{bhaskara1983} and
\emph{possibility measures}.  The latter \Citep{cooman2001,zadeh1978a} are set
functions $\Pi$ satisfying $\Pi\left(\bigcup_{i\in I}A_i\right)=\sup_{i\in
  I}\Pi(A_i)$ for any family $(A_i)_{i\in I}$ of subsets of $\values$. $\Pi$
is a coherent \emph{upper} probability if and only if $\Pi(\values)=1$.
\par
Finally, we can consider a particular instance of a completely monotone
coherent lower prevision that allows us to model complete ignorance, the
so-called \emph{vacuous lower prevision}. It is given by
\begin{equation*}
  \lpr_\values(f)=\inf_{x\in\values}f(x),
\end{equation*}
for all gambles $f$ on $\values$. It corresponds to the set of
almost-desirable gambles $\desirs=\cone_+=\set{f}{f\geq0}$, and to the set
$\solp=\linprevs(\gambles)$ of all coherent previsions on $\gambles$. If we
have no information at all about the values that $\rv$ takes in $\values$, we
have no reason to reject any coherent prevision $\pr$, and this leads to the
vacuous lower prevision as a belief model. More generally, we can consider a
vacuous lower prevision relative to some subset $A$ of $\values$, which is
given by
\begin{equation*}
  \lpr_A(f)=\inf_{x \in A}f(x).
\end{equation*}
A vacuous lower prevision relative to a set $A$ is the adequate belief model
when we know that the random variable $\rv$ assumes values in $A$, and nothing
else.  The restriction to events of a vacuous upper prevision is a
(zero-one-valued) possibility measure.

\subsection{Incomparability and indifference}
\label{sec:indifference}
We claimed in the Introduction that Bayesian belief models do not take
indecision seriously, and that we therefore need to look at a larger class of
belief models that do not have this defect. Here, we present a better
motivation for this claim.
\par
Consider two gambles $f$ and $g$ on $\values$. We say that a subject
\emph{almost-prefers} $f$ to $g$, and denote this as $f\pref g$, whenever she
accepts to exchange $g$ for $f$ in return for any (strictly) positive amount
of utility. Given this definition, it is straightforward to check that we can
express this in terms of the three equivalent belief models $\desirs$, $\lpr$
and $\solp$ of the previous sections by
\begin{align*}
  f\pref g
  &\asa f-g\in\desirs\\
  &\asa\lpr(f-g)\geq0\\
  &\asa(\forall\pr\in\solp)\pr(f)\geq\pr(g).
\end{align*}
The binary relation $\pref$ is a partial pre-order on $\gambles(\values)$,
i.e., it is reflexive and transitive.\footnote{The binary relation $\pref$ is
  actually a \emph{vector ordering} on the linear space $\gambles(\values)$,
  because it is compatible with the addition of gambles, and the scalar
  multiplication of gambles with non-negative real numbers.}  Observe also
that $f\pref g\asa f-g\pref0$ and that $f\pref0\asa f\in\desirs$, so $f$ is
almost-preferred to $g$ if and only if $f-g$ is almost-preferred to the zero
gamble, which in turn is equivalent to the fact that our subject
\emph{almost-accepts} $f-g$, i.e., that $f-g$ is almost-desirable to her.
\par
Unless our subject's lower prevision $\lpr$ is actually a (precise) prevision
$\pr$ (meaning that $\desirs$ is the semi-space $\set{f}{\pr(f)\geq0}$, and
that $\solp=\{\pr\}$), this ordering is not linear, or total: it does not hold
for all gambles $f$ and $g$ that $f\pref g$ or $g\pref f$. When, therefore,
both $f\not\pref g$ and $g\not\pref f$, we say that both gambles are
\emph{incomparable}, or that the subject is undecided about choosing between
$f$ and $g$, and we write this as $f\incomp g$.
\par
It is instructive to see why the relation $\incomp$ is non-empty unless $\lpr$
is a precise prevision $\pr$. If $\lpr$ is not precise (but coherent), there
is some gamble $h$ such $\lpr(h)<\upr(h)$. Let $x$ be any real number such
that $\lpr(h)<x<\upr(h)$. In this case, the subject does not express a
willingness to buy $h$ for the price $x$, because $x$ is strictly greater than
her supremum acceptable price $\lpr(h)$ for buying $h$. Nor does she express a
willingness to sell $h$ for a price $x$, because $x$ is strictly smaller than
her infimum acceptable price $\upr(h)$ for selling $h$. But there is more.
Consider the gambles $f:=h-x$ (buying $h$ for a price $x$) and $g:=x-h$
(selling $h$ for a price $x$). Then it follows from the coherence of $\lpr$
that
\begin{equation*}
  \lpr(f-g)=2\lpr(h-x)=2[\lpr(h)-x]<0
  \text{ and }
  \lpr(g-f)=2\lpr(x-h)=2[x-\upr(h)]<0,
\end{equation*}
so $f\incomp g$: our subject is also undecided in the choice between buying
$h$ for $x$ or selling $h$ for that price.
\par
We say that our subject is \emph{indifferent} between $f$ and $g$, and denote
this as $f\indif g$ whenever both $f\pref g$ and $g\pref f$. This means that
$\lpr(f-g)=\lpr(g-f)=0$, or equivalently, $\pr(f)=\pr(g)$ for all $\pr$ in
$\solp$. Clearly, $\indif$ is an equivalence relation (a reflexive,
symmetrical and transitive binary relation) on $\gambles(\values)$. It is
important to distinguish between incomparability and indifference.
Indifference between gambles $f$ and $g$ represents strong behavioural
dispositions: it means that our subject almost-accepts to exchange $f$ for $g$
and \textit{vice versa}; on the other hand, incomparability has no behavioural
implications, it merely records the absence of a(n expressed) behavioural
disposition to choose between $f$ and $g$.

\section{Monoids of transformations}
\label{sec:transformations}
Symmetry is generally characterised mathematically as invariance under certain transformations. In this section, we provide the necessary mathematical apparatus that will allow us to describe and characterise symmetry for the belief models we are interested in.

\subsection{Transformations and lifting}
We are interested in models for beliefs that concern a random variable $\rv$.  So let us begin by concentrating on transformations of the set of possible values $\values$ for $\rv$.  A \emph{transformation} of $\values$ is defined mathematically as a map $T:\values\to\values\colon x\mapsto Tx$. At this point, we do not require that such a map $T$ should be \emph{onto} (or surjective), i.e., that $T(\values):=\set{Tx}{x\in\values}$ should be equal to $\values$. Neither do we require that $T$ should be \emph{one-to-one} (or injective), meaning that $Tx=Ty$ implies $x=y$ for all $x$ and $y$ in $\values$. A transformation of $\values$ that is both onto and one-to-one will be called a \emph{permutation} of $\values$, but we shall in the sequel also need to consider transformations of $\values$ that are not permutations.
\par
Suppose we have two transformations, $T$ and $S$, of $\values$ that are of interest. Then there is no real reason why we shouldn't also consider the combined action of $T$ and $S$ on $\values$, leading to new transformations $ST:=S\circ T$ and $TS:=T\circ S$, defined by $(ST)x:=S(Tx)$ and similarly $TSx:=T(Sx)$ for all $x$ in $\values$.  And of course, we could also consider in a similar way $TST$ and $STS$, or for that matter $TTTSST$, which we shall also write as $T^3S^2T$. So it is natural in this context to consider a set $\transfos$ of transformations of $\values$ that is closed under composition, i.e.,
\begin{equation}\label{eq:semigroup}
  (\forall T,S\in\transfos)(TS\in\transfos)\tag{SG}
\end{equation}
Such a set is called a \emph{semigroup of transformations}.\footnote{A semigroup is defined as a set with a   binary operation that is internal and associative. Composition of maps is always an associative binary   operation, and~\eqref{eq:semigroup} guarantees that it is internal in $\transfos$.}  If moreover the semigroup $\transfos$ contains the identity map $\iden_\values$, defined by $\iden_\values x:=x$ for all $x$ in $\values$, it is called a \emph{monoid}. As the identity map leaves all elements of $\values$ unchanged, it has no implications as far as symmetry and invariance are concerned, and we can therefore in what follows assume without loss of generality that any $\transfos$ we consider actually contains $\iden_\values$ (is a monoid).
\par
A monoid $\transfos$ is \emph{Abelian} if $ST=TS$ for all $T$ and $S$ in $\transfos$. An important example of an Abelian monoid is the following.  Consider a single transformation $T$ of $\values$, and the Abelian monoid $\transfos_T$ generated by $T$, given by
\begin{equation*}
  \transfos_T:=\set{T^n}{n\geq0},
\end{equation*}
where $T^0:=\iden_\values$ is the identity map on $\values$, $T^1:=T$ and for $n\geq2$,
\begin{equation*}
  T^n:=\underset{\text{$n$ times}}{\underbrace{T\circ T\circ\dots\circ T}}.
\end{equation*}
\par
A monoid $\transfos$ of transformations is called \emph{left-} (respectively~\emph{right-})\emph{cancellable} when for every transformation $T$ in $\transfos$ there is some $S$ in $\transfos$ such that $ST=\iden_\values$ (respectively~$TS=\iden_\values$). This transformation $S$ is then called a \emph{left-} (respectively~\emph{right-})\emph{inverse} of $T$. If $\transfos$ is both left- and right-cancellable, then the left-and right-inverses of $T$ are unique and coincide for any $T$ in $\transfos$, and $\transfos$ is called a \emph{group}. Any element of $\transfos$ is then a permutation of $\values$.
\par
For our purposes here, we generally only need to assume that $\transfos$ is a monoid, because there interesting (and relevant) situations where $\transfos$ is not a group; this is for instance the case for the Abelian monoid of the \emph{shift} transformations of the set of natural numbers $\nats$:
\begin{equation}\label{eq:shift}
  \transfos_\theta:=\set{\theta^n}{n\geq0},
\end{equation}
where $\theta(m)=m+1$, and $\theta^n(m)=m+n$ for all natural numbers $m$ and $n$. Another important example is the monoid $\transfos_\values$ of all transformations of $\values$, which is generally not Abelian, nor a group.
\par
Since we are also concerned with gambles $f$ on $\values$, we need a way to turn a transformation of $\values$ into a transformation of $\gambles(\values)$. This is done by the procedure of \emph{lifting}: given any gamble $f$ on $\values$, we shall denote by $T^tf$ the gamble $f\circ T$, i.e.,
\begin{equation*}
  T^tf(x):=f(Tx),
\end{equation*}
for all $x$ in $\values$. For an event $A$, $T^tI_A=I_{T^{-1}(A)}$, where $T^{-1}(A):=\set{x\in\values}{Tx\in   A}$ is the so-called \emph{inverse image} of $A$ under $T$. On the other hand, given a constant $\mu$, we have $T^t \mu=\mu$ for any transformation $T$.
\par
The following observation is quite important. Consider two transformations $T$ and $S$ on $\values$. Then for any gamble $f$ on $\values$ we see that
\begin{equation*}
  (ST)^tf=f\circ(S\circ T)=(f\circ S)\circ T=(S^tf)\circ T=T^t(S^tf),
\end{equation*}
so $(ST)^t=T^tS^t$, and lifting reverses the order of application of the transformations: for $x$ in $\values$, $STx$ means that $T$ is applied first to $x$, and then $S$ to $Tx$. For $f$ in $\gambles(\values)$, $(ST)^tf$ means that $S^t$ is applied first to $f$ and then $T^t$ to $S^tf$.
\par
Any transformation $T$ of $\values$ can therefore be lifted to a transformation $T^t$ of $\gambles(\values)$, and we denote the corresponding set of liftings by $\transfos^t$. $\transfos^t$ is then a monoid of transformations of $\gambles(\values)$. Lifting preserves the most common properties of semigroups, taking into account the above-mentioned order-inversion: being a monoid, being Abelian, and being a group are preserved under lifting. But being left-cancellable is turned into being right-cancellable, and \textit{vice   versa}. Lifting also has the interesting property that it turns a transformation $T$ on $\values$ into a \emph{linear} transformation $T^t$ of the linear space $\gambles(\values)$: for any pair of gambles $f$ and $g$ on $\values$ and any real numbers $\lambda$ and $\mu$, we have
\begin{equation*}
  T^t(\lambda f+\mu g)=\lambda T^tf+\mu T^tg.
\end{equation*}

\subsection{Invariant (sets of) gambles}
We now turn to the important notions of invariance under transformations. We start with the invariance of a set of gambles, because that is the most general notion, from which all other notions of invariance can be derived. If $\domain$ is a set of gambles on $\values$, and $T$ any transformation of $\transfos$, then we denote by
\begin{equation*}
  T^t\domain:=\set{T^tf}{f\in\domain}
\end{equation*}
the direct image of the set $\domain$ under $T^t$, and we say that \emph{$\domain$ is $\transfos$-invariant} if
\begin{equation*}
  (\forall T^t\in\transfos^t)(T^t\domain\subseteq\domain),
\end{equation*}
i.e., if all transformations in $\transfos^t$ are \emph{internal} in $\domain$.\footnote{So $\transfos^t$ is a   monoid of transformations of $\domain$.}
\par
A gamble $f$ on $\values$ is called \emph{$\transfos$-invariant} if the singleton $\{f\}$ is, i.e., if $T^tf=f$ for all transformations $T$ in the monoid $\transfos$.  We call an event $A$ $\transfos$-invariant if its indicator $I_A$ is, i.e., if $T^{-1}(A)=A$ for all $T$ in $\transfos$.
\par
Let us denote by $\invars_\transfos$ the set of all $\transfos$-invariant events. It is easy to check that $\invars_\transfos$ is an \emph{ample field}, i.e., it contains $\emptyset$ and $\values$, and it is closed under arbitrary unions and complementation, and therefore also under arbitrary intersections.  For any $x$ in $\values$, we shall call
\begin{equation*}
  \atom{x}{\transfos}
  :=\bigcap\set{A}{\text{$A\in\invars_\transfos$ and $x\in A$}}
\end{equation*}
the \emph{$\transfos$-invariant atom containing $x$}. It is the smallest $\transfos$-invariant event that contains $x$. Any $\transfos$-invariant event $A$ is a union of $\transfos$-invariant atoms: $A=\bigcup_{x\in   A}\atom{x}{\transfos}$. We shall denote by $\atoms_\transfos$ the set of all invariant atoms: $\atoms_\transfos:=\set{\atom{x}{\transfos}}{x\in\values}$.  It is a partition of $\values$. A gamble $f$ on $\values$ is $\transfos$-invariant if and only if it is constant on the $\transfos$-invariant atoms of $\values$.
\par
Of course, the bigger the set of transformations $\transfos$, the smaller the number of $\transfos$-invariant events (or, equivalently, the bigger the atoms $\atom{x}{\transfos}$).  The following proposition relates the $\transfos$-invariant atoms $\atom{x}{\transfos}$ to the images of $x$ under the transformations in $\transfos$.

\begin{proposition}\label{prop:invariant-atoms}
  Let $\transfos$ be a monoid of transformations of $\values$, and let $x$ be any element of $\values$. In   general we have that $\set{Tx}{T\in\transfos}\subseteq\atom{x}{\transfos}$.  If $\transfos$ is   left-cancellable, then $\atom{x}{\transfos}=\set{Tx}{T\in\transfos}$.
\end{proposition}

\begin{proof}
  Fix $x$ in $\values$.  Let $\transfos(x):=\set{Tx}{T\in\transfos}$ for brevity of notation.  Consider any   $T$ in $\transfos$. Since $\atom{x}{\transfos}$ is $T$-invariant, we have that   $T^{-1}(\atom{x}{\transfos})=\atom{x}{\transfos}$. Since $x\in\atom{x}{\transfos}$ because $\transfos$ is a   monoid, we infer from this equality that $Tx\in\atom{x}{\transfos}$. Hence indeed   $\transfos(x)\subseteq\atom{x}{\transfos}$.
  \par
  To prove the converse inequality, assume that $\transfos$ is left-cancellable. Consider any $S$ in   $\transfos$. If we can prove that $\transfos(x)$ is $S$-invariant, meaning that   $S^{-1}(\transfos(x))=\transfos(x)$, then the proof is complete, since then $\transfos(x)$ will be   $\transfos$-invariant, and since this set contains $x$ [because $\iden_\values\in\transfos$], it must   include the smallest $\transfos$-invariant set $\atom{x}{\transfos}$ that contains $x$.  So we set out to   prove that $S^{-1}(\transfos(x))=\transfos(x)$.  Consider any $y$ in $\values$. First assume that   $y\in\transfos(x)$. Then there is some $T$ in $\transfos$ such that $y=Tx$, whence $Sy=STx\in\transfos(x)$,   since $ST\in\transfos$.  Conversely, assume that $y\in S^{-1}(\transfos(x))$, or equivalently, that   $Sy\in\transfos(x)$, then there is some $T$ in $\transfos$ such that $Sy=Tx$, and since $\transfos$ is   assumed to be left-cancellable, there is some $S'$ in $\transfos$ such that $S'S=\iden_\values$, whence   $\transfos(x)\ni S'Tx=S'Sy=y$, since $S'T\in\transfos$.
\end{proof}

An important special case is the following. Consider a transformation $T$ of $\values$, and the Abelian monoid $\transfos_T=\set{T^n}{n\geq0}$ generated by $T$. Then a set of gambles $\domain$ is $\transfos_T$-invariant if and only if $T^t\domain\subseteq\domain$, and we simply say that $\domain$ is \emph{$T$-invariant}. Similarly, a gamble $f$ is $\transfos_T$-invariant if and only if $T^tf=f$, and we say that $f$ is \emph{$T$-invariant}. In what follows, we shall always use the phrase `$T$-invariant' for `$\transfos_T$-invariant'. Also $\invars_T$ is the set of $T$-invariant events, and it is an ample field whose atoms are denoted by $\atom{x}{T}$.  With this notation, we have for an arbitrary monoid $\transfos$ that $\invars_\transfos=\bigcap_{T\in\transfos}\invars_T$.
\par
For instance, the particular case of the shift transformations of $\nats$ given by Eq.~\eqref{eq:shift} concerns the Abelian monoid generated by $\theta$.  Here, the only $\theta$- (or shift-)invariant events are $\emptyset$ and $\nats$, and consequently a gamble $f$ on $\nats$ is $\theta$-invariant if and only if it is constant. This also shows that the equality in the first part of Proposition~\ref{prop:invariant-atoms} need not hold when the monoid of transformations $\transfos$ is not left-cancellable: in the present case, we have that $\transfos_\theta(m)=\set{\theta^n(m)}{n\geq0}=\set{n\in\nats}{n\geq m}$ is strictly included in the invariant atom $\atom{m}{\theta}=\nats$ for all $m\geq1$.
\par
Another interesting case is that of $\transfos_\values$, the class of all transformations of $\values$. This a monoid, but it is not generally a group, nor Abelian. Moreover, it is not generally left-cancellable.  We have, for any element $x$ of $\values$ that $\set{Tx}{T\in\transfos_\values}=\values$, and from Proposition~\ref{prop:invariant-atoms} we deduce in a trivial manner that $\atom{x}{\transfos_\values}=\values$: the only invariant events under all transformations of $\values$ are $\emptyset$ and $\values$.  This shows that the left-cancellability condition in the second part of Proposition~\ref{prop:invariant-atoms} is not generally necessary.

\section{Symmetry and invariance for belief models}
\label{sec:belief-model-invariance}
We now have the necessary mathematical tools for studying the issue of symmetry in relation to the belief models discussed in Section~\ref{sec:ip-models}. We shall see that for these coherent sets of almost-desirable gambles, there is an important distinction between the concepts `symmetry of models' (which we shall call weak invariance) and `models of symmetry' (which we shall call strong invariance). Let us first turn to the discussion of symmetrical belief models.

\subsection{Weak invariance: symmetry of models}
\label{sec:symmetry-of-models}
Consider a monoid $\transfos$ of transformations of $\values$.  We want to express that a belief model about the value that the random variable $\rv$ assumes in $\values$, exhibits a symmetry that is characterised by the transformations in $\transfos$. Thus, the notion of (weak) invariance of belief models that we are about to introduce is in a sense a purely mathematical one: it expresses that these belief models are left invariant under the transformations in $\transfos$.

\begin{definition}[Weak invariance]\label{def:weak-invariance}
  A coherent set of almost-desirable gambles $\desirs$ is called \emph{weakly} \emph{$\transfos$-invariant} if   it is $\transfos$-invariant as a set of gambles, i.e., if $T^t\desirs\subseteq\desirs$ for all $T$ in   $\transfos$.
\end{definition}
\noindent
Why don't we require equality rather than the weaker requirement of set inclusion in this definition? In linear algebra, invariance of a subset of a linear space with respect to a linear transformation of that space is generally defined using only the inclusion. If we recall from Section~\ref{sec:transformations} that lifting turns any transformation $T$ of $\values$ into a linear transformation $T^t$ of the linear space $\gambles(\values)$, we see that our definition of invariance is just a special case of a notion that is quite common in the mathematical literature.
\par
A few additional comments are in order. First of all, any coherent set of almost-desirable gambles is weakly $\iden_\values$-invariant, so we may indeed always assume without loss of generality that $\transfos$ is at least a monoid (contains $\iden_\values$).
\par
Secondly, we have given an invariance definition for almost-desirability, but the definition for coherent sets of really desirable gambles $\rdesirs$ is completely analogous: for all $T$ in $\transfos$, $T^t\rdesirs\subseteq\rdesirs$. Observe that if $\rdesirs$ is weakly $\transfos$-invariant then the associated set of almost-desirable gambles $\desirs_\rdesirs$, given by~\eqref{eq:ad}, is weakly $\transfos$-invariant as well.
\par
Thirdly, if $\transfos$ is a group (or at least left-cancellable), then the weak invariance condition is actually equivalent to $T^t\desirs=\desirs$ for all $T$ in $\transfos$: given a transformation $T$ in $\transfos$ and its (left-)inverse $S\in\transfos$, consider $f\in\desirs$; then $T^t(S^tf)=(ST)^tf=f$, so there is a gamble $g=S^tf$, which belongs to $\desirs$ by weak invariance, such that $f=T^tg$; this means that $f\in T^t\desirs$, so $\desirs\subseteq T^t\desirs$ as well.
\par
In summary, weak invariance is a mathematical notion that states that a subject's behavioural dispositions, as represented by a belief model $\desirs$, are invariant under certain transformations. If we posit that a subject's dispositions are in some way a reflection of the evidence available to her, we see that weak invariance is a way to model `symmetry of evidence'.  The following examples try to argue that if there is `symmetry of evidence', then corresponding belief models should at least be weakly invariant.

\begin{shifts}
  Suppose our subject is completely ignorant about the value of a random variable $\rv$ that assumes only   non-negative integer values, so $\values=\nats$. If her belief model is to be a reflection of the available   evidence (none), we should like it to be weakly invariant with respect to the shift transformations in   $\transfos_\theta$ ( which is an Abelian monoid, but not a group). Indeed, if she is ignorant about $\rv$,   she is also ignorant about $\theta(\rv)=\rv+1$, apart from the fact that she knows that $\theta(\rv)$ cannot   assume the value $0$, whereas $\rv$ can.  Therefore, if our subject almost-accepts a gamble $f$, she should   almost-accept $\theta^tf$: $\theta^tf(\rv)=f(\rv+1)$ may assume the same values as $f(\rv)$, apart from the   value $f(0)$, and because of her ignorance, our subject has no reason to treat the shifted gamble   differently.  $\blacklozenge$
\end{shifts}

\begin{dice}
  Let us go back to the die example. Suppose that whatever evidence our subject has about the outcome $\rv$ of   rolling the die, is left invariant by permutations $\pi$ of $\values_6=\{1,\dots,6\}$. Assume that our   subject almost-accepts a gamble $f$, meaning that she is willing to accept the uncertain reward   $f(\rv)+\epsilon$ for any $\epsilon>0$. But since the evidence gives our subject no reason to distinguish   between the random variables $\rv$ and $\pi(\rv)$, she should also be willing to accept the uncertain reward   $f(\pi(\rv))+\epsilon$ for any $\epsilon>0$, or in other words, she should almost-accept the gamble $\pi^t   f$.
\end{dice}

\par
We now investigate the corresponding notions for weak invariance for the equivalent belief models: coherent lower previsions and weak*-closed convex sets of coherent previsions. In order to do this, it is convenient to define the transformation of a (lower) prevision under a transformation $T$ on $\values$, by lifting $T$ to yet a higher level.

\begin{definition}[Transformation of a functional]\footnote{We use the same
    notation $T$ for the transformation of $\values$ and for the corresponding
    transformation of a functional, first of all because we do not want to
    overload the mathematical notation, and also because, in contrast with
    lifting only once, lifting twice preserves the order of application of
    transformations.} Let $T$ be a transformation of $\values$ and let
  $\Lambda$ be a real-valued functional defined on a $T$-invariant set of gambles   $\domain\subseteq\gambles(\values)$. Then the transformation $T\Lambda$ of $\Lambda$ is the real-valued   functional defined on $\domain$ by $T\Lambda:=\Lambda\circ T^t$, or equivalently, by   $T\Lambda(f):=\Lambda(T^tf)=\Lambda(f\circ T)$ for all gambles $f$ in $\domain$.
\end{definition}

\begin{theorem}\label{theo:weak-invariance}
  Let $\lpr$ be a coherent lower prevision on $\gambles(\values)$, $\desirs$ a coherent set of   almost-desirable gambles, and $\solp$ a weak*-closed convex set of coherent previsions on   $\gambles(\values)$. Assume that these belief models are equivalent, in the sense that they correspond to   one another using the bijective relations in Table~\ref{tab:equivalent-models}.  Then the following   statements are equivalent.
  \begin{enumerate}[1.]
  \item $\desirs$ is weakly $\transfos$-invariant, in the sense that $T^t\desirs\subseteq\desirs$ for all $T$     in $\transfos$.
  \item $\lpr$ is \emph{weakly $\transfos$-invariant}, in the sense that $T\lpr\geq\lpr$ for all $T$ in     $\transfos$, or equivalently $\lpr(T^tf)\geq\lpr(f)$ for all $T$ in $\transfos$ and $f$ in     $\gambles(\values)$;
  \item $\solp$ is \emph{weakly $\transfos$-invariant}, in the sense that $T\solp\subseteq\solp$ for all $T$     in $\transfos$, or equivalently, $T\pr\in\solp$ for all $\pr$ in $\solp$ and all $T$ in     $\transfos$.\footnote{This shows that our notion of a weakly invariant belief model corresponds to       \citegen{pericchi1991} notion of a `reasonable (or invariant) class of priors', rather than a `class of       reasonable (or invariant) priors', the latter being what our notion of strong invariance will correspond       to. On the other hand, \citet[Definition~3.5.1]{walley1991} defines a $\transfos$-invariant lower       prevision $\lpr$ as one for which $\lpr(T^tf)=\lpr(f)$ for all $T\in\transfos$ and all gambles $f$, so       he requires equality rather than inequality, as we do here.}
  \end{enumerate}
\end{theorem}

\begin{proof}
  We give a circular proof.  Assume that $\desirs$ is weakly $\transfos$-invariant.  Consider any $T$ in   $\transfos$ and $f$ in $\domain$, and observe that for the corresponding lower prevision $\lpr$
  \begin{equation*}
    \lpr(T^tf)
    =\max\set{\mu}{T^tf-\mu\in\desirs}
    \geq\max\set{\mu}{f-\mu\in\desirs}
    =\lpr(f),
  \end{equation*}
  where the inequality follows from the invariance assumption on $\desirs$.  This shows that the first   statement implies the second.
  \par
  Next, assume that $\lpr$ is weakly $\transfos$-invariant, and consider any $T$ in $\transfos$ and $\pr$ in   the corresponding $\solp=\solp(\lpr)=\set{\pr}{(\forall f)\pr(f)\geq\lpr(f)}$. Then for any gamble $f$ on   $\values$ we have that $T\pr(f)=\pr(T^tf)\geq\lpr(T^tf)\geq\lpr(f)$, where the second inequality follows for   the invariance assumption on $\lpr$. This tells us that indeed $T\pr\in\solp(\lpr)$, so the second statement   implies the third.
  \par
  Finally, assume that $\solp$ is weakly $\transfos$-invariant.  Consider any $T$ in $\transfos$ and any   gamble $f$ in the corresponding $\desirs=\desirs_\solp=\set{f}{(\forall\pr\in\solp)\pr(f)\geq0}$. Then we   have for any $\pr$ in $\solp$ that $\pr(T^tf)=T\pr(f)\geq0$, since $T\pr$ belongs to $\solp(\lpr)$ by the   invariance assumption on $\solp$.  Consequently $T^tf\in\desirs$, which proves that the third statement   implies the first.
\end{proof}

A coherent prevision $\pr$ on $\gambles(\values)$ is weakly $\transfos$-invariant if and only if $T\pr=\pr$ for all $T$ in $\transfos$.  This is easiest to prove by observing that $\solp(\pr)=\{\pr\}$.\footnote{See   Proposition~\ref{prop:weak-invariance-previsions} for a more direct proof.}  So for coherent previsions, we have an equality in the weak invariance condition. As we argued before, we generally won't have such an equality for arbitrary monoids $\transfos$, but the following corollary gives another sufficient condition on $\transfos$.

\begin{corollary}
  If the monoid $\transfos$ is left-cancellable, then the first weak invariance condition in   Theorem~\ref{theo:weak-invariance} becomes $T^t\desirs=\desirs$ for all $T$ in $\transfos$. If $\transfos$   is right-cancellable, then the second and third weak invariance conditions become $T\lpr=\lpr$ and   $T\solp=\solp$ for all $T$ in $\transfos$.\footnote{The reason for the difference in terms of left- versus     right-cancellability lies of course in the fact that in the first condition, we work with transformations     $T^t$ of gambles, and in the second and third condition we work with transformations $T$ of functionals,     which are liftings of the former; simply recall that lifting reverses the order of application of     transformations.}
\end{corollary}

\begin{proof}
  We have already proven the first statement near the beginning of Section~\ref{sec:symmetry-of-models}. To   prove the second statement, it suffices to show that when $\transfos$ is right-cancellable,   $\transfos$-invariance implies that $\lpr\geq T\lpr$ and $\solp\subseteq T\solp$ for all $T$ in   $\transfos$. Consider any transformation $T$ in the monoid $\transfos$, and let $R$ be a right-inverse for   $T$, i.e., $TR=\iden_\values$. Consider a gamble $h$ on $\values$, then   $\lpr(h)=\lpr((TR)^th)=\lpr(R^t(T^th))\geq\lpr(T^th)$, where the inequality follows from the weak invariance   of $\lpr$. So indeed, $\lpr\geq T\lpr$.  Similarly, consider $\pr$ in $\solp$. Then $R\pr\in\solp$ by weak   invariance, and for any gamble $f$ on $\values$, $T(R\pr)(f)=R\pr(T^tf)=\pr(R^t(T^tf))=\pr(f)$ since   $R^t(T^tf)=(TR)^tf=f$.  So there is a $\apr=R\pr$ in $\solp$ such that $\pr=T\apr$, meaning that $\pr\in   T\solp$. So indeed $\solp\subseteq T\solp$.
\end{proof}

We see from the definition that if a coherent set of almost-desirable gambles $\desirs$ (or a coherent lower prevision, or a weak*-closed convex set of coherent previsions) is weakly $\transfos$-invariant, it is also weakly $\transfos'$-invariant for any sub-monoid of transformations $\transfos'\subseteq\transfos$. Hence, as we add transformations, the collection of weakly invariant belief models will not increase. The limit case is when we consider the class $\transfos_\values$ of all transformations on $\values$.  The following theorem shows that the vacuous belief models are the only ones that are \emph{completely weakly invariant}, i.e., weakly $\transfos_\values$-invariant.

\begin{theorem}\label{theo:vacuous-invariant}
  Let $\transfos_\values$ be the monoid of all transformations of $\values$.  Then the vacuous coherent set of   almost-desirable gambles $\cone_+$ (or equivalently, the vacuous lower prevision $\lpr_\values$, or   equivalently, the weak*-closed convex set of all coherent previsions $\linprevs(\values)$) is the only   coherent set of almost-desirable gambles (coherent lower prevision, weak*-closed convex set of coherent   previsions) that is weakly $\transfos_\values$-invariant.
\end{theorem}

\begin{proof}
  We give the proof for coherent sets of almost-desirable gambles. It is obvious that $\cone_+$ is   $\transfos_\values$-invariant. So, consider any $\transfos_\values$-invariant coherent set of   almost-desirable gambles $\desirs$.  It follows from coherence [axiom~(M2)] that $\cone_+\subseteq\desirs$.   Assume \emph{ex absurdo} that $\cone_+\subset\desirs$ and let $f$ be any gamble in   $\desirs\setminus\cone_+$. This means that there is some $x_0$ in $\values$ such that $f(x_0)<0$. Consider   the transformation $T_{x_0}$ of $\values$ that maps all elements of $\values$ to $x_0$, then   $T_{x_0}^tf=f(x_0)$ and it follows from the $T_{x_0}$-invariance of $\desirs$ that the constant gamble   $f(x_0)\in\desirs$, which violates coherence axiom~(M1), so $\desirs$ cannot be coherent, a   contradiction.\footnote{A similar argument tells us that the same result holds for complete weak invariance     of coherent sets of really desirable gambles, where now the axiom~(D1) will be violated.}
\end{proof}

This result also tells us in particular that the vacuous belief model is always $\transfos$-invariant for any monoid of transformations $\transfos$.  This implies that for any monoid of transformations $\transfos$, there always are $\transfos$-invariant belief models.
\par
What are the behavioural consequences of weak invariance with respect to a monoid of transformations $\transfos$? It seems easiest to study this in terms of coherent lower previsions. First of all, we have that for any gamble $f$ on $\values$ and any $T$ in $\transfos$, our subject's supremum buying price $\lpr(T^tf)$ for the transformed gamble $T^tf$ should not be strictly smaller that her supremum price $\lpr(f)$ for buying $f$ itself.
\par
But there is also a more interesting consequence. Indeed, it follows from the coherence of $\lpr$ that
\begin{equation*}
  \lpr(f-T^tf)\leq\lpr(f)-\lpr(T^tf)\leq0.
\end{equation*}
\citet[Section~3.8.1]{walley1991} suggests that a subject \emph{strictly prefers} a gamble $f$ to a gamble $g$, which we denote as $f\spref g$, if $f>g$, or also if she accepts to pay some (strictly) positive price for exchanging $g$ with $f$, so if $\lpr(f-g)>0$. This means that weak $\transfos$-invariance implies that
\begin{equation*}
  \text{$f\not\spref T^tf$ for all $f$ in $\gambles(\values)$ and all $T$ in
    $\transfos$ such that $f\not>T^tf$}
\end{equation*}
which models that our subject \emph{has no reason} (or disposition) \emph{to strictly prefer any gamble $f$ to   any of its transformations $T^tf$ that it doesn't strictly dominate}.

\subsection{Strong invariance: models of symmetry}
\label{sec:strong-invariance}
Next, suppose that our subject believes that the (phenomenon underlying the) random variable $\rv$ is subject to symmetry with respect to the transformations $T$ in $\transfos$, so that she has \emph{reason not to   distinguish} between a gamble $f$ and its transformation $T^tf$. Let us give an example to get a more intuitive understanding of what this means.

\begin{dice-con}
  Again, let us go back to the die example. Consider the gambles $I_{\{x\}}$, for   $x\in\values_6:=\{1,\dots,6\}$. Since our subject believes the die (and the rolling mechanism behind it) to   be symmetrical, she will be willing to exchange any gamble $I_{\{x\}}$ for any other gamble $I_{\{y\}}$ in   return for any strictly positive amount of utility: $I_{\{x\}}-I_{\{y\}}$ should therefore be   almost-desirable to her, or in other words, in terms of her lower prevision $\lpr$:
  \begin{equation*}
    \lpr(I_{\{x\}}-I_{\{y\}})\geq0\text{ for all $x$ and $y$ in $\values_6$}.
  \end{equation*}
  This is equivalent to stating that $I_{\{x\}}-\pi^tI_{\{x\}}$ should be almost-desirable, or that   $\lpr(I_{\{x\}}-\pi^tI_{\{x\}})\geq0$ for all $x\in\values_6$ and all permutations $\pi$ of $\values_6$. Now   the only coherent lower prevision that satisfies these requirements is the uniform (precise) prevision,   which assigns precise probability $\frac{1}{6}$ to each event $\{x\}$ [simply observe that for any coherent   prevision $\pr$ in $\solp(\lpr)$ it follows from these requirements that   $\pr(I_{\{x\}})=\pr(I_{\{y\}})$]. $\blacklozenge$
\end{dice-con}

Let us now try and formalise the intuitive requirements in this example into a more formal definition. We stated above that if our subject believes that the (phenomenon underlying the) random variable $\rv$ is subject to symmetry with respect to the transformations $T$ in $\transfos$, then she has \emph{reason not to   distinguish} between a gamble $f$ and its transformation $T^tf$.  Suppose she has the gamble $f$ in her possession, then she should be willing to exchange this for the gamble $T^tf$ in return for any strictly positive price, and \textit{vice versa}. This means that she should almost-accept both $f-T^tf$ and $T^tf-f$, or in the language of Section~\ref{sec:indifference}, that she is \emph{indifferent between $f$ and $T^tf$}: $f\indif T^tf$.  If $\desirs$ is her coherent set of almost-desirable gambles, this means that
\begin{equation*}
  \text{$f-T^tf\in\desirs$ and $T^tf-f\in\desirs$ for all $f$ in
    $\gambles(\values)$ and all $T$ in $\transfos$}.
\end{equation*}
If we define
\begin{equation*}
  \desirs_\transfos
  :=\set{f-T^tf}{f\in\gambles(\values),T\in\transfos}
  =\set{T^tf-f}{f\in\gambles(\values),T\in\transfos},
\end{equation*}
this leads to the following definition.

\begin{definition}\label{def:strong-invariance}
  A coherent set of almost-desirable gambles $\desirs$ is called \emph{strongly $\transfos$-invariant} if   $f-T^tf\in\desirs$ for all $f$ in $\gambles(\values)$ and all $T$ in $\transfos$, or equivalently, if   $\desirs_\transfos\subseteq\desirs$.
\end{definition}
\noindent
The following theorem gives equivalent characterisations of strong invariance in terms of the alternative types of belief models.

\begin{theorem}\label{theo:strong-invariance}
  Let $\lpr$ be a coherent lower prevision on $\gambles(\values)$, $\desirs$ a coherent set of   almost-desirable gambles, and $\solp$ a weak*-closed convex set of coherent previsions on   $\gambles(\values)$. Assume that these belief models are equivalent, in the sense that they correspond to   one another using the bijective relations in Table~\ref{tab:equivalent-models}.  Then the following   statements are equivalent:
  \begin{enumerate}[1.]
  \item $\desirs$ is strongly $\transfos$-invariant, in the sense that $\desirs_\transfos\subseteq\desirs$;
  \item $\lpr$ is \emph{strongly $\transfos$-invariant}, in the sense that $\lpr(f-T^tf)\geq0$ and     $\lpr(T^tf-f)\geq0$, and therefore $\lpr(f-T^tf)=\lpr(T^tf-f)=0$ for all $f$ in $\gambles(\values)$ and     $T$ in $\transfos$;
  \item $\solp$ is \emph{strongly $\transfos$-invariant}, in the sense that $T\pr=\pr$ for all $\pr$ in     $\solp$ and all $T$ in $\transfos$.\footnote{So strongly invariant belief models correspond to the       \citegen{pericchi1991} notion of a `class of reasonable (or invariant) priors'.}
  \end{enumerate}
\end{theorem}

\begin{proof}
  We give a circular proof. Assume that $\desirs$ is strongly $\transfos$-invariant, and consider any gamble   $f$ on $\values$ and any $T$ in $\transfos$. Then we find for the associated coherent lower prevision $\lpr$   that $\lpr(f-T^tf)=\max\set{s}{f-T^tf-s\in\desirs}\geq0$, and similarly that $\lpr(T^tf-f)\geq0$. But since   $\lpr$ is coherent, we find that also $\lpr(f-T^tf)=-\upr(T^tf-f)\leq-\lpr(T^tf-f)\leq0$ and similarly   $\lpr(T^tf-f)=-\upr(f-T^tf)\leq-\lpr(f-T^tf)\leq0$, whence indeed $\lpr(f-T^tf)=\lpr(T^tf-f)=0$, so the   first statement implies the second.
  \par
  Next, assume that $\lpr$ is strongly $\transfos$-invariant and consider any $\pr$ in the associated set of   dominating coherent previsions $\solp=\set{\pr}{(\forall f)(\pr(f)\geq\lpr(f))}$ and any $T$ in   $\transfos$. Then for any gamble $f$ on $\values$ we see that $\pr(f-T^tf)\geq0$ and $\pr(T^tf-f)\geq0$, and   since $\pr$ is a coherent prevision, this implies that $\pr(T^tf)=\pr(f)$, so indeed $T\pr=\pr$.  Hence, the   second statement implies the third.
  \par
  Finally, assume that $\solp$ is strongly $\transfos$-invariant, and consider any gamble $f$ on $\values$ and   any $T$ in $\transfos$. Then for all $\pr$ in $\solp$ we have that $\pr(f-T^tf)=\pr(T^tf-f)=0$, so both   $f-T^tf$ and $T^tf-f$ belong to the associated set of almost-desirable gambles   $\desirs=\set{g}{(\forall\pr\in\solp)\pr(g)\geq0}$. This tells us that the third statement implies the   first.
\end{proof}

Let us now study in more detail the relationship between weak and strong invariance. First of all, strong invariance implies weak invariance, but generally not the other way around. It is easiest to see this using weak*-closed convex sets of coherent previsions $\solp$. If $\solp$ is strongly $\transfos$-invariant, we have that $T\pr=\pr$ and consequently $T\pr\in\solp$ for all $\pr$ in $\solp$, so $\solp$ is also weakly $\transfos$-invariant. To see that the converse doesn't generally hold, consider the set of all coherent previsions $\linprevs(\values)$ (the vacuous belief model), which is weakly invariant with respect to any monoid of transformations, but not necessarily strongly so, as, unless $\values$ contains only one element, we can easily find transformations $T$ and coherent previsions $\pr$ such that $T\pr$ is different from $\pr$ (also see Theorem~\ref{theo:no-complete-strong-invariance} below).
\par
But the theorem above, when interpreted well, also tells us a number of very interesting things on this issue. First of all, we see that a coherent prevision $\pr$ on $\gambles(\values)$ is strongly $\transfos$-invariant if and only if it is weakly $\transfos$-invariant, so both notions of invariance coincide for coherent previsions. \emph{So anyone who insists on modelling beliefs with Bayesian belief models   (coherent previsions) only, cannot distinguish between the two types of invariance.} This confirms in general what we claimed in the Introduction about Bayesian belief models. From now on, we shall therefore no longer distinguish between strong and weak invariance for coherent previsions, and simply call them \emph{invariant}.
\par
Furthermore, we see that a coherent lower prevision $\lpr$ is strongly $\transfos$-invariant if and only if all its dominating coherent lower previsions are, or equivalently, if all its dominating coherent previsions, i.e., all the coherent previsions in $\solp(\lpr)$, are $\transfos$-invariant.  Or even stronger, it is easy to see that a coherent lower prevision is strongly invariant if and only if it is a lower envelope of some (not necessarily weak*-closed nor convex) set of invariant coherent previsions.
\par
The notions of weak and strong invariance, and the motivation for introducing them, are tailored to the direct behavioural interpretation of lower previsions, or the equivalent belief models. But what happens if we give a lower prevision $\lpr$ a Bayesian sensitivity analysis interpretation? We then hold that there is some actual precise coherent prevision $\pr_a$ modelling the subject's uncertainty about the random variable $\rv$, that we have only imperfect information about in the sense that we only know that $\pr_a\geq\lpr$, or equivalently, that $\pr_a\in\solp(\lpr)$. Assume that we want the imperfect model $\lpr$ to capture that there is `symmetry of evidence' with respect to a monoid of transformations $\transfos$. The actual model $\pr_a$ then should be weakly $\transfos$-invariant, but since this is a (precise) coherent prevision, we can not distinguish between weak and strong invariance, and it should therefore simply be $\transfos$-invariant: $T\pr_a=\pr_a$ for all $T\in\transfos$. Since $\solp(\lpr)$ is interpreted as the set of candidate models for $\pr_a$, all of the coherent previsions $\pr$ in $\solp(\lpr)$ must be $\transfos$-invariant too, or equivalently $\lpr$ must be \emph{strongly} $\transfos$-invariant. A completely analogous course of reasoning shows that if we want $\lpr$ to capture `evidence of symmetry', $\lpr$ must be strongly $\transfos$-invariant as well. So in contradistinction with the direct behavioural interpretation, \emph{on a Bayesian sensitivity analysis   interpretation of $\lpr$, we cannot distinguish between `symmetry of evidence' and `evidence of symmetry',   and strong invariance is the proper symmetry property to use in both cases}.\footnote{See   \cite[Section~9.5]{walley1991} for related comments about the difference between permutability and   exchangeability. These notions will be briefly discussed in Section~\ref{sec:exchangeability}.}
\par
As is the case for weak invariance, a belief model that is strongly $\transfos$-invariant, is also strongly $\transfos'$-invariant for any sub-monoid $\transfos'\subseteq\transfos$. But in contrast with weak invariance, given any monoid $\transfos$, there do not always exist coherent belief models that are strongly invariant with respect to $\transfos$. This is an immediate consequence of the following theorem, which makes an even stronger claim: it is totally \emph{irrational} to require \emph{complete} strong invariance, i.e., strong invariance with respect to the monoid $\transfos_\values$ of all transformations of $\values$.

\begin{theorem}\label{theo:no-complete-strong-invariance}
  Assume that $\values$ contains more than one element.  Then any belief model that is strongly   $\transfos_{\values}$-invariant incurs a sure loss.
\end{theorem}

\begin{proof}
  We shall give a proof for lower previsions. Assume \emph{ex absurdo} that $\lpr$ avoids sure loss, so   $\solp(\lpr)$ is non-empty. Consider any $\pr$ in $\solp(\lpr)$ and any non-constant gamble $f$ on $\values$   [there is at least one such gamble because $\values$ contains more than one element].  This implies that   there are (different) $x_1$ and $x_2$ in $\values$ such that $f(x_1)\not=f(x_2)$.  For any $y$ in $\values$,   consider the transformation $T_y$ that maps all elements of $\values$ to $y$. Then we find that   $T_y^tf=f(y)$, whence $\pr(f(y)-f)\geq\lpr(f(y)-f)\geq0$ and $\pr(f-f(y))\geq\lpr(f-f(y))\geq0$, since   $\lpr$ is by assumption in particular strongly $T_y$-invariant. Consequently $\pr(f)=f(y)$. But this holds   in particular for $y=x_1$ and for $y=x_2$, so we infer that $f(x_1)=\pr(f)=f(x_2)$, a contradiction.
\end{proof}

In fact, we easily see in this proof that given the transformation $T_y$ that maps all elements of $\values$ to $y$, the only strongly $T_y$-invariant belief model that avoids sure loss is the constant prevision on $y$. Consequently, if we consider a monoid $\transfos$ that includes two different constant transformations, any belief model that is strongly $\transfos$-invariant incurs a sure loss.

As a result, we see that there are monoids $\transfos$ for which there are no strongly invariant coherent (lower) previsions. Under which conditions, then, are there strongly $\transfos$-invariant coherent (lower) previsions? It seems easiest, and yields most insight, if we look at this problem in terms of sets of almost-desirable gambles: indeed if we consider a coherent lower prevision $\lpr$ on $\gambles(\values)$, then it is strongly $\transfos$-invariant if and only if for its associated set of almost-desirable gambles $\desirs_\lpr=\set{f\in\gambles(\values)}{\lpr(f)\geq0}$ we have that $\desirs_\transfos\subseteq\desirs_\lpr$.  We can consider $\desirs_\transfos$ itself as a set of almost-desirable gambles, but at this point, we do not know whether $\desirs_\transfos$ is coherent, or whether it even avoids sure loss.  Interestingly, the set of coherent previsions that is associated with $\desirs_\transfos$ is given by
\begin{align*}
  \solp(\desirs_\transfos) &=\set{\pr\in\linprevs(\values)}
  {(\forall g\in\desirs_\transfos)(\pr(g)\geq0)}\\
  &=\set{\pr\in\linprevs(\values)} {(\forall f\in\gambles(\values))(\forall T\in\transfos)     (\pr(f)=\pr(T^tf))}.
\end{align*}
So $\solp(\desirs_\transfos)$ is precisely the convex and weak*-closed set of all $\transfos$-invariant coherent previsions, and $\lpr$ is strongly $\transfos$-invariant if and only if $\solp(\lpr)\subseteq\solp(\desirs_\transfos)$, or in other words, if and only if all coherent previsions that dominate $\lpr$ are $\transfos$-invariant. So there are strongly $\transfos$-invariant coherent lower previsions if and only if $\solp(\desirs_\transfos)\not=\emptyset$, i.e., if there are $\transfos$-invariant coherent previsions, and in this case the lower envelope of $\solp(\desirs_\transfos)$ is the point-wise smallest strongly $\transfos$-invariant coherent lower prevision.
\par
In summary, we see that there are $\transfos$-invariant coherent previsions if and only if the set of almost-desirable gambles $\desirs_\transfos$ avoids sure loss,\footnote{Also see   \citegen[Lemma~3.3.2]{walley1991} Separation Lemma.} which, taking into account~\eqref{eq:avoids-sure-loss-desirs}, is equivalent\footnote{Observe that the set $\desirs_\transfos$ is   a \emph{cone}, i.e., closed under scalar multiplication with non-negative real numbers.} to the condition\footnote{The same condition was derived by \citet[Theorem~3.5.2 and Corollary~3.5.4]{walley1991}   using an argument that works directly with coherent lower previsions. Although our argument strongly plays   on the connection between the three equivalent types of belief models of Table~\ref{tab:equivalent-models},   we believe that it produces more insight, once this connection is fully understood.\label{fn:walley-result}}
\begin{equation}\label{eq:existence}
  \sup\sum_{k=1}^n\left[f_k-T_k^tf_k\right]\geq0
  \quad\text{for all $n\geq0$, $f_1$, \dots, $f_n$ in $\gambles(\values)$ and
    $T_1$, \dots, $T_n$ in $\transfos$.}
\end{equation}
In that case, the natural extension $\edesirs_\transfos:=\edesirs^m_{\desirs_\transfos}$ of $\desirs_\transfos$ to a coherent set of almost-desirable gambles is given by\footnote{Again, observe that   $\desirs_\transfos$ is a cone.}
\begin{equation}\label{eq:invnatex-desirs}
  \edesirs_\transfos
  =\bigcap_{\epsilon>0}
  \set{f\in\gambles(\values)}{f-\epsilon
    \geq\sum_{k=1}^n\left[f_k-T_k^tf_k\right]
    \text{ for some $n\geq0$, $f_k\in\gambles(\values)$, $T_k\in\transfos$}}
\end{equation}
This is the smallest coherent and strongly $\transfos$-invariant set of almost-desirable gambles, or in other words, the belief model that represents evidence of symmetry involving the monoid $\transfos$. The corresponding lower prevision, defined by\footnote{It is easy to see that   $\solp(\desirs_\transfos)=\solp(\edesirs_\transfos)$.}
\begin{align}
  \lnex_\transfos(f) &=\min\set{\pr(f)}{\pr\in\solp(\desirs_\transfos)}
  \label{eq:first-equality}\\
  &=\max\set{\mu\in\reals}{f-\mu\in\edesirs_\transfos}
  \label{eq:second-equality}
\end{align}
is then, by virtue of Eq.~\eqref{eq:first-equality} [see also Theorem~\ref{theo:dominating-invariant-previsions} further on], the point-wise smallest (most conservative) strongly $\transfos$-invariant coherent lower prevision on $\gambles(\values)$, and if we combine Eqs.~\eqref{eq:invnatex-desirs} and~\eqref{eq:second-equality}, we find that\footnote{Again,   \citet[Theorem~3.5.2 and Corollary~3.5.4]{walley1991} proves the same result in a different manner, see also   footnote~\ref{fn:walley-result}.\label{fn:walley-result-too}}
\begin{equation}\label{eq:lnex-transfos}
  \lnex_\transfos(f)
  =\sup\set{\inf\left[f-\sum_{k=1}^n\left[f_k-T_k^tf_k\right]\right]}
  {n\geq0, f_k\in\gambles(\values),T_k\in\transfos}.
\end{equation}
Remember that this lower prevision is only well-defined (assumes finite real values) whenever the condition~\eqref{eq:existence} is satisfied. Taking into account Theorem~\ref{theo:dominating-invariant-previsions} further on, we deduce that a coherent (lower) prevision is (strongly) $\transfos$-invariant if and only if it dominates $\lnex_\transfos$. Also, $\lnex_\transfos$ is the belief model we should use if nothing else but the evidence of symmetry is given. Finally, this formula for the \emph{lower} prevision is constructive, but usually the existence of invariant previsions (on infinite spaces) is proven in a non-constructive (Hahn--Banach) way; see Section~\ref{sec:shift-invariance}, and also \cite{agnew1938} and \citet[Section~2.1.3(8)]{bhaskara1983}. So we cannot usually get to the coherent invariant previsions by construction, but we can always construct their lower envelope explicitly!
\par
We shall have much more to say about the existence of strongly invariant belief models in Section~\ref{sec:strongly-invariant-lpr}, where we show that this existence is guaranteed in particular if the monoid $\transfos$ is Abelian, or if it is a finite group. The following counterexample tells us that there is no such guarantee for infinite groups.

\begin{example}[Permutation invariance on the natural numbers]
  Consider the set $\permuts_\nats$ of all permutations of the set of natural numbers $\nats$. We show that   there are no (strongly) $\permuts_\nats$-invariant coherent (lower) previsions on $\gambles(\nats)$ by   showing that the condition~\eqref{eq:existence} doesn't hold. Indeed, consider the partition of $\nats$ made   up of the sets
  \begin{equation*}
    R_3^r=\set{3n+r}{n\in\nats},\quad r=0,1,2,
  \end{equation*}
  and any permutations $\pi_r$ for $r=0,1,2$ such that for all $n\in\nats$, $\pi_r(n)\in R_3^r$ if and only if   $n\notin R_3^r$ [for instance, let $\pi_r$ be involutive and such that it assigns the first element of   $R_3^r$ to the first of $(R_3^r)^c$, the second element of $R_3^r$ to the second of $(R_3^r)^c$, etc.]   Consider the gamble $G=\sum_{r=0}^2[I_{R_3^r}-\pi_r^tI_{R_3^r}]$ on $\nats$, then we are done if we can show   that $\sup G<0$. Indeed, if $n\in R_3^r$ then $G(n)=1+0+0-(1+1+0)=-1$ for $r=0,1,2$, so $\sup   G=-1$. $\blacklozenge$
\end{example}

These results expose another fundamental difference between weak and strong invariance: while strong invariance with respect to a greater number of transformations means that we must refine our beliefs (i.e, it make them more precise), this is not the case with weak invariance.
\par
On the other hand, strong invariance is preserved by dominating lower previsions: if $\lpr_1$ is a coherent lower prevision that is strongly $\transfos$-invariant and $\lpr_2$ is a coherent lower prevision that dominates $\lpr_1$, then $\lpr_2$ is also strongly $\transfos$-invariant. It indeed seems reasonable that, if a subject has evidence of symmetry, and she has some additional information that allows her to make her judgements more precise, she can add assessments while still preserving strong invariance. But a similar result does not hold for weak invariance: since the vacuous lower prevision is weakly $\transfos_\values$-invariant, this would mean that any lower prevision should be weakly $\transfos_\values$-invariant, \textit{quod non}.
\par
In summary, there is an important conceptual difference between weak and strong invariance. Weakly invariant belief models capture in particular that a subject \emph{has no reason to} strictly prefer a gamble $f$ to its transformation $T^tf$ whenever $f\not>T^tf$. Strong invariance captures that a subject \emph{has reason not   to} distinguish between, i.e., to be indifferent between, the gambles $f$ and $T^tf$. And it is only if you insist on using Bayesian belief models always that you must infer indifference from having no reason to (strictly) prefer. This is of particular relevance for belief models that try to represent a subject's complete ignorance, as we now proceed to show.

\section{Modelling complete ignorance}\label{sec:complete-ignorance}
Suppose our subject is completely ignorant about the value that $\rv$ assumes in $\values$. Then she has no relevant information that would allow her to favour one possible value of $\rv$ over another. This implies that the corresponding belief model should be symmetric in the possible values of $\rv$, or in other words it should be weakly invariant with respect to the group $\permuts_\values$ of all permutations of $\values$. This leads to a form of \citegen[Section~5.5.1]{walley1991} Symmetry Principle.

\begin{SP}
  If a subject is completely ignorant about the value of a random variable $\rv$ in $\values$, then her   corresponding belief model should be weakly invariant with respect to the group $\permuts_\values$ of all   permutations of $\values$.
\end{SP}
\noindent
We have mentioned before that the appropriate belief model for complete ignorance about $\rv$ seems to be the vacuous lower prevision $\lpr_\values$.  But SP by itself is not sufficient to single out this lower prevision: if, for instance, $\values$ is finite, then the uniform precise prevision $\pr^u$, given by
\begin{equation*}
  \pr^u(f)=\frac{1}{\abs{\values}}\sum_{x\in\values}f(x)
\end{equation*}
for each gamble $f$ on $\values$, which assigns equal probability mass $1/\abs{\values}$ to each element of $\values$, is also weakly permutation invariant. We shall also see in Examples~\ref{ex:two-elements} and~\ref{ex:three-elements} of Section~\ref{sec:finite} that there may be many more coherent lower previsions that share the same weak permutation invariance property. If, however, we strengthen the Symmetry Principle to require weak invariance with respect to \emph{all transformations}, and not just all permutations, then Theorem~\ref{theo:vacuous-invariant} tells us that the vacuous lower prevision $\lpr_\values$ is indeed the only coherent lower prevision that is compatible with the following

\begin{SSP}
  If a subject is completely ignorant about the value of a random variable $\rv$ in $\values$, then her   corresponding belief model should be weakly invariant with respect to the monoid $\transfos_\values$ of all   transformations of $\values$.
\end{SSP}

\citet[Section~5.5.1 and note~7 on p.~526]{walley1991} has shown that for random variables $\rv$ taking values in a finite set $\values$, the vacuous lower prevision $\lpr_\values$ is the only coherent lower prevision that is compatible with SP and the so-called\footnote{For additional discussion of this principle, see also   \citet{walley1996b,walley1999}.}

\begin{EP}
  Consider a random variable $\rv$, and consider a set of possible values $A$ for $\rv$. Then the (lower)   probability assigned to the event $A$, i.e., the lower probability that $\rv\in A$, should not depend on the   set $\values$ of all possible values for $\rv$ in which $A$ is embedded.
\end{EP}
\noindent
So under coherence, SSP is equivalent to SP and EP taken together. Under coherence, it is also equivalent to the following rationality principle, as we shall shortly see.

\begin{RPIR}
  If you have two \emph{different} gambles $f$ and $g$ on a random variable $\rv$ that you are completely   ignorant about, then if $f\not\geq g$ you have no reason to prefer $f$ to $g$.
\end{RPIR}
\noindent
Indeed, the only coherent belief model that is compatible with this principle, is the vacuous one. We shall argue in terms of real desirability models\footnote{A similar argument can be given for almost-desirability   models $\desirs$ and lower previsions $\lpr$, using for preference   \citegen[Sections~3.7.7--3.7.9]{walley1991} corresponding notion of \emph{strict preference}, which   corresponds to the present argument by using $\desirs_\lpr^+\cup\{0\}$ as a coherent set of really desirable   gambles.}  $\rdesirs$ (see Section~\ref{sec:really-desirable}). Say that a subject (really) \emph{prefers} $f$ to $g$ whenever $f\not=g$ and $f-g\in\rdesirs$, i.e., she accepts to exchange $g$ for $f$. Then RPIR implies that for all $f\not=0$, $f\not\geq0$ implies that $f\not\in\rdesirs$, or equivalently, by contraposition, that $f\in\rdesirs$ implies $f\geq0$. Hence $\rdesirs=\cone_+$ is the vacuous belief model.
\par
In summary, we have the following equivalences, under coherence, and the only belief model that is compatible with these three equivalent rationality requirements, is the vacuous one:
\begin{equation*}
  \text{SSP}\asa\text{SP+EP}\asa\text{RPIR}.
\end{equation*}
\par
RPIR is a revised version of the Principle of Insufficient Reason (PIR), which states that if you are completely ignorant about the value of a random variable $X$, then you have no reason to distinguish between the different possible values, and therefore should consider all these values to have equal probability. Indeed, from a historical point of view, the PIR was used extensively by Laplace (see for instance \citet{howie2002}) to justify using a uniform probability for modelling complete ignorance.
\par
We are of course aware that our reformulation RPIR of Laplace's PIR is quite unusual and has little or no historical grounds, which is why we refer to it as a \emph{revised}, or perhaps better, improved principle. It might have been preferable to call RPIR the `Principle of Insufficient Reason to Prefer', but we decided against that for aesthetical reasons.
\par
We think that RPIR is reasonable, but that PIR isn't. Indeed, one of the reasons for the critical attitudes of many researchers towards `Bayesian methods' and inverse probability in the nineteenth and early twentieth century seem to lie in the indiscriminate use by many of Laplace's PIR in order to obtain uniform prior probabilities that can be plugged into Bayes's formula.\footnote{An interesting historical discussion of such   attitudes can be found in \citet{howie2002} and \citet{zabell1989a}.} And by `indiscriminate use' we mean precisely the confusion that exists between symmetry of evidence and evidence of symmetry: we have argued that it is only evidence of symmetry that justifies using strongly invariant belief models (and in many cases, such as permutation invariance for finite spaces, strong invariance singles out the uniform probability as the only compatible belief model, see also Section~\ref{sec:finite}). If there is only symmetry of evidence, we should use weakly invariant belief models, and in the special case of complete ignorance, vacuous ones. Of course, as we said in the Introduction and proved in the previous section, for precise previsions (Bayesian belief models) there is no difference between weak and strong invariance, so if you insist on using a Bayesian belief model, symmetry of evidence leads you to a (strongly) invariant one!  The problem with the PIR, therefore, is that the belief model is only allowed to be precise: there would be fewer or no difficulties if in its formulation we just replaced `probability' with `lower and upper probability', for instance.

\section{Weakly invariant lower previsions}
\label{sec:weakly-invariant-lpr} Let us now turn to a more involved mathematical study of the invariance of coherent lower previsions. So far, we have only looked at coherent lower previsions that were defined on all gambles. But of course, it will usually happen that our subject specifies a supremum acceptable buying price $\lpr(f)$ for only a limited number of gambles $f$, say those in a subset $\domain$ of $\gambles(\values)$. And then we can ask ourselves whether such an assessment can be coherently extended to a weakly, or to a strongly, $\transfos$-invariant lower prevision on all gambles. We shall address these, and related, problems in this and the following section. Let us begin here with weak invariance. The following definition generalises the already established notion of weak invariance to lower previsions defined on any $\transfos$-invariant domain, that are not necessarily coherent (they may even incur a sure loss).\footnote{Our notion of weak invariance for a lower prevision is weaker than   \citegen[Section~3.5.1]{walley1991} corresponding notion of invariance, which requires equality, and has the   drawback that it is not preserved by natural extension.}

\begin{definition}[Weak invariance]
  A lower prevision $\lpr$ defined on a set of gambles $\domain\subseteq\gambles(\values)$ is called weakly   \emph{$\transfos$-invariant} if
  \begin{enumerate}[(W1)]
  \item $T^tf\in\domain$ for all $f$ in $\domain$ and $T$ in $\transfos$, i.e., $\domain$ is     $\transfos$-invariant;
  \item $\lpr(T^tf)\geq\lpr(f)$ for all $f$ in $\domain$ and $T$ in $\transfos$, i.e., all $T\lpr$ point-wise     dominate $\lpr$.
  \end{enumerate}
\end{definition}

As before, if $\transfos$ is right-cancellable (and in particular if it is a group), the inequality in the invariance definition is actually an equality: consider a gamble $f$ in $\domain$, a transformation $T$ in $\transfos$ and its right-inverse $R$, we have $\lpr(f)=\lpr((TR)^tf)=\lpr(R^t(T^tf))\geq\lpr(T^tf)$ in addition to $\lpr(T^tf)\geq\lpr(f)$.
\par
Next, because taking convex combinations, lower envelopes, limits inferior and superior preserves inequalities, it is easy to see that convex combinations, lower envelopes and point-wise limits of weakly invariant lower previsions are also weakly invariant. Observe by the way that the same operations also preserve coherence.
\par
The following proposition looks at weak invariance for (precise) previsions.

\begin{proposition}\label{prop:weak-invariance-previsions}
  Let $\pr$ be a prevision, i.e., a self-conjugate lower prevision, defined on a negation-invariant domain   $\domain=-\domain$. Assume that $\domain$ is also $\transfos$-invariant. Then $\pr$ is weakly   $\transfos$-invariant if and only if $\pr(T^tf)=\pr(f)$ for all $T$ in $\transfos$ and all $f$ in $\domain$.
\end{proposition}

\begin{proof}
  It is clear that the condition is sufficient. To show that it is also necessary, assume that $\pr$ is   $\transfos$-invariant, and consider any $T$ in $\transfos$ and any gamble $f$ in $\domain$. Then it follows   from the $\transfos$-invariance of $\pr$ that on the one hand $\pr(T^tf)\geq\pr(f)$, and on the other hand,   since $-f\in\domain$ and $T^t(-f)=-T^tf\in\domain$, that $\pr(-T^tf)=\pr(T^t(-f))\geq\pr(-f)$, or   equivalently, using the self-conjugacy of $\pr$, that $\pr(f)\geq\pr(T^tf)$.
\end{proof}

We study next whether a weakly invariant lower prevision $\lpr$ with domain $\domain$ can be extended to a coherent weakly invariant lower prevision on the set of all gambles, or more generally, whether there is a coherent weakly invariant lower prevision on all gambles that dominates $\lpr$. We already know from the material in Section~\ref{sec:natex-lpr} that a necessary condition for this is that $\lpr$ should avoid sure loss. Indeed, if $\lpr$ incurs sure loss then it has no dominating coherent lower prevision, let alone a weakly invariant one. The perhaps surprising result we prove next is that avoiding sure loss is also sufficient, and that all we have to do is consider the natural extension $\lnex_\lpr$ of $\lpr$, as it preserves weak invariance.  This natural extension is automatically guaranteed to be the point-wise smallest weakly $\transfos$-invariant coherent lower prevision that dominates $\lpr$.\footnote{This result is   mentioned, with only a hint at the proof, by \citet[Theorem~3.5.2]{walley1991}.}

\begin{theorem}[Natural extension preserves weak invariance]
  \label{theo:natext-preserves-invariance}
  The natural extension $\lnex_\lpr$ of a weakly $\transfos$-invariant lower prevision $\lpr$ on a set of   gambles $\domain$ that avoids sure loss is still weakly $\transfos$-invariant, i.e., for all gambles $f$ on   $\values$ and all $T$ in $\transfos$,
  \begin{equation*}
    T\lnex_\lpr(f)=\lnex_\lpr(T^tf)\geq\lnex_\lpr(f).
  \end{equation*}
  Consequently, $\lnex_\lpr$ is the point-wise smallest weakly $\transfos$-invariant coherent lower prevision   on $\gambles$ that dominates $\lpr$ on its domain $\domain$.
\end{theorem}

\begin{proof}
  Consider any gamble $f$ on $\values$ and any $T$ in $\transfos$.  From the definition~\eqref{eq:natex-lpr}   of natural extension, and the fact that $T^t\domain\subseteq\domain$, we get
  \begin{align}
    \lnex_\lpr(T^tf) &=\sup_{\substack{\lambda_k\geq0,f_k\in\domain\\k=1\dots,n,n\geq0}}     \set{\alpha}{T^tf-\alpha\geq\sum_{k=1}^n
      \lambda_k\left[f_k-\lpr(f_k)\right]}\notag\\
    &\geq\sup_{\substack{\lambda_k\geq0,g_k\in\domain\\k=1\dots,n,n\geq0,}}     \set{\alpha}{T^tf-\alpha\geq\sum_{k=1}^n       \lambda_k\left[T^tg_k-\lpr(T^tg_k)\right]}.\label{eq:first-inequality}
  \end{align}
  Now it follows from the $T$-invariance of $\lpr$ that $\lpr(T^tg_k)\geq\lpr(g_k)$, whence
  \begin{equation*}
    \sum_{k=1}^n\lambda_k\left[T^tg_k-\lpr(T^tg_k)\right]
    \leq T^t\sum_{k=1}^n\lambda_k\left[g_k-\lpr(g_k)\right],
  \end{equation*}
  and consequently $f-\alpha\geq\sum_{k=1}^n\lambda_k[g_k-\lpr(g_k)]$ implies that
  \begin{equation*}
    T^tf-\alpha\geq T^t\sum_{k=1}^n\lambda_k\left[g_k-\lpr(g_k)\right]
    \geq\sum_{k=1}^n\lambda_k\left[T^tg_k-\lpr(T^tg_k)\right].
  \end{equation*}
  So we may infer from the inequality~\eqref{eq:first-inequality} that
  \begin{equation*}
    \lnex_\lpr(T^tf)
    \geq\sup_{\substack{\lambda_k\geq0,g_k\in\domain\\k=1\dots,n,n\geq0,}}
    \set{\alpha}{f-\alpha\geq\sum_{k=1}^n\lambda_k\left[g_k-\lpr(g_k)\right]}
    =\lnex_\lpr(f),
  \end{equation*}
  which completes the proof.
\end{proof}

Hence, if we start out with a lower prevision $\lpr$ on $\domain$ that is weakly $\transfos$-invariant and already coherent, then its natural extension $\lnex_\lpr$ is the smallest coherent and weakly $\transfos$-invariant lower prevision on all gambles that agrees with $\lpr$ on $\domain$. As we shall show further on, this result does not carry over to strong invariance.

\section{Strongly invariant lower previsions}
\label{sec:strongly-invariant-lpr}
We now turn to the study of strong invariance for lower previsions on general domains.

\subsection{Definition and immediate properties}
The following definition generalises the notion of strong invariance introduced in Section~\ref{sec:strong-invariance} to lower previsions that needn't be coherent, nor defined on all of $\gambles(\values)$.

\begin{definition}[Strong invariance]
  A lower prevision $\lpr$ defined on a set of gambles $\domain\subseteq\gambles(\values)$ is called   \emph{strongly $\transfos$-invariant} if
  \begin{enumerate}[(S1)]
  \item $T^tf-f\in\domain$ and $f-T^tf\in\domain$ for all $f$ in $\domain$ and all $T\in\transfos$;
  \item $\lpr(T^tf-f)\geq0$ and $\lpr(f-T^tf)\geq0$ for all $f$ in $\domain$ and all $T\in\transfos$.
  \end{enumerate}
\end{definition}

As is the case for weak invariance, it is easy to see that strong $\transfos$-invariance is preserved under convex combinations, lower envelopes, and point-wise limits, simply because all these operations preserve inequalities.

\begin{proposition}\label{prop:strong-implies-weak}
  A strongly $\transfos$-invariant coherent lower prevision on a $\transfos$-invariant domain is also weakly   $\transfos$-invariant.
\end{proposition}

\begin{proof}
  First of all, the coherence and strong invariance of $\lpr$ imply that $0\leq   \lpr(T^tf-f)\leq\lpr(T^tf)-\lpr(f)$, whence $\lpr(T^tf)\geq\pr(f)$ and similarly, we derive from   $\lpr(f-T^tf)\geq0$ that $\lpr(f)\geq\lpr(T^tf)$.  So we see that $\lpr$ is also weakly   $\transfos$-invariant (with equality).
\end{proof}
\noindent
To see that a converse result does not generally hold, so weak invariance is actually weaker than strong invariance, consider the vacuous lower prevision $\lpr_\values$ on $\gambles(\values)$ and the transformation $T_{x_0}$ that maps all elements $x$ of $\values$ to $x_0$. Then, for any gamble $f$ such that $\inf f<f(x_0)$ we have $\lpr_\values(f-T_{x_0}^tf)<0$. Hence, $\lpr_\values$ is not strongly $T_{x_0}$-invariant but Theorem~\ref{theo:vacuous-invariant} implies that it is weakly $T_{x_0}$-invariant. If we consider a finite space $\values$ and the vacuous lower prevision $\lpr_\values$ on $\gambles(\values)$ and the class $\permuts_\values$ of all permutations of $\values$, we can see that weak invariance (with equality) does not imply strong invariance.
\par
So weak invariance is indeed a weaker notion than strong invariance. The following theorem expresses the main difference between these two concepts: while the former means that the set of coherent previsions $\solp(\lpr)$ is invariant, the latter means that every element of this set is invariant.

\begin{theorem}\label{theo:dominating-invariant-previsions}
  Let $\domain$ be a negation invariant and $\transfos$-invariant set of gambles such that $T^tf-f$ is in   $\domain$ for all $f$ in $\domain$ and $T$ in $\transfos$.
  \begin{enumerate}[1.]
  \item A coherent prevision $\pr$ on $\domain$ is weakly $\transfos$-invariant if and only if it is strongly     $\transfos$-invariant.  In either case we simply call it \emph{$\transfos$-invariant}.
  \item A coherent lower prevision $\lpr$ on $\domain$ is strongly $\transfos$-invariant if and only if all     its dominating coherent previsions are (strongly) $\transfos$-invariant on $\domain$.
  \end{enumerate}
\end{theorem}

\begin{proof}
  We start with the first statement. We only need to prove the direct implication, so assume that $\pr$ is   weakly $\transfos$-invariant, and consider any $f$ in $\domain$. Then from the assumption and   Proposition~\ref{prop:weak-invariance-previsions} we get $\pr(T^tf)=\pr(f)$, and it follows from the   linearity of $\pr$ that indeed $\pr(T^tf-f)=\pr(f-T^tf)=0$.
  \par
  We now turn to the second statement. Since any coherent lower prevision is the lower envelope of its   dominating coherent previsions, the converse implications follow at once, since taking a lower envelope   preserves strong invariance. To prove the direct implication, assume that $\lpr$ is strongly   $\transfos$-invariant, and consider any coherent prevision $\pr$ in $\solp(\lpr)$. For any $T$ in   $\transfos$ and any $f$ in $\domain$ we then find that
  \begin{equation*}
    0\leq\lpr(f-T^tf)\leq\pr(f-T^tf)=-\pr(T^tf-f)\leq-\lpr(T^tf-f)\leq0,
  \end{equation*}
  whence indeed $\pr(f)=\pr(T^tf)$.
\end{proof}

\subsection{Strongly invariant natural extension}
\label{sec:invnatex-strong} 
We have shown when studying weak invariance that for any weakly $\transfos$-invariant lower prevision $\lpr$ on some domain $\domain$ that avoids sure loss, there is a point-wise smallest weakly invariant coherent lower prevision defined on all gambles that dominates it: its natural extension $\lnex_\lpr$. Let us now investigate whether something similar can be done for the notion of strong invariance.  The question then is: Consider a monoid $\transfos$ of transformations of $\values$ and a lower prevision $\lpr$ on $\domain$ that avoids sure loss, are there strongly $\transfos$-invariant coherent lower previsions on all $\gambles(\values)$ that dominate $\lpr$, and if so, what is the point-wise smallest such lower prevision?  Let us denote, as before, by
\begin{equation*}
  \desirs_\lpr=\set{f\in\gambles(\values)}{\lnex_\lpr(f)\geq0}
\end{equation*}
the set of almost-desirable gambles associated with $\lpr$, and by
\begin{equation*}
  \solp(\lpr)=\set{\pr\in\linprevs(\values)}
  {(\forall f\in\domain)(\pr(f)\geq\lpr(f))}
\end{equation*}
its set of dominating coherent previsions, then clearly a coherent lower prevision $\alpr$ on $\gambles(\values)$ is strongly $\transfos$-invariant and dominates $\lpr$ if and only if $\solp(\alpr)\subseteq\solp(\lpr)\cap\solp(\desirs_\transfos)$, or equivalently, $\desirs_\lpr\cup\desirs_\transfos\subseteq\desirs_\alpr$. So there are strongly $\transfos$-invariant coherent (lower) previsions that dominate $\lpr$ if and only if $\solp(\lpr)\cap\solp(\desirs_\transfos)\not=\emptyset$, or equivalently, if the set of almost-desirable gambles $\desirs_\lpr\cup\desirs_\transfos$ avoids sure loss, and in this case the lower envelope of $\solp(\lpr)\cap\solp(\desirs_\transfos)$, or equivalently, the lower prevision associated with the natural extension of the set of almost-desirable gambles $\desirs_\lpr\cup\desirs_\transfos$, is the smallest such lower prevision. In the language of coherent lower previsions, this leads to the following theorem.\footnote{\citet[Theorems~3.5.2 and~3.5.3]{walley1991} proves similar results involving   Eqs.~\eqref{eq:existence-with-transfos} and~\eqref{eq:natex-with-transfos} for what we call weakly   $\transfos$-invariant $\lpr$ that avoid sure loss, in a different manner.  See also   footnotes~\ref{fn:walley-result} and~\ref{fn:walley-result-too}.}

\begin{theorem}[Strongly invariant natural extension]\label{theo:invnatext}
  Consider a lower prevision $\lpr$ on $\domain$ that avoids sure loss, and a monoid $\transfos$ of   transformations of $\values$. Then there are strongly $\transfos$-invariant coherent (lower) previsions on   $\gambles(\values)$ that dominate $\lpr$ on $\domain$ if and only if
  \begin{equation}\label{eq:existence-with-transfos}
    \unex_\lpr\left(\sum_{k=1}^n\left[f_k-T_k^tf_k\right]\right)\geq0
    \quad\text{for all $n\geq0$, $f_1$, \dots, $f_n$ in $\gambles(\values)$
      and $T_1$, \dots, $T_n$ in $\transfos$},
  \end{equation}
  or equivalently, if
  \begin{equation}\label{eq:existence-with-lpr}
    \unex_\transfos\left(\sum_{k=1}^n\lambda_k\left[f_k-\lpr(f_k)\right]\right)\geq0
    \quad\text{for all $n\geq0$, and $f_1$, \dots, $f_n$ in $\domain$.}
  \end{equation}
  In that case the smallest coherent and strongly $T$-invariant lower prevision on $\gambles(\values)$ that   dominates $\lpr$ on its domain $\domain$ is given by
  \begin{align}
    \lnex_{\lpr,\transfos}(f) &=\sup\set{\lnex_\lpr \left(f-\sum_{k=1}^n\left[f_k-T_k^tf_k\right]\right)}     {n\geq0, f_k\in\gambles(\values),T_k\in\transfos}
    \label{eq:natex-with-transfos}\\
    &=\sup\set{\lnex_\transfos \left(f-\sum_{k=1}^n\lambda_k\left[f_k-\lpr(f_k)\right]\right)}     {n\geq0,f_k\in\domain,\lambda_k\geq0}
    \label{eq:natex-with-lpr}
  \end{align}
  for all gambles $f$ on $\values$; and $\solp(\lnex_{\lpr,\transfos})$ is the set of all   $\transfos$-invariant coherent previsions that dominate $\lpr$ on $\domain$.
\end{theorem}

\begin{proof}
  We already know that there is a dominating coherent (lower) prevision if and only if   $\solp(\lpr)\cap\solp(\desirs_\transfos)$ is non-empty.  Let us show that this is equivalent to the   conditions~\eqref{eq:existence-with-transfos} and~\eqref{eq:existence-with-lpr}. To see the equivalence   between these two conditions, it suffices to notice [use Eq.~\eqref{eq:natex-lpr} and the fact that   $\unex_\lpr(h)=-\lnex_\lpr(-h)$] that condition~\eqref{eq:existence-with-transfos} is equivalent to
  \begin{multline}\label{eq:asl-intermediate}
    \sup\left[\sum_{k=1}^n\left[f_k-T_k^tf_k\right]+
      \sum_{j=1}^{m}\left[g_j-\lpr(g_j)\right]\right]\geq0\\
    \text{for all $n,m\geq0$, $f_k\in\gambles(\values)$, $T_k\in\transfos$, $g_j\in\domain$},
  \end{multline}
  and that this is in turn [use Eq.~\eqref{eq:lnex-transfos} and the fact that   $\unex_\transfos(h)=-\lnex_\transfos(-h)$] equivalent to condition~\eqref{eq:existence-with-lpr}.
  But, considering condition~\eqref{eq:avoids-sure-loss-desirs}, we see that   condition~\eqref{eq:asl-intermediate} holds if and only if the set of almost-desirable gambles   $\desirs_\lpr\cup\desirs_\transfos$ avoids sure loss, or equivalently, if the corresponding set of coherent   previsions $\solp(\lpr)\cap\solp(\desirs_\transfos)$ is non-empty.
  \par
  We now prove the validity of the expression~\eqref{eq:natex-with-lpr} for the lower envelope   $\lnex_{\lpr,\transfos}$ of the set of coherent previsions $\solp(\lpr)\cap\solp(\desirs_\transfos)$. The   proof for the expression~\eqref{eq:natex-with-transfos} is analogous. We know from the material in   Section~\ref{sec:ip-models} that this lower envelope is also the coherent lower prevision associated with   the natural extension of the set of almost-desirable gambles $\desirs_\lpr\cup\desirs_\transfos$, so we get   by applying Eq.~\eqref{eq:natex-lpr} with $\desirs=\desirs_\lpr\cup\desirs_\transfos$ that
  \begin{align*}
    \lnex_{\lpr,\transfos}(f)
    &=\sup_{\substack{\lambda_k\geq0,g_k\in\desirs_\lpr\\
        k=1,\dots,n,n\geq0}}\quad
    \sup_{\substack{\mu_\ell\geq0,h_\ell\in\desirs_\transfos\\
        \ell=1,\dots,m,m\geq0}} \inf\left[f-\sum_{k=1}^n\lambda_kg_k
      -\sum_{\ell=1}^m\mu_\ell h_\ell\right]\\
    &=\sup_{\substack{\lambda_k\geq0,g_k\in\desirs_\lpr\\
        k=1,\dots,n,n\geq0}}\quad
    \sup_{\substack{\mu_\ell\geq0,h_\ell\in\desirs_\transfos\\
        \ell=1,\dots,m,m\geq0}} \inf\left[\left(f-\sum_{k=1}^n\lambda_kg_k\right)
      -\sum_{\ell=1}^m\mu_\ell h_\ell\right]\\
    &=\sup_{\substack{\lambda_k\geq0,g_k\in\desirs_\lpr\\
        k=1,\dots,n,n\geq0}} \lnex_\transfos\left(f-\sum_{k=1}^n\lambda_kg_k\right)
    =\sup_{\substack{\lambda_k\geq0,f_k\in\domain\\
        k=1,\dots,n,n\geq0}} \lnex_\transfos\left(f-\sum_{k=1}^n\lambda_k[f_k-\lpr(f_k)]\right),
  \end{align*}
  for every gamble $f$ on $\values$, also taking into account the definition~\eqref{eq:lnex-transfos} of   $\lnex_\transfos$.
\end{proof}

In conclusion, whenever the equivalent conditions~\eqref{eq:existence-with-transfos} and~\eqref{eq:existence-with-lpr} are satisfied for a lower prevision $\lpr$ that avoids sure loss, then (and only then) the functional $\lnex_{\lpr,\transfos}$, defined by Eqs.~\eqref{eq:natex-with-transfos} and~\eqref{eq:natex-with-lpr}, is the point-wise smallest coherent and strongly $\transfos$-invariant lower prevision that dominates $\lpr$. We shall call $\lnex_{\lpr,\transfos}$ the \emph{strongly   $\transfos$-invariant natural extension} of $\lpr$, as it is the belief model that the assessments captured in $\lpr$ lead to if in addition a (so-called structural)\footnote{Structural assessments are discussed in   general in \citet[Chapter~9]{walley1991}.} assessment of symmetry involving the monoid $\transfos$ is made.

\subsection{The existence of strongly invariant coherent (lower) previsions}
There is a beautiful and surprisingly simple argument to show that for some types of monoids $\transfos$, there always are strongly $\transfos$-invariant lower previsions that dominate a given lower prevision that is weakly $\transfos$-invariant and avoids sure loss.  It is based on the combination of a number of ideas in the literature: (i) \citet[Section~2]{agnew1938} constructed some specific type of Minkowski functional and used this together with a Hahn--Banach extension result to prove the existence of linear functionals that are invariant with respect to certain groups of permutations; (ii) \citet[Theorem~3]{day1942} showed, in a discussion of ergodic theorems, that a similar construction always works for Abelian semigroups of transformations; (iii) with crucially important insight, \citet[Theorems~3.5.2 and~3.5.3]{walley1991} recognised that the Minkowski functional in the existence proofs of \citeauthor{agnew1938}, and \citeauthor{day1942}, is actually what we have called a strongly invariant lower prevision, and he used the ideas behind this construction to introduce what we shall call \emph{mixture lower previsions} in Section~\ref{sec:mixture-lower-previsions}; (iv) in another seminal discussion of mean ergodic theorems, \citet{alaoglu1940} show that (Moore--Smith-like) convergence of convex mixtures of linear transformations is instrumental in characterising ergodicity; and (v) \citet[Section~2.1.3]{bhaskara1983} use so-called Banach limits to generate shift-invariant probability charges. In this and the next section, we combine and extend these ideas to prove more general existence results for (strongly) invariant coherent (lower) previsions, and to investigate their relation to (generalised) Banach limits (Section~\ref{sec:shift-invariance}). As we shall see in Section~\ref{sec:mixture-lower-previsions}, \citegen[Section~3.5]{walley1991} results can then be derived from our more general treatment.
\par
Consider a monoid $\transfos$ of transformations of $\values$.  We can, as before, consider the set of lifted transformations $\transfos^t$ as a monoid of linear transformations of the linear space $\gambles(\values)$. A \emph{convex combination} $T^*$ of elements of $\transfos^t$ is a linear transformation of $\gambles(\values)$ of the form
\begin{equation*}
  T^*=\sum_{k=1}^n\lambda_kT_k^t,
\end{equation*}
where $n\geq1$, $\lambda_1$, \dots, $\lambda_n$ are non-negative real numbers that sum to one, and of course $T^*f=\sum_{k=1}^n\lambda_kT_k^tf$. We denote by $\transfos^*$ the set of all convex combinations of elements of $\transfos^t$.  We have of course for any two elements $T^*_1=\sum_{k=1}^m\lambda_kU_k^t$ and $T^*_2=\sum_{k=1}^n\mu_kV_k^t$ of $\transfos^*$ that their composition
\begin{equation*}
  T_2^*T_1^*
  =\sum_{k=1}^n\mu_kV_k^t\left(\sum_{\ell=1}^m\lambda_\ell U_\ell^t\right)
  =\sum_{k=1}^n\sum_{\ell=1}^m\mu_k\lambda_\ell V_k^tU_\ell^t
  =\sum_{k=1}^n\sum_{\ell=1}^m\lambda_\ell\mu_k(U_\ell V_k)^t
\end{equation*}
again belongs to $\transfos^*$. This implies that $\transfos^*$ is a monoid of linear transformations of $\gambles(\values)$ as well. We can now introduce invariance definitions involving transformations in $\transfos^*$ in precisely the same way as we defined them for $\transfos$ (or actually $\transfos^t$).  We can also define, for any real functional $\Lambda$ and $T^*\in\transfos^*$, the transformed functional $T^*\Lambda$ as $\Lambda\circ T^*$. We then have the following result.

\begin{proposition}\label{prop:convexification}
  The following statements hold, where $f$ is a gamble on $\values$, $\domain$ is a convex set of gambles on   $\values$, and $\lpr$ is a coherent lower prevision on $\domain$:
  \begin{enumerate}[1.]
  \item $f$ is $\transfos$-invariant if and only if $f$ is $\transfos^*$-invariant;
  \item $\domain$ is $\transfos$-invariant if and only if $\domain$ is $\transfos^*$-invariant;
  \item $\lpr$ is weakly $\transfos$-invariant if and only if $\lpr$ is weakly $\transfos^*$-invariant;
  \item $\lpr$ is strongly $\transfos$-invariant if and only if $\lpr$ is strongly $\transfos^*$-invariant.
  \end{enumerate}
\end{proposition}

\begin{proof}
  It suffices of course to prove the direct implications. Consider an arbitrary   $T^*=\sum_k\lambda_kT_k\in\transfos^*$. For the first statement, let $f$ be $\transfos$-invariant, then   $T^*f=\sum_k\lambda_kT_k^tf=\sum_k\lambda_kf=f$, where the second equality follows from the   $\transfos$-invariance of $f$. So $f$ is $\transfos^*$-invariant. For the second statement, let $\domain$ be   $\transfos$-invariant and let $f\in\domain$, then $T^*f=\sum_k\lambda_kT_k^tf\in\domain$, because   $T_k^tf\in\domain$ for all $k$ by the $\transfos$-invariance of $\domain$ and because $\domain$ is   convex. So $\domain$ is $\transfos^*$-invariant. For the third statement, assume that $\lpr$ is weakly   $\transfos$-invariant. For any $f\in\domain$,
  \begin{equation*}
    \lpr(T^*f)
    =\lpr\left(\sum_k\lambda_kT_k^tf\right)
    \geq\sum_k\lambda_k\lpr(T_k^tf)
    \geq\sum_k\lambda_k\lpr(f)
    =\lpr(f),
  \end{equation*}
  where the first inequality follows from the coherence of $\lpr$, and the second from the weak   $\transfos$-invariance of $\lpr$. Hence $\lpr$ is weakly $\transfos^*$-invariant. For the last statement,   assume that $\lpr$ is strongly $\transfos$-invariant. For any $f\in\domain$,
  \begin{equation*}
    \lpr\left(\sum_k\lambda_kT_k^tf-f\right)
    =\lpr\left(\sum_k\lambda_k(T_k^tf-f)\right)
    \geq\sum_k\lambda_k\lpr(T_k^tf-f)
    \geq0,
  \end{equation*}
  where the first inequality follows from the coherence of $\lpr$, and the second from the strong   $\transfos$-invariance of $\lpr$. Similarly $\lpr(f-\sum_k\lambda_kT_k^tf)\geq0$.  Hence $\lpr$ is strongly   $\transfos^*$-invariant.
\end{proof}

We now define the following binary relation $\succeeds$ on $\transfos^*$: for $T_1^*$ and $T_2^*$ in $\transfos^*$ we say that $T_2^*$ \emph{is a successor} of $T_1^*$, and we write $T_2^*\succeeds T_1^*$, if and only if there is some $T^*$ in $\transfos^*$ such that $T_2^*=T^*T_1^*$. Clearly $\succeeds$ is a reflexive and transitive relation, because $\transfos^*$ is a monoid. We say that $\transfos^*$ has the \emph{Moore--Smith property}, or is \emph{directed by $\succeeds$}, if any two elements of $\transfos^*$ have a common successor, i.e., for any $T_1^*$ and $T_2^*$ in $\transfos^*$ there is some $T^*$ in $\transfos^*$ such that $T^*\succeeds T_1^*$ and $T^*\succeeds T_2^*$. It is not difficult to see that if $\transfos$ is Abelian, or a finite group, then $\transfos^*$ is directed by the successor relation. This need not hold if $\transfos$ is an infinite group or a finite monoid, however.
\par
Now, given a \emph{net} $\alpha$ on $\transfos^*$, i.e., a mapping $\alpha\colon\transfos^*\to\mathbb{R}$, we can take the \emph{Moore--Smith limit} of $\alpha$ with respect to the directed set $(\transfos^*,\succeeds)$ \citep[Section~I, p.~103]{moore1922}, which, if it exists, is uniquely defined as the real number $a$ such that, for every $\epsilon>0$, there is a $T^*_\epsilon$ in $\transfos^*$, such that $|\alpha(T^*)-a|<\epsilon$ for all $T^*\succeeds T^*_\epsilon$.  The Moore--Smith limit $a$ of $\alpha$ is denoted by $\lim_{T^*\in\transfos^*}\alpha(T^*)$. This limit always exists if $\alpha$ is non-decreasing and bounded from above, or if $\alpha$ is non-increasing and bounded from below.

\begin{theorem}\label{theo:Moore-Smith}
  Let $\lpr$ be a coherent and weakly $\transfos$-invariant lower prevision on $\gambles(\values)$, and assume   that $\transfos^*$ has the Moore--Smith property. Then for any gamble $f$ on $\values$ the Moore--Smith   limit $\lim_{T^*\in\transfos^*}\lpr(T^*f)$ converges to a real number $\alpr_{\lpr,\transfos}(f)$. Moreover,   $\alpr_{\lpr,\transfos}$ is the point-wise smallest strongly $\transfos$-invariant coherent lower prevision   on $\gambles(\values)$ that dominates $\lpr$ on $\gambles(\values)$, and
  \begin{equation}\label{eq:Q-functional}
    \alpr_{\lpr,\transfos}(f)
    =\sup\set{\lpr(T^*f)}{T^*\in\transfos^*}
    =\sup\set{\lpr\left(\frac{1}{n}\sum_{k=1}^{n}T_k^tf\right)}
    {n\geq1, T_1,\dots,T_n\in\transfos}.
  \end{equation}
\end{theorem}

\begin{proof}
  First, fix $f$ in $\gambles(\values)$. Consider $T_1^*$ and $T_2^*$ in $\transfos^*$, and assume that   $T_2^*\succeeds T_1^*$. This means that there is some $T^*$ in $\transfos^*$ such that $T_2^*=T^*T_1^*$, and   consequently we find that
  \begin{equation*}
    \lpr(T_2^*f)=\lpr(T^*(T_1^*f))\geq\lpr(T_1^*f),
  \end{equation*}
  where the inequality follows from the fact that $\lpr$ is in particular weakly $T^*$-invariant [observe that   $\gambles(\values)$ is convex and that $\lpr$ is weakly $\transfos$-invariant, and apply   Proposition~\ref{prop:convexification}].  This means that the net $\lpr(T^*f)$, $T^*\in\transfos^*$ is   non-decreasing.  Since this net is moreover bounded from above [by $\sup f$, since $\lpr$ is coherent], it   converges to a real number $\alpr_{\lpr,\transfos}(f)$, and clearly
  \begin{equation}\label{eq:ms-intermediate}
    \alpr_{\lpr,\transfos}(f)
    =\lim_{T^*\in\transfos^*}\lpr(T^*f)
    =\sup\set{\lpr(T^*f)}{T^*\in\transfos^*}.
  \end{equation}
  This tells us that the net of coherent lower previsions $T^*\lpr$, $T^*\in\transfos^*$ converges point-wise   to the lower prevision $\alpr_{\lpr,\transfos}$, so $\alpr_{\lpr,\transfos}$ is a coherent lower prevision   as well [taking a point-wise limit preserves coherence].  Since $\iden_\values^t\in\transfos^*$, it follows   from Eq.~\eqref{eq:ms-intermediate} that $\alpr_{\lpr,\transfos}(f)\geq\lpr(\iden_\values^tf)=\lpr(f)$, so   $\alpr_{\lpr,\transfos}$ dominates $\lpr$ on $\gambles(\values)$. We now show that $\alpr_{\lpr,\transfos}$   is strongly $\transfos$-invariant.\footnote{The idea for this part of the proof is due to \citet[Point (iv)     of the proof of Theorem~3.5.3]{walley1991}.}  Consider any $f$ in $\gambles(\values)$ and $T$ in   $\transfos$. Then for any $n\geq1$, $T^*_n:=\frac{1}{n}\sum_{k=1}^n(T^k)^t$ belongs to $\transfos^*$, and it   follows from the coherence of $\lpr$ that
  \begin{multline*}
    \lpr(T^*_n(f-T^tf)) =\frac{1}{n}\lpr(T^tf-(T^{n+1})^tf)
    \geq\frac{1}{n}\inf\left[T^tf-(T^{n+1})^tf\right]\\
    =-\frac{1}{n}\sup\left[(T^{n+1})^tf-T^tf\right] \geq-\frac{2}{n}\sup\lvert f\rvert,
  \end{multline*}
  and consequently
  \begin{equation*}
    \alpr_{\lpr,\transfos}(f-T^tf)
    \geq\sup\set{-\frac{2}{n}\sup\lvert f\rvert}{n\geq1}=0.
  \end{equation*}
  A similar argument can be given for $\alpr_{\lpr,\transfos}(T^tf-f)\geq0$, so $\alpr_{\lpr,\transfos}$ is   indeed strongly $\transfos$-invariant.
  \par
  Next, consider any strongly $\transfos$-invariant and coherent lower prevision $\alpr$ on   $\gambles(\values)$, and assume that it dominates $\lpr$. Then we get for any gamble $f$ on $\values$ and   any $T^*$ in $\transfos^*$:
  \begin{equation*}
    \alpr(f)
    =\alpr(f-T^*f+T^*f)
    \geq\alpr(f-T^*f)+\alpr(T^*f)
    \geq\alpr(T^*f)\geq\lpr(T^*f),
  \end{equation*}
  where the first inequality follows from the coherence of $\alpr$, the second inequality from its strong   $\transfos$-invariance [use Proposition~\ref{prop:convexification}], and the last inequality from the fact   that $\alpr$ dominates $\lpr$. We then deduce from Eq.~\eqref{eq:ms-intermediate} that $\alpr$ dominates   $\alpr_{\lpr,\transfos}$. So $\alpr_{\lpr,\transfos}$ is indeed the point-wise smallest strongly   $\transfos$-invariant coherent lower prevision on $\gambles(\values)$ that dominates $\lpr$ on   $\gambles(\values)$.
  \par
  Finally, let us prove the second equality in Eq.~\eqref{eq:Q-functional}.  Consider a gamble $f$ and any   $\epsilon>0$. Then, by Eq.~\eqref{eq:ms-intermediate}, there is some $T^*$ in $\transfos^*$ such that   $\alpr_{\lpr,\transfos}(f)\leq\lpr(T^*f)+\frac{\epsilon}{2}$.  For this $T^*$, there are $n\geq1$, $T_1$,   \dots, $T_n$ in $\transfos$ and $\lambda_1,\dots,\lambda_n\geq 0$ that sum to one, such that   $T^*=\sum_{k=1}^{n}\lambda_kT_k^t$. Let $\rho_1$, \dots, $\rho_n$ be non-negative rational numbers   satisfying $\abs{\rho_i-\lambda_i}\leq\frac{\epsilon}{2n\sup\abs{f}}$ such that moreover   $\sum_{i=1}^{n}\rho_i=1$.\footnote{To see that such rational numbers exist, it suffices to consider     non-negative rational numbers $\rho_1,\dots,\rho_{n-1}$ such that $0\leq\rho_i\leq\lambda_i\leq1$ and     $\abs{\rho_i-\lambda_i}\leq\frac{\epsilon}{2n^2\sup\abs{f}}$ for $i=1,\dots,n-1$, and to let     $\rho_n:=1-\sum_{i=1}^{n-1}\rho_i\geq1-\sum_{i=1}^{n-1}\lambda_i =\lambda_n\geq0$.  Then $\rho_n\in[0,1]$,     and for $n$ big enough, and unless we are in the trivial case where $\lambda_i=1$ for some $i$, we get     $\abs{\rho_n-\lambda_n}\leq\frac{\epsilon}{2n\sup\abs{f}}$.} Now it follows from the coherence of $\lpr$   that
  \begin{equation*}
    \lpr(T^* f)
    =\lpr\left(\sum_{i=1}^{n}\lambda_iT_i^tf\right)
    \leq\lpr\left(\sum_{i=1}^{n}\rho_iT_i^tf\right)
    -\lpr\left(\sum_{i=1}^{n}(\rho_i-\lambda_i)T_i^tf\right),
  \end{equation*}
  and also
  \begin{multline*}
    \lpr\left(\sum_{i=1}^{n}(\rho_i-\lambda_i)T_i^tf\right) \geq\sum_{i=1}^{n}\lpr((\rho_i-\lambda_i)T_i^tf)
    \geq\sum_{i=1}^{n}\inf(\rho_i-\lambda_i)T_i^tf\\
    \geq\sum_{i=1}^{n}-\frac{\epsilon}{2n\sup\abs{f}}\sup\abs{f} =-\frac{\epsilon}{2},
  \end{multline*}
  whence
  \begin{equation*}
    \alpr_{\lpr,\transfos}(f)
    \leq\lpr(T^*f)+\frac{\epsilon}{2}
    \leq\lpr\left(\sum_{i=1}^{n}\rho_iT_i^tf\right)+\epsilon,
  \end{equation*}
  and consequently
  \begin{equation*}
    \alpr_{\lpr,\transfos}(f)
    =\sup\set{\lpr\left(\sum_{i=1}^{n}\rho_iT_i^tf\right)}
    {n\geq1,T_1,\dots,T_n\in\transfos,\rho_1,\dots,\rho_n\in\mathbb{Q}^+,
      \sum_{i=1}^{n}\rho_i=1},
  \end{equation*}
  where $\mathbb{Q}^+$ denotes the set of non-negative rational numbers. Now, it is easy to see [just consider   the least common multiple of the denominators of $\rho_1$, \dots, $\rho_n$] that this supremum coincides   with the right-hand side of Eq.~\eqref{eq:Q-functional}.
\end{proof}

This result allows us to establish the following corollary. It gives a sufficient condition for the existence of strongly $\transfos$-invariant lower previsions dominating a given coherent lower prevision $\lpr$. The smallest such lower prevision reflects how initial behavioural dispositions, reflected in $\lpr$, are modified (strengthened) to $\lnex_{\lpr,\transfos}$ when we add the extra assessment of strong invariance with respect to a monoid $\transfos$ of transformations.

\begin{corollary}[Strongly invariant natural extension]\label{cor:invnatex}
  Let\/ $\transfos$ be a monoid of transformations of\/ $\values$ and let\/ $\lpr$ be a weakly   $\transfos$-invariant lower prevision on some set of gambles $\domain$, that avoids sure loss. Assume that   $\transfos^*$ has the Moore--Smith property.  Then there are strongly $\transfos$-invariant coherent lower   previsions on $\gambles(\values)$ that dominate $\lpr$ on $\gambles(\values)$, and the smallest such lower   prevision, which is called the \emph{strongly $\transfos$-invariant natural extension} of $\lpr$, is given   by $\lnex_{\lpr,\transfos}=\alpr_{\lnex_\lpr,\transfos}$. Moreover, for every $\transfos$-invariant gamble   $f$ we have that $\lnex_{\lpr,\transfos}(f)=\lnex_\lpr(f)$.
\end{corollary}

\begin{proof}
  The first part of the proof follows at once from the observation that a coherent lower prevision $\alpr$ on   $\gambles(\values)$ dominates $\lpr$ on $\domain$ if and only if it dominates $\lnex_\lpr$ on all   gambles. For the second part of the proof, simply observe that if $f$ is a $\transfos$-invariant gamble,   then $T^*f=f$ and therefore $\lnex_\lpr(T^*f)=\lnex_\lpr(f)$ for all $T^*$ in $\transfos^*$.
\end{proof}

Let us show in particular how this result applies when we consider the monoid $\transfos_T$ generated by a single transformation $T$:

\begin{corollary}\label{cor:invnatex-cyclic}
  Let $T$ be a transformation of $\values$ and consider the Abelian monoid $\transfos_T=\set{T^n}{n\geq0}$.   Then for any weakly $T$-invariant lower prevision $\lpr$ on some set of gambles $\domain$ that avoids sure   loss, there are strongly $T$-invariant coherent (lower) previsions on $\gambles(\values)$ that dominate   $\lpr$, and the point-wise smallest such lower prevision $\lnex_{\lpr,T}$ is given by
  \begin{equation*}
    \lnex_{\lpr,T}(f)
    =\lim_{n\to\infty}\lnex_\lpr
    \left(\frac{1}{n}\sum_{k=0}^{n-1}(T^k)^tf\right)
    =\sup_{n\geq1}\lnex_\lpr
    \left(\frac{1}{n}\sum_{k=0}^{n-1}(T^k)^tf\right).
  \end{equation*}
\end{corollary}

\begin{proof}
  The existence of strongly $T$-invariant coherent (lower) previsions on $\gambles(\values)$ that dominate   $\lpr$ follows from Corollary~\ref{cor:invnatex}, and the fact that for any Abelian monoid $\transfos$,   $\transfos^*$ has the Moore--Smith property. It also follows from this corollary that for any gamble $f$ on   $\values$,
  \begin{equation*}
    \lnex_{\lpr,T}(f)
    =\sup\set{\lnex_\lpr(T^*f)}{T^*\in\transfos_T^*}
    \geq\sup_{n\geq1}\lnex_\lpr
    \left(\frac{1}{n}\sum_{k=0}^{n-1}(T^k)^tf\right).
  \end{equation*}
  To prove the converse inequality, fix any $T^*$ in $\transfos_T^*$ and any gamble $f$ on $\values$.  Then   there is some $N\geq1$ and non-negative $\lambda_0$, \dots, $\lambda_{N-1}$ that sum to one, such that   $T^*=\sum_{k=0}^{N-1}\lambda_k(T^k)^t$. Consider the element   $S_M^*=\frac{1}{M}\sum_{\ell=0}^{M-1}(T^{\ell})^t$ of $\transfos^*$, where $M$ is any natural number such   that $M\geq N$.  Observe that
  \begin{equation*}
    S_M^*T^*
    =\frac{1}{M}\sum_{\ell=0}^{M-1}(T^{\ell})^t
    \left(\sum_{k=0}^{N-1}\lambda_k(T^{k})^t\right)
    =\sum_{\ell=0}^{M-1}\sum_{k=0}^{N-1}\frac{\lambda_k}{M}(T^{k+\ell})^t
    =\sum_{m=0}^{M+N-2}\mu_m(T^m)^t,
  \end{equation*}
  where we let, for $0\leq m\leq M+N-2$,
  \begin{equation*}
    \mu_m
    :=\sum_{k=0}^{N-1}\sum_{\ell=0}^{M-1}\frac{\lambda_k}{M}\delta_{m,k+\ell}
    =\begin{cases}
      \sum_{k=0}^{m} \frac{\lambda_k}{M} &\text{ if } 0 \leq m
      \leq N-2 \\
      \frac{1}{M} &\text{ if } N-1 \leq m \leq M-1 \\
      \sum_{k=m-M+1}^{N-1} \frac{\lambda_k}{M} &\text{ if } M \leq m
      \leq M+N-2.
    \end{cases}
  \end{equation*}
  This tells us that $\mu_m=\frac{1}{M}$ for $N-1\leq m\leq M-1$, and $0\leq\mu_m\leq\frac{1}{M}$ for all   other $m$.  If we let $\delta_m:=\mu_m-\frac{1}{N+M-1}$, it follows at once that
  \begin{equation*}
    \abs{\delta_m}\leq
    \begin{cases}
      \dfrac{N-1}{M(M+N-1)}&\text{if $N-1\leq m\leq M-1$}\\
      \dfrac{1}{M+N-1}&\text{if $0\leq m\leq N-2$ or $M \leq m\leq M+N-2$}
    \end{cases}
  \end{equation*}
  Consequently, it follows from the weak $\transfos$-invariance and the coherence of $\lnex_\lpr$ that
  \begin{multline*}
    \lnex_{\lpr}(T^*f)\\
    \begin{aligned}
      &\leq\lnex_{\lpr}(S_M^*T^*f)\\
      &=\lnex_{\lpr}\left(S_{M+N-1}^*f +\sum_{m=0}^{M+N-2}\delta_m(T^m)^tf\right)       \leq\lnex_{\lpr}(S_{M+N-1}^*f)
      +\sum_{m=0}^{M+N-2}\abs{\delta_m}\sup\abs{f}\\
      &\leq\lnex_{\lpr}(S_{M+N-1}^*f)+\sup\abs{f}
      \left[\dfrac{N-1}{M(M+N-1)}(M-N+1) +\dfrac{1}{M+N-1}(2N-2)\right]\\
      &=\lnex_{\lpr}(S_{M+N-1}^*f)+\sup\abs{f} \dfrac{(N-1)(3M-N+1)}{M(M+N-1)}.
    \end{aligned}
  \end{multline*}
  Recall that $f$ and $T^*$, and therefore also $N$ are fixed. Consider any $\epsilon>0$, then there is some   $M_\epsilon\geq N$ such that $\sup\abs{f} \frac{(N-1)(3M-N+1)}{M(M+N-1)}<\epsilon$ for all $M\geq   M_\epsilon$, whence
  \begin{equation*}
    \lnex_{\lpr}(T^*f)
    \leq\lnex_{\lpr}(S_{M_\epsilon+N-1}^*f)+\epsilon
    \leq\sup_{n\geq 1}\lnex_{\lpr}(S_n^*f)+\epsilon.
  \end{equation*}
  Since this holds for all $\epsilon>0$, we get $\lnex_{\lpr}(T^*f)\leq\sup_{n\geq     1}\lnex_{\lpr}(S_n^*f)$. Taking the supremum over all $T^*$ in $\transfos^*$ leads to the desired   inequality.
\end{proof}

\subsection{Mixture lower previsions}\label{sec:mixture-lower-previsions}
The condition established in Theorem~\ref{theo:Moore-Smith} is fairly general, and guarantees for instance the existence of $\transfos$-invariant coherent previsions whenever the monoid $\transfos$ is Abelian, or a finite group. In case $\transfos^*$ is not directed, however, as may happen for instance for groups $\transfos$ that are not finite nor Abelian, there may still be $\transfos$-invariant coherent previsions, as we shall see in Example~\ref{ex:not-directed} below. So we see that the directedness of $\transfos^*$ is not a necessary condition for the existence of $\transfos$-invariant coherent previsions.
\par
But consider a weakly $\transfos$-invariant lower prevision $\lpr$ defined on some domain $\domain$, that avoids sure loss. Even if $\transfos^*$ is not directed,\footnote{This is the general situation that   \citet[Section~3.5]{walley1991} considers, and he doesn't discuss the directedness of $\transfos^*$. He does   consider the special case that $\transfos$ is Abelian for which he proves that the existence of invariant   coherent previsions is guaranteed. The results in this section were first proven by him.} we may still associate with $\lpr$ a lower prevision $\alpr_{\lpr,\transfos}$ on $\gambles(\values)$ through Eq.~\eqref{eq:Q-functional}:
\begin{equation*}
  \alpr_{\lpr,\transfos}(f)
  =\sup_{T^*\in\transfos^*}\lnex_\lpr(T^*f)
  =\sup\set{\lnex_\lpr\left(\frac{1}{n}\sum_{k=1}^{n}T_k^tf\right)}
  {n\geq1, T_1,\dots,T_n\in\transfos},
\end{equation*}
where we have replaced the Moore--Smith limit by a supremum (with which it would coincide in case $\transfos^*$ were directed), and where $\lnex_\lpr$ is the natural extension of $\lpr$ to all gambles. We shall call this lower prevision the \emph{mixture lower prevision} associated with the weakly invariant $\lpr$. The supremum in this expression is finite, since it is dominated by $\sup f$. This mixture lower prevision is not necessarily coherent, but it is still strongly $\transfos^*$-invariant.\footnote{Simply   observe that the relevant part (near the end) of the proof of Theorem~\ref{theo:Moore-Smith} is not based on   the directedness of $\transfos^*$.} Moreover, this mixture lower prevision dominates $\lnex_\lpr$, and therefore also $\lpr$ [observe that $\lnex_\lpr$ is weakly invariant because $\lpr$ is]; and if there are $\transfos$-invariant coherent previsions, it is dominated by the strongly $\transfos$-invariant natural extension $\lnex_{\lpr,\transfos}$ of $\lpr$.\footnote{To prove that the mixture lower prevision dominates   $\lpr$, consider $T^*=\iden_\values$ in its definition. To prove that it is dominated by the strongly   invariant natural extension, take $f_k=f/n$ in the expression~\eqref{eq:natex-with-transfos} for this   natural extension.} This shows that $\solp(\alpr_{\lpr,\transfos})=\solp(\lnex_{\lpr,\transfos})$, since all coherent previsions that dominate the strongly $\transfos$-invariant $\alpr_{\lpr,\transfos}$ are necessarily $\transfos$-invariant. And clearly then, if this mixture lower prevision is coherent, it coincides with the strongly invariant natural extension. So we see that the mixture lower prevision, even if it is not coherent, still allows us to characterise all $\transfos$-invariant coherent previsions. In particular, there are such invariant coherent previsions if and only if it avoids sure loss.

\begin{example}[Directedness is not necessary]\label{ex:not-directed}
  Let us consider the space $\values_3:=\{1,2,3\}$, and let $T_1$ and $T_2$ be the transformations of   $\values$ given by $T_1(1)=1$, $T_1(2)=2$, $T_1(3)=2$ and $T_2(1)=1$, $T_2(2)=3$, $T_2(3)=3$,   respectively. Since $T_1T_1=T_1$, $T_2T_2=T_2$, $T_2T_1=T_2$ and $T_1T_2=T_1$, we deduce that the set of   transformations $\transfos=\{\iden_\values,T_1,T_2\}$ is a monoid.  Let $P_{\{1\}}$ be the coherent   prevision on $\gambles(\values)$ given by $P_{\{1\}}(f)=f(1)$ for any gamble $f$, i.e., all of whose   probability mass lies in $1$.  Then we have $P_{\{1\}}(f)=P_{\{1\}}(T_1^tf)=P_{\{1\}}(T_2^tf)$ for any   gamble $f$, so $P_{\{1\}}$ is $\transfos$-invariant. Let us show that $\transfos^*$ does not have the   Moore--Smith property.
  \par
  Consider $T_1^*$ and $T_2^*$ in $\transfos^*$ given by $T_1^*=\lambda T_1^t+(1-\lambda)T_2^t$ and $T_2^*=\mu   T_1^t+(1-\mu)T_2^t$, with $\lambda\neq\mu$.  Let $T^*$ be another element of $\transfos^*$, so there are   non-negative $\alpha_1$, $\alpha_2$ and $\alpha_3$ such that $\alpha_1+\alpha_2+\alpha_3=1$ and   $T^*=\alpha_1\iden_\values^t+\alpha_2 T_1^t+\alpha_3T_2^t$. Now,
  \begin{align*}
    T^* T_1^* &=\alpha_1\lambda\iden_\values^tT_1^t
    +\alpha_1(1-\lambda)\iden_\values^tT_2^t\\
    &\qquad+\alpha_2\lambda T_1^tT_1^t+\alpha_2(1-\lambda)T_1^tT_2^t
    +\alpha_3\lambda T_2^tT_1^t+\alpha_3(1-\lambda)T_2^tT_2^t\\
    &=\alpha_1\lambda T_1^t+\alpha_1(1-\lambda)T_2^t+\alpha_2\lambda T_1^t
    +\alpha_2(1-\lambda)T_2^t+\alpha_3\lambda T_1^t+\alpha_3(1-\lambda)T_2^t\\
    &=\lambda T_1^t+(1-\lambda)T_2^t=T_1^*.
  \end{align*}
  Similarly, $T^* T_2^*=T_2^*$ for any $T^* \in \transfos^*$. This means that $T_1^*$ is the only possible   successor of $T_1^*$, and $T_2^*$ is the only possible successor of $T_2^*$. Hence, $\transfos^*$ cannot   have the Moore--Smith property. Nevertheless, there is a $\transfos$-invariant coherent prevision   $\pr_{\{1\}}$.
  \par
  Let us consider the vacuous, and therefore weakly $\transfos$-invariant and coherent, lower prevision   $\lpr_{\values_3}$ on $\gambles(\values_3)$, and the mixture lower prevision   $\alpr_{\lpr_{\values_3},\transfos}$ that corresponds with it. It is easy to show that for any gamble $f$,   $\alpr_{\lpr_{\values_3},\transfos}(f)=\min\{f(1),\max\{f(2),f(3)\}\}$ and this lower prevision avoids sure   loss, and is therefore strongly $\transfos$-invariant, but it is not coherent [it is not super-additive]. It   is easy to see that $\pr_{\{1\}}$ is the only coherent prevision that dominates   $\alpr_{\lpr_{\values_3},\transfos}$, and is therefore the only $\transfos$-invariant coherent   prevision. $\blacklozenge$
\end{example}

\subsection{Invariance and Choquet integration}
Until now, we have explored the relation between coherence and (weak or strong) invariance. To complete this section, we intend to explore this relation for the particular case of the $n$-monotone lower previsions and probabilities introduced near the end of Section~\ref{sec:linear-previsions}.
\par
Consider an $n$-monotone lower probability $\lpr$ defined on a lattice of events $\domain$ containing $\emptyset$ and $\values$. Then its natural extension to all events coincides with its inner set function $\lpr_*$, which is given by $\lpr_*(A)=\sup\set{\lpr(B)}{B\in\domain, B\subseteq A}$.  Furthermore, the natural extension to all gambles is given by the Choquet integral with respect to $\lpr_*$:
\begin{equation*}
  \lnex_\lpr(f)
  =(C)\int_\values f\dif\lpr_*
  :=\inf f+(R)\int_{\inf f}^{\sup f}
  \lpr_*(\set{x\in\values}{f(x)\geq\alpha})\dif\alpha
\end{equation*}
for all gambles $f$ on $\values$, where the integral on the right-hand side is a Riemann integral. This natural extension (and therefore also the inner set function) is still $n$-monotone \Citep{cooman2005d,cooman2005b}. Since we have proven in Theorem~\ref{theo:natext-preserves-invariance} that natural extension preserves weak invariance, we can deduce that the inner set function of a $n$-monotone weakly invariant coherent lower probability, and the associated Choquet functional, are still weakly invariant, $n$-monotone and coherent. We now show that weak invariance of the inner set function and the associated Choquet integral is still guaranteed if the lower probability $\lpr$ is not coherent or 2-monotone, but only monotone. In what follows, it is important to remember that for a transformation $T$ of $\values$ and a subset $A$ of $\values$, $T^tI_A=I_{T^{-1}(A)}$.

\begin{proposition}\label{prop:Choquet-weak}
  Let $\lpr$ be a weakly $\transfos$-invariant monotone lower probability, defined on a $\transfos$-invariant   lattice of events $\domain$ that contains $\emptyset$ and $\values$, and such that $\lpr(\emptyset)=0$ and   $\lpr(\values)=1$. Then
  \begin{enumerate}[1.]
  \item the inner set function $\lpr_*$ of $\lpr$ is weakly $\transfos$-invariant; and
  \item the Choquet integral with respect to $\lpr_*$ is weakly $\transfos$-invariant.
  \end{enumerate}
\end{proposition}

\begin{proof}
  To prove the first statement, consider any $A\subseteq\values$, and let $B\in\domain$ be a any subset of   $A$. Then for any $T$ in $\transfos$, $T^{-1}(B)\in\domain$ and $T^{-1}(B)=\set{x}{Tx\in     B}\subseteq\set{x}{Tx\in A}=T^{-1}(A)$, whence $\lpr(B)\leq\lpr(T^{-1}(B))\leq\lpr_*(T^{-1}(A))$, where   the first inequality follows from the weak invariance of $\lpr$, and the second from the fact that $\lpr_*$   is monotone and coincides with $\lpr$ on its domain, because $\lpr$ is assumed to be monotone. Consequently   $\lpr_*(A)=\sup_{B\in\domain,B\subseteq A}\lpr(B)\leq\lpr_*(T^{-1}(A))$.  Hence, $\lpr_*$ is also weakly   $\transfos$-invariant.
  \par
  To prove the second statement, let $f$ be any gamble on $\values$. Define, for any $\alpha$ in $\reals$, the   level set $f_\alpha:=\set{x}{f(x)\geq\alpha}$. Then by the first statement,
  \begin{equation*}
    \lpr_*(f_\alpha)
    \leq\lpr_*(T^{-1}(f_\alpha))
    =\lpr_*(\set{x}{Tx\in f_\alpha})
    =\lpr_*(\set{x}{f(Tx)\geq\alpha})
    =\lpr_*((T^tf)_\alpha).
  \end{equation*}
  Hence,
  \begin{multline*}
    (C)\int f\dif P_*
    =\inf f+(R)\int_{\inf f}^{\sup f}\lpr_*(f_\alpha)\dif\alpha\\
    \leq\inf f+(R)\int_{\inf f}^{\sup f}\lpr_*((T^tf)_\alpha)\dif\alpha =(C)\int T^tf\dif\lpr_*,
  \end{multline*}
  also taking into account for the last equality that $\lpr_*((T^tf)_\alpha)=1$ for all $\alpha$ in $[\inf   f,\inf T^tf)$, and that $\lpr_*(f_\alpha)=0$ for all $\alpha$ in $(\sup T^tf,\sup f]$.
\end{proof}

As we said before, natural extension does not preserve strong invariance in general, and a simple example shows that this continues to hold in particular for $n$-monotone lower previsions: the unique coherent lower prevision defined on $\{\emptyset,\values\}$ is trivially completely monotone and strongly invariant with respect to any monoid of transformations $\transfos$, but its natural extension, the vacuous lower prevision $\lpr_\values$ (which is completely monotone), is not strongly $\transfos$-invariant unless in the trivial case that $\transfos=\{\iden_\values\}$.
\par
It is nonetheless interesting that if we restrict ourselves to coherent previsions (which constitute a particular instance of completely monotone lower previsions), natural extension from events to gambles does preserve strong invariance. This is a consequence of the following theorem.

\begin{theorem}\label{theo:natural-strong}
  Let $\lpr$ be a coherent lower prevision on $\gambles(\values)$ and let $\transfos$ be a monoid of   transformations on $\values$. Then $\lpr$ is strongly $\transfos$-invariant if and only if any $\pr$ in   $\solp(\lpr)$, its restriction to events is (weakly) $\transfos$-invariant, in the sense that   $\pr(T^{-1}(A))=\pr(A)$ for all $A\subseteq\values$ and all $T\in\transfos$.
\end{theorem}

\begin{proof}
  We start with the direct implication. If $\lpr$ is strongly $\transfos$-invariant, then any $\pr$ in   $\solp(\lpr)$ is $\transfos$-invariant by Theorem~\ref{theo:dominating-invariant-previsions}.  Hence, given   $A\subseteq\values$ and $T\in\transfos$, we get $\pr(A)=\pr(T^{-1}(A))$.
  \par
  Conversely, consider $\pr$ in $\solp(\lpr)$. Recall that a coherent prevision on all events has only one   coherent extension from all events to all gambles, namely its natural extension, or Choquet functional; see   \Citep{cooman2005b}. So for any gamble $f$ on $\values$ and any $T$ in $\transfos$, taking into account that   $\pr$ is assumed to be invariant on events, and that $T^{-1}(f_\alpha)=(T^tf)_\alpha$ [see the proof of   Proposition~\ref{prop:Choquet-weak}], we get
  \begin{align*}
    \pr(f) &=(C)\int_\values f\dif\pr
    =\inf f+(R)\int_{\inf f}^{\sup f}\pr(f_\alpha)\dif\alpha\\
    &=\inf f+(R)\int_{\inf f}^{\sup f}\pr(T^{-1}(f_\alpha))\dif\alpha
    =\inf f+(R)\int_{\inf f}^{\sup f}\pr((T^tf)_\alpha)\dif\alpha\\
    &=(C)\int_\values T^tf\dif\pr =\pr(T^tf).
  \end{align*}
  Hence, $\pr$ is strongly $\transfos$-invariant and, applying   Theorem~\ref{theo:dominating-invariant-previsions}, so is the lower envelope $\lpr$ of $\solp(\lpr)$.
\end{proof}
\noindent
We see that, although the condition of strong invariance cannot be considered for lower probabilities, in the sense that $I_A-T^tI_A$ will not be in general the indicator of an event, it is still to some extent characterised by behaviour on events. Moreover, we may deduce the following result.

\begin{corollary}\label{co:nex-strong}
  Let $\lpr$ be a strongly $\transfos$-invariant lower prevision on a $\transfos$-invariant set of gambles   $\domain$ that includes all indicators of events. Assume that $\lpr$ avoids sure loss. Then its natural   extension to all gambles is strongly $\transfos$-invariant, and coincides therefore with the strongly   invariant natural extension of $\lpr$.
\end{corollary}

\begin{proof}
  Since $\lpr$ avoids sure loss, $\solp(\lpr)$ is non-empty. Since $\lpr$ is strongly invariant on a domain   that includes all events, any element $\pr$ of $\solp(\lpr)$ is (strongly) invariant on all events. Hence,   by the previous theorem, $\pr$ is also (strongly) invariant on all gambles, since a coherent prevision on   all events has only one coherent extension from all events to all gambles (namely its natural extension, or   Choquet functional).  Therefore, the natural extension of $\lpr$ is a lower envelope of invariant coherent   previsions, and is therefore strongly invariant.
\end{proof}
\noindent This result provides further insight into the existence problem for strongly invariant coherent lower previsions. The existence of strongly invariant coherent lower previsions on all gambles is equivalent to the existence of invariant coherent previsions on all gambles, which in turn is equivalent to the existence of invariant coherent previsions on all events (or in other words, invariant finitely additive probabilities). And it is the impossibility of satisfying invariance with finitely additive probabilities in some cases (for instance for the class $\transfos_\values$ of all transformations) that prevents the existence of coherent strongly invariant belief models.
\par
We also infer that if the restriction $\alpr$ of a coherent lower prevision $\lpr$ on $\gambles(\values)$ to gambles of the type $I_A-T^tI_A$ and $T^tI_A-I_A$, involving only indicators of events, is strongly invariant, then $\lpr$ is strongly invariant on all of $\gambles(\values)$: it will dominate the natural extension $\lnex_\alpr$ of $\alpr$, which is strongly invariant by Corollary~\ref{co:nex-strong}, and consequently it will also be strongly invariant.
\par
We can also deduce the following result. Recall that a linear lattice of gambles $\domain$ is a set of gambles that is at once a lattice of gambles and a linear subspace of $\gambles(\values)$. If in addition $\domain$ contains all constant gambles, then for any coherent prevision $\pr$ defined on $\domain$, its natural extension to all gambles \citep[Theorem~3.1.4]{walley1991} is given by the \emph{inner extension} $\pr_*(f):=\sup\set{\pr(g)}{g\in\domain,g\leq f}$.  Let us denote by $\pr^*$ the conjugate upper prevision of $\pr_*$.

\begin{corollary}
  Let $\transfos$ be a monoid of transformations of $\values$, and let $\lpr$ be a strongly   $\transfos$-invariant lower prevision on a linear lattice of gambles $\domain$ that contains all constant   gambles. The natural extension $\lnex_\lpr$ of $\lpr$ to all gambles is strongly $\transfos$-invariant if   and only if for any coherent prevision $\pr$ on $\domain$ that dominates $\lpr$, we have $\pr_*(A\setminus   T^{-1}(A))=\pr^*(A\setminus T^{-1}(A))=\pr_*(T^{-1}(A)\setminus A)=\pr^*(T^{-1}(A)\setminus A)$ for all   $A\subseteq\values$ and all $T\in\transfos$.
\end{corollary}

\begin{proof}
  It follows from \citet[Theorem~3.4.2]{walley1991} that $\lnex_\lpr$ is the lower envelope of the coherent   lower previsions $\pr_*$, where $\pr$ is any coherent prevision on $\domain$ that dominates $\lpr$ on   $\domain$. But then, clearly, $\lnex_\lpr$ will be strongly $\transfos$-invariant if and only if all the   $\pr_*$ are. Consider any such $\pr_*$. By Theorem~\ref{theo:natural-strong}, $\pr_*$ is strongly invariant   if and only if for all $A\subseteq\values$ and $T\in\transfos$:
  \begin{equation*}
    \apr(A)=\apr(T^{-1}(A))\quad\text{for all $\apr$ in $\solp(\pr_*)$}
  \end{equation*}
  which is obviously equivalent to $\pr_*(I_A-T^tI_A)=\pr_*(T^tI_A-I_A)=0$.  Now observe that   $I_A-T^tI_A=I_A-I_{T^{-1}(A)}=I_{A\setminus T^{-1}(A)}-I_{T^{-1}(A)\setminus A}$, and that the functions   $I_{A\setminus T^{-1}(A)}$ and $-I_{T^{-1}(A)\setminus A}$ are comonotone.  Since $\pr$ is a coherent   prevision on $\domain$, it is completely monotone.  Hence, its inner extension $\pr_*$ is coherent and   completely monotone on all gambles, and therefore comonotone additive \Citep{cooman2005b}. This means that
  \begin{multline*}
    \pr_*(I_A-T^tI_A) =\pr_*(I_{A\setminus T^{-1}(A)}-I_{T^{-1}(A)\setminus A})
    =\pr_*(I_{A\setminus T^{-1}(A)})+\pr_*(-I_{T^{-1}(A)\setminus A})\\
    =\pr_*(I_{A\setminus T^{-1}(A)})-\pr^*(I_{T^{-1}(A)\setminus A}) =\pr_*(A\setminus     T^{-1}(A))-\pr^*(T^{-1}(A)\setminus A)
  \end{multline*}
  and similarly $\pr_*(T^tI_A-I_A)=\pr_*(T^{-1}(A)\setminus A)-\pr^*(A\setminus T^{-1}(A))$. The rest of the   proof is now immediate.
\end{proof}

\section{Shift-invariance and its generalisations}
\label{sec:shift-invariance}
\subsection{Strongly shift-invariant coherent lower previsions on $\gambles(\nats)$} Let us consider, as an example, the case of the shift-invariant, i.e., $\transfos_\theta$-invariant, coherent previsions on $\gambles(\nats)$. These are usually called \emph{Banach limits} in the literature, see for instance, \citet[Section~2.1.3]{bhaskara1983} or \citet[Sections~2.9.5 and~3.5.7]{walley1991}. We know from Corollary~\ref{cor:invnatex} that there are always Banach limits that dominate a given weakly shift-invariant lower prevision---so we know that there actually are Banach limits. Let us denote by $\linprevs_\theta(\nats)$ the set of all Banach limits. We also know that a coherent lower prevision on $\gambles(\nats)$ is strongly shift-invariant if and only if it is a lower envelope of such Banach limits.  The smallest strongly shift-invariant coherent lower prevision $\lnex_\theta$ on $\gambles(\nats)$ is the lower envelope of all Banach limits, and it is given by:\footnote{See also \citet[Section~3.5.7]{walley1991}. The expression on the   right hand side is not a limit inferior!}
\begin{equation}\label{eq:shift-invariant-natex}
  \lnex_\theta(f)
  =\sup_{\substack{m_1,\dots m_n\geq0\\n\geq0}}\,
  \inf_{k\geq0}\,\frac{1}{n}\sum_{\ell=1}^nf(k+m_\ell)
  =\lim_{n\to\infty}\,\inf_{k\geq0}\,\frac{1}{n}\sum_{\ell=k}^{k+n-1}f(\ell),
\end{equation}
for any gamble $f$ on $\nats$ (or in other words, for any bounded sequence $f(n)_{n\in\nats}$ of real numbers). The first equality follows from Corollary~\ref{cor:invnatex}, and the second from Corollary~\ref{cor:invnatex-cyclic}. $\lnex_\theta(f)$ is obtained by taking the infimum sample mean of $f$ over `moving windows' of length $n$, and then letting the window length $n$ go to infinity. Since this is the lower prevision on $\gambles(\nats)$ that can be derived \emph{solely} using considerations of coherence and the evidence of shift-invariance, we believe that this $\lnex_\theta$ is a natural candidate for a `\emph{uniform distribution}' on $\nats$. It is the belief model to use if we only have evidence of shift-invariance, as all other strongly shift-invariant coherent lower previsions will point-wise dominate $\lnex_\theta$, and will therefore represent stronger behavioural dispositions than warranted by the mere evidence of shift-invariance.\footnote{But this belief model has the important defect that, like the lower   prevision $\lsamp_\theta$ defined further on, it is not fully conglomerable; see   \citet[Section~6.6.7]{walley1991} and observe that the counterexample that Walley gives for $\lsamp_\theta$,   also applies to $\lnex_\theta$. Walley's remark there that his example shows that there are no (what we   call) fully conglomerable (strongly) shift-invariant (lower) previsions that dominate $\lsamp_\theta$, can   be extended in a straightforward manner to $\lnex_\theta$ to show that \emph{there are no fully     conglomerable (strongly) shift-invariant (lower) previsions}.\label{fn:congomerability}}
\par
We could also sample $f$ over the set $\{1,\dots,n\}$ leading to a coherent `sampling' prevision
\begin{equation*}
  S_n(f)=\frac{1}{n}\sum_{\ell=0}^{n-1}f(\ell),
\end{equation*}
but the problem here is that for any given $f$ the sequence of sampling averages $\samp_n(f)$ is not guaranteed to converge. Taking the limits inferior of such sequences (one for each gamble $f$), however, yields a coherent lower prevision\footnote{A limit inferior of a sequence of coherent lower previsions is   always coherent, see \citet[Corollary~2.6.7]{walley1991}.}  $\lsamp_\theta$ given by
\begin{equation*}
  \lsamp_\theta(f)=\liminf_{n\to\infty}\samp_n(f)
  =\liminf_{n\to\infty}\frac{1}{n}\sum_{\ell=0}^{n-1}f(\ell)
\end{equation*}
for any gamble $f$ on $\nats$. For any event $A\subseteq N$, or equivalently, any zero-one-valued sequence, we have that $\samp_n(A)=\frac{1}{n}\abs{A\cap\{0,\dots,n-1\}}$ is the `relative frequency' of ones in the sequence $I_A(n)$ and
\begin{equation*}
  \lsamp_\theta(A)=\liminf_{n\to\infty}\samp_n(A)
  =\liminf_{n\to\infty}\frac{1}{n}\abs{A\cap\{0,\dots,n-1\}}.
\end{equation*}
Let $\usamp_\theta$ denote the conjugate of $\lsamp_\theta$, given by $\usamp_\theta(f)=\limsup_n S_n(f)$. Those events $A$ for which $\lsamp_\theta(A)=\usamp_\theta(A)$ have a `limiting relative frequency' equal to this common value. It is not difficult to show that the coherent `limiting relative frequency' lower prevision $\lsamp_\theta$ is actually also strongly shift-invariant.\footnote{The following simple proof is   due to \citet[Section~3.5.7]{walley1991}. Observe that $\samp_n(\theta^tf-f)=[f(n)-f(0)]/n\to0$ as   $n\to\infty$, so $\lsamp_\theta(\theta^tf-f)=\usamp_\theta(\theta^tf-f)=0$.} This implies that all the coherent previsions that dominate $\lsamp_\theta$ are strongly shift-invariant. But it is easy to see (see Example~\ref{ex:banach-limits} below) that $\lnex_\theta$ is strictly dominated by $\lsamp_\theta$, so there are Banach limits that do not dominate $\lsamp_\theta$.

\begin{proposition}\label{prop:banach-limits}
  Let $\banlim$ be any Banach limit on $\gambles(\nats)$, let $f$ be any gamble on $\nats$. Then the following   statements hold.
  \begin{enumerate}[1.]
  \item $\liminf_{n\to\infty}f(n)\leq\lnex_\theta(f)\leq\lsamp_\theta(f)\leq     \usamp_\theta(f)\leq\unex_\theta(f)\leq\limsup_{n\to\infty}f(n)$.
  \item If\/ $\lim_{n\to\infty}f(n)$ exists, then
    \begin{equation*}
      \lnex_\theta(f)=\lsamp_\theta(f)
      =\unex_\theta(f)=\usamp_\theta(f)
      =\banlim(f)=\lim_{n\to\infty}f(n).
    \end{equation*}
  \item If $f$ is $\theta^m$-invariant (has period $m\geq1$), then
    \begin{equation*}
      \lnex_\theta(f)=\lsamp_\theta(f)
      =\unex_\theta(f)=\usamp_\theta(f)
      =\banlim(f)=\frac{1}{m}\sum_{r=1}^{m-1}f(r).
    \end{equation*}
  \item If $f$ is zero except in a finite number of elements of\/ $\nats$, then     $\lnex_\theta(f)=\lsamp_\theta(f)=\unex_\theta(f)=\usamp_\theta(f) =L(f)=0$. In particular, this holds for     the indicator of any finite subset $A$ of\/ $\nats$.
  \end{enumerate}
\end{proposition}

\begin{proof}
  We begin with the first statement. By conjugacy, we can concentrate on the lower previsions. We have already   argued that $\lsamp_\theta$ is a strongly shift-invariant coherent lower prevision, so $\lsamp_\theta$ will   dominate the smallest strongly shift-invariant coherent lower prevision $\lnex_\theta$. So it remains to   prove that $\lnex_\theta$ dominates the limit inferior. Consider the first equality in   Eq.~\eqref{eq:shift-invariant-natex}. Fix the natural numbers $n\geq1$, $m_1$, \dots $m_n$. We can assume   without loss of generality that the $m_1$ is the smallest of all the $m_\ell$. Observe that
  \begin{equation*}
    \inf_{k\geq0}\,\frac{1}{n}\sum_{\ell=1}^nf(k+m_\ell)
    \geq\inf_{k\geq0}\,\min_{\ell=1}^nf(k+m_\ell)
    =\min_{\ell=1}^n\,\inf_{k\geq m_\ell}f(k)
    =\inf_{k\geq m_1}f(k),
  \end{equation*}
  and therefore
  \begin{equation*}
    \lnex_\theta(f)
    \geq\sup_{m_1\geq0}\,\inf_{k\geq m_1}f(k)
    =\liminf_{n\to\infty}f(n).
  \end{equation*}
  \par
  The second statement is an immediate consequence of the first, and the third follows easily from the   definition of $\lnex_\theta$ and $\unex_\theta$.  Finally, the fourth statement follows at once from the   second.
\end{proof}

\begin{example}[Not all Banach limits dominate $\lsamp_\theta$]\label{ex:banach-limits} Consider the event
  \begin{equation*}
    A=\set{n^2+k}{n\geq1,k=0,\dots,n-1}.
  \end{equation*}
  Then $A$ has `limiting relative frequency' $\lsamp_\theta(A)=\usamp_\theta(A)=1/2$, whereas   $\lnex_\theta(A)=0$ and $\unex_\theta(A)=1$. This shows that $\lsamp_\theta$ strictly dominates   $\lnex_\theta$, so not all Banach limits dominate $\lsamp_\theta$.
  \par
  Indeed, for the limiting relative frequency, consider the subsequence $S_{m^2-1}(A)$, $m\geq2$ of   $S_n(A)$. Then
  \begin{equation*}
    S_{m^2-1}(A)
    =\frac{1}{m^2-1}\abs{A\cap\{0,\dots,m^2-2\}}
    =\frac{1+2+\dots+m-1}{m^2-1}
    =\frac{\frac{1}{2}m(m-1)}{m^2-1}
    =\frac{1}{2}\frac{m}{m+1},
  \end{equation*}
  so this subsequence converges to $\frac{1}{2}$. Now the `integer intervals' $[m^2-1,(m+1)^2-1]$, $m\geq1$   cover the set of all natural numbers, and as $n$ varies over such an interval, $S_n(A)$ starts at   $S_{m^2-1}(A)=\frac{1}{2}\frac{m}{m+1}<\frac{1}{2}$, increases to   $S_{m^2+m}(A)=\frac{1}{2}\frac{m^2+m}{m^2+m}=\frac{1}{2}$, and then again decreases to   $S_{(m+1)^2-1}(A)=\frac{1}{2}\frac{m+1}{m+2}<\frac{1}{2}$. Both the lower and upper bounds converge to   $\frac{1}{2}$ as $m\to\infty$, and therefore the sequence $S_n(A)$ converges to $\frac{1}{2}$ as well.
  \par
  To calculate $\lnex_\theta(A)$, we consider the second equality in Eq.~\eqref{eq:shift-invariant-natex}. Fix   $n\geq1$ and let $k=n^2+n$, then $k+n-1=(n+1)^2-2$, so
  \begin{equation*}
    \frac{1}{n}\sum_{\ell=k}^{k+n-1}I_A(\ell)
    =\frac{1}{n}\sum_{\ell=n^2+n}^{(n+1)^2-2}I_A(\ell)
    =0,
  \end{equation*}
  whence $\inf_{k\geq0}\frac{1}{n}\sum_{\ell=k}^{k+n-1}I_A(\ell)=0$ for all $n\geq1$, and therefore   $\lnex_\theta(A)=0$. To calculate $\unex_\theta(A)$, fix $n\geq1$ and let $k=n^2$ then
  \begin{equation*}
    \frac{1}{n}\sum_{\ell=k}^{k+n-1}I_A(\ell)
    =\frac{1}{n}\sum_{\ell=n^2}^{n^2+n-1}I_A(\ell)=1,
  \end{equation*}
  whence $\sup_{k\geq0}\frac{1}{n}\sum_{\ell=k}^{k+n-1}I_A(\ell)=1$ for all $n\geq1$, and therefore   $\unex_\theta(A)=1$. $\blacklozenge$
\end{example}

In an interesting paper, \citet{kadane1995} study candidates for the `uniform distribution' on $\nats$. They consider, among others, all the finitely additive probabilities (or equivalently, all coherent previsions) that coincide with the limiting relative frequency on all events for which this limit exists. One could also consider as such candidates the coherent previsions that dominate the sampling lower prevision $\lsamp_\theta$, which have the benefit of being strongly shift-invariant. But, we actually believe that \emph{all} Banach limits (or actually, their lower envelope) are good candidates for being called `uniform distributions on $\nats$' and not just the ones that dominate $\lsamp_\theta$. \citeauthor*{kadane1995} also propose to consider other coherent previsions, and their idea is to consider the `residue sets', which are the subsets
\begin{equation*}
  R_m^r=\set{km+r}{k\geq0}=\set{\ell\in\nats}{\ell=r\mod m}
\end{equation*}
of $\nats$, where $m\geq 1$ and $r=1,\dots,m-1$. These sets are $\theta^m$-invariant, so we already know from Proposition~\ref{prop:banach-limits} that $\lnex_\theta(R_m^r)=\lsamp_\theta(R_m^r)=\usamp_\theta(R_m^r) =\unex_\theta(R_m^r)=\frac{1}{m}$ for all $m\geq 1$ and $r=1,\dots,m-1$. Now what \citeauthor*{kadane1995} do, is consider the set of all coherent previsions (finitely additive probabilities in their paper, but that is equivalent) that extend the probability assessments $\pr(R_m^r)=1/m$ for all events $R_m^r$. In other words, they consider the natural extension $\lnex_{\mathrm{res}}$ of all such assessments, i.e., the lower envelope of all such coherent previsions. It is not difficult to prove that this natural extension is given by\footnote{See \Citet{cooman2005e} for a proof.}
\begin{equation*}
  \lnex_{\mathrm{res}}(f)
  =\lim_{m\to\infty}\frac{1}{m}\sum_{r=0}^{m-1}\inf_{k\in\nats}f(km+r).
\end{equation*}
This coherent lower prevision is completely monotone [as a point-wise limit of completely monotone lower previsions, even (natural extensions to gambles of so-called) belief functions \citep{shafer1976}], and weakly shift-invariant [since the natural extension of any weakly shift-invariant lower prevision is]. Since the assessments $\pr(R_m^r)=\frac{1}{m}$ coincide with the values given by $\lnex_\theta$, we see that $\lnex_\theta$ will point-wise dominate the natural extension $\lnex_{\mathrm{res}}$ of these assessments to all gambles.  But as we shall shortly prove in Example~\ref{ex:residue-sets}, $\lnex_{\mathrm{res}}$ is not strongly shift-invariant, meaning that among the coherent previsions that extend these assessments, there also are coherent previsions that are not Banach limits (not shift-invariant).

\begin{example}\label{ex:residue-sets}
  Here we show by means of a counterexample that $\lnex_{\mathrm{res}}$ is not strongly shift-invariant. Let   $B_m:=\{0,\dots,m-1\}$ and $A:=\bigcup_{m\geq1}\{m\}\times B_m$, and consider the map
  \begin{equation*}
    \phi\colon A\to\nats\colon(m,r)\mapsto\phi(m,r):=\frac{m(m-1)}{2}+r+1.
  \end{equation*}
  It is easy to see that $\phi$ is a bijection (one-to-one and onto). Also define the map
  \begin{equation*}
    \kappa\colon A\to\nats\colon(m,r)\mapsto\kappa(m,r):=Nm\phi(m,r)+r.
  \end{equation*}
  for some fixed $N\geq2$. We consider the strict order $<$ on $A$ induced by the bijection $\phi$, i.e.,   $(m,r)<(m',r')$ if and only if $\phi(m,r)<\phi(m',r')$ [if and only if $m<m$, or $m=m'$ and $r<r'$, so $<$   is the lexicographic order].  Then $\kappa$ is an increasing map with respect to this order. To see this,   assume that $(m,r)<(m',r')$.  If $m<m'$, then
  \begin{multline*}
    \kappa(m,r) =Nm\phi(m,r)+r
    <Nm\phi(m',0)+r\\
    <Nm'\phi(m',0)+0 \leq Nm'\phi(m',r')+r' =\kappa(m,r').
  \end{multline*}
  If on the other hand $m=m'$ and $r<r'$, then $\kappa(m,r)=Nm\phi(m,r)+r<Nm\phi(m,r')+r'=\kappa(m,r')$.
  \par
  Moreover, given $(m,r)<(m',r')$, we see that $\kappa(m',r')-\kappa(m,r)\geq N$.  Indeed, since $\kappa$ is   increasing, it suffices to prove this for consecutive pairs in the order $<$ we have defined on $A$. There   are only two possible expressions of consecutive pairs $(m,r)$ and $(m',r')$: either we have   $(m',r')=(m,r+1)$, and then we get
  \begin{equation*}
    \kappa(m,r+1)-\kappa(m,r)=Nm[\phi(m,r+1)-\phi(m,r)]+1=Nm+1\geq N;
  \end{equation*}
  or we have $r=m-1,(m',r')=(m+1,0)$, and then we get
  \begin{align*}
    \kappa(m+1,0)-\kappa(m,m-1)
    &=Nm[\phi(m+1,0)-\phi(m,m-1)]+N\phi(m+1,0)-(m-1)\\
    &=Nm+N\phi(m+1,0)-(m-1)\geq Nm\geq N,
  \end{align*}
  taking into account that $\phi(m+1,0)\geq m-1$ by definition of $\phi$.
  \par
  Consider the set $C=\kappa(A)^c$. Then $\lnex_{\mathrm{res}}(C)   =\lim_{m\to\infty}\frac{1}{m}\sum_{r=0}^{m-1}\inf_{k\in\nats}I_C(km+r)$.  Since for every $m\in\nats$ and   $r\in B_m$ the value $\kappa(m,r)=Nm\phi(m,r)+r$ does not belong to $C$, we deduce that   $\frac{1}{m}\sum_{r=0}^{m-1}\inf_{k\in\nats}I_C(km+r)=0$ for all $m$, and consequently   $\lnex_{\mathrm{res}}(C)=0$.
  \par
  On the other hand, $\lnex_\theta(C)   =\lim_{n\to\infty}\inf_{k\geq0}\frac{1}{n}\sum_{\ell=k}^{k+n-1}I_C(\ell)$.  Since by construction any two   elements in $\kappa(A)$ differ in at least $N$ elements, we deduce that   $\inf_{k\geq0}\frac{1}{n}\sum_{\ell=k}^{k+n-1}I_C(\ell)\geq 1-\frac{2}{N+1}$, and this for all   $n\in\nats$. This implies that $\lnex_\theta(C)\geq1-\frac{2}{N+1}>0$. Hence, $\lnex_{\mathrm{res}}$ is   strictly smaller than the smallest strongly shift-invariant natural extension $\lnex_\theta$, and therefore   not strongly shift-invariant.  $\blacklozenge$
\end{example}

\subsection{Strong $T$-invariance}
Now consider an arbitrary non-empty set $\values$. Also consider a transformation $T$ of $\values$ and the Abelian monoid $\transfos_T=\set{T^n}{n\geq0}$ generated by $T$. We shall characterise the strongly $T$-invariant coherent lower previsions on $\gambles(\values)$ using the Banach limits on $\gambles(\nats)$.
\par
First of all, consider any coherent lower prevision $\lpr$ on $\gambles(\values)$, and any gamble $f$ on $\values$. Define the gamble $f_\lpr$ on $\nats$ as
\begin{equation}\label{eq:bhaskara-gamble}
  f_\lpr(n):=\lpr((T^t)^nf)=\lpr(f\circ T^n).
\end{equation}
[This is indeed a gamble, as for all $n$ we deduce from the coherence of
$\lpr$ that $f_\lpr(n)=\lpr(f\circ T^n)\leq\sup[f\circ T^n]\leq\sup f$ and
similarly $f_\lpr(n)\geq\inf f$.] On the one hand
$(T^tf)_\lpr(n)=\lpr(T^tf\circ T^n))=\lpr(T^t(f\circ T^n))=f_{T\lpr}(n)$ and
on the other hand $(T^tf)_\lpr(n)=\lpr(f\circ
T^{n+1})=f_\lpr(n+1)=f_\lpr(\theta n)$, so
\begin{equation}\label{eq:bhaskara-lpr}
  (T^tf)_\lpr=f_{T\lpr}=\theta^t f_\lpr,
\end{equation}
and this observation allows us to establish a link between the transformation $T$ on $\values$ and the shift transformation $\theta$ on $\nats$. This makes us think of the following trick, inspired by what \citet[Section~2.1.3(9)]{bhaskara1983} do for probability charges, rather than coherent lower previsions.  Let $\banlim$ be any shift-invariant coherent prevision on $\gambles(\nats)$, or in other words, a Banach limit on $\gambles(\nats)$.  Define the real-valued functional $\lpr_\banlim$ on $\gambles(\values)$ by $\lpr_\banlim(f):=\banlim(f_\lpr)$. We show that this functional has very special properties.

\begin{proposition}\label{prop:bhaskara-banlim}
  Let $\banlim$ be a shift-invariant coherent prevision on $\gambles(\nats)$, let $\lpr$ be a coherent lower   prevision on $\gambles(\values)$, and let $T$ be a transformation of $\values$. Then the following   statements hold.
  \begin{enumerate}[1.]
  \item $\lpr_\banlim$ is a weakly $T$-invariant coherent lower prevision on $\gambles(\values)$ (with     equality).
  \item If $\lpr$ dominates a weakly $T$-invariant coherent lower prevision $\alpr$ on $\gambles(\values)$,     then $\lpr_\banlim$ dominates $\alpr$.
  \item If $\lpr=\pr$ is a coherent prevision, then $\pr_\banlim$ is a (strongly) $T$-invariant coherent     prevision on $\gambles(\values)$.
  \item If $\alpr$ is a weakly $T$-invariant coherent lower prevision on $\gambles(\values)$, then the     (strongly) $T$-invariant coherent prevision $\pr_\banlim$ dominates $\alpr$ for any $\pr$ in     $\solp(\alpr)$.
  \item If $\lpr=\pr$ is a $T$-invariant coherent prevision, then $\pr_\banlim=\pr$.
  \end{enumerate}
\end{proposition}

\begin{proof}
  We first prove the first statement. Consider gambles $f$ and $g$ on $\values$. Since $\inf f\leq f_\lpr$, it   follows from the coherence of $\banlim$ that $\inf f\leq\banlim(f_\lpr)=\lpr_\banlim(f)$.  Moreover, we have   for any $n$ in $\nats$ that
  \begin{equation*}
    (f+g)_\lpr(n)=\lpr((f+g)\circ T^n)=\lpr(f\circ T^n+g\circ T^n)
    \geq\lpr(f\circ T^n)+\lpr(g\circ T^n)=f_\lpr(n)+g_\lpr(n),
  \end{equation*}
  where the inequality follows from the coherence [super-additivity] of $\lpr$. Since $\banlim$ is coherent,   we see that $\lpr_\banlim(f+g)\geq\banlim(f_\lpr)+\banlim(g_\lpr) =\lpr_\banlim(f)+\lpr_\banlim(g)$.   Finally, for any $\lambda\geq0$, we have that $(\lambda f)_\lpr(n)=\lpr((\lambda f)\circ T^n)   =\lpr(\lambda(f\circ T^n))=\lambda\lpr(f\circ T^n)=\lambda f_\lpr(n)$, since $\lpr$ is coherent.   Consequently $\lpr_\banlim(\lambda f)=\banlim(\lambda f_\lpr)   =\lambda\banlim(f_\lpr)=\lambda\lpr_\banlim(f)$, since $\banlim$ is coherent. This proves that   $\lpr_\banlim$ is a coherent lower prevision on $\gambles(\values)$ [because (P1)--(P3) are satisfied].  To   show that it is weakly $T$-invariant, recall that $(T^tf)_\lpr=\theta^t f_\lpr$, whence
  \begin{equation*}
    \lpr_\banlim(T^tf)=\banlim((T^tf)_{\lpr})
    =\banlim(\theta^tf_\lpr)=\banlim(f_\lpr)=\lpr_\banlim(f),
  \end{equation*}
  since $\banlim$ is shift-invariant.
  \par
  To prove the second statement, assume that $\lpr$ dominates the weakly $T$-invariant coherent lower   prevision $\alpr$ on $\gambles(\values)$. Then for any gamble $f$ on $\values$, we see that
  \begin{equation*}
    f_\lpr(n)=\lpr(f\circ T^n)\geq\alpr(f\circ T^n)\geq\alpr(f),
  \end{equation*}
  where the last inequality follows from the weak $T$-invariance of $\alpr$.  Consequently, since $\banlim$ is   coherent, we get $\lpr_\banlim(f)=\banlim(f_\lpr)\geq\alpr(f)$.
  \par
  The third statement follows immediately from the first and the fact that $\pr_\banlim$ is a self-conjugate   coherent lower prevision (and therefore a coherent prevision) because $\pr$ and $\banlim$ are.
  \par
  The fourth statement follows at once from the second and the third. The fifth is an immediate consequence of   the definition of $\pr_\banlim$.
\end{proof}

We can use the results in this proposition to characterise all strongly $T$-invariant coherent lower previsions using Banach limits on $\gambles(\nats)$.

\begin{theorem}\label{theo:all-invariants-are-banach}
  Let $\lpr$ be a weakly $T$-invariant coherent lower prevision defined on some $T$-invariant domain   $\domain$, that avoids sure loss. Then the set of all $T$-invariant coherent previsions on   $\gambles(\values)$ that dominate $\lpr$ on $\domain$ is given by
  \begin{equation*}
    \set{P_\banlim}
    {\text{$\pr\in\solp(\lpr)$ and $\banlim\in\linprevs_\theta(\nats)$}},
  \end{equation*}
  so the smallest strongly $T$-invariant coherent lower prevision $\lnex_{\lpr,T}$ on $\gambles(\values)$ that   dominates $\lpr$, i.e., the strongly $T$-invariant natural extension of $\lpr$, is the lower envelope of   this set, and also given by
  \begin{equation*}
    \lnex_{\lpr,T}(f)
    =\inf_{\pr\in\solp(\lpr)}\lnex_\theta(f_\pr)
    =\inf_{\pr\in\solp(\lpr)}\sup_{n\geq1}\inf_{k\geq0}
    \left[\frac{1}{n}\sum_{\ell=k}^{k+n-1}\pr((T^\ell)^tf)\right]
  \end{equation*}
  for any gamble $f$ on $\values$. As a consequence, the set $\linprevs_T(\values)$ of all $T$-invariant   coherent previsions on $\gambles(\values)$ is given by
  \begin{equation*}
    \linprevs_T(\values)
    =\set{P_\banlim}
    {\text{$\pr\in\linprevs$ and $\banlim\in\linprevs_\theta(\nats)$}}.
  \end{equation*}
  This tells us that all $T$-invariant coherent previsions can be constructed using Banach limits on   $\gambles(\nats)$. The smallest strongly $T$-invariant coherent lower prevision $\lnex_{T}$ on   $\gambles(\values)$ is the lower envelope of this set, and also given by
  \begin{equation*}
    \lnex_{T}(f)
    =\inf_{\pr\in\linprevs(\values)}\lnex_\theta(f_\pr)
    =\inf_{\pr\in\linprevs(\values)}\sup_{n\geq1}\inf_{k\geq0}
    \left[\frac{1}{n}\sum_{\ell=k}^{k+n-1}\pr((T^\ell)^tf)\right]
  \end{equation*}
  for any gamble $f$ on $\values$.
\end{theorem}

\begin{proof}
  First of all, a coherent prevision $\pr$ on $\gambles(\values)$ belongs to $\solp(\lpr)$, i.e., dominates   $\lpr$ on its domain $\domain$, if and only if $\pr$ dominates the natural extension $\lnex_\lpr$ on all   gambles.  Moreover, $\lnex_\lpr$ is weakly $T$-invariant by   Theorem~\ref{theo:natext-preserves-invariance}. Now consider any $\pr\in\solp(\lpr)$. Use the above   observations together with Proposition~\ref{prop:bhaskara-banlim} [statements 3 and 4] to show that for any   Banach limit $\banlim$ on $\gambles(\nats)$, $\pr_\banlim$ is a $T$-invariant coherent prevision that   dominates $\lpr$. Conversely, if $\pr$ is a $T$-invariant coherent prevision on $\gambles(\values)$ that   dominates $\lpr$ on $\domain$, then by Proposition~\ref{prop:bhaskara-banlim} [statement 5],   $\pr=\pr_\banlim$ for any Banach limit $\banlim$ on $\gambles(\nats)$. This shows that   $\set{\pr_\banlim}{\pr\in\solp(\lpr), \banlim\in\linprevs_\theta(\nats)}$ is indeed the set of $T$-invariant   coherent previsions on $\gambles(\values)$ that dominate $\lpr$ on $\domain$. Consequently, $\lnex_{\lpr,T}$   is the lower envelope of this set, whence for any gamble $f$ on $\values$
  \begin{align*}
    \lnex_{\lpr,T}(f) &=\inf_{\pr\in\solp(\lpr)} \inf_{\banlim\in\linprevs_\theta(\nats)}\pr_\banlim(f)     =\inf_{\pr\in\solp(\lpr)}
    \inf_{\banlim\in\linprevs_\theta(\nats)}\banlim(f_\pr)\\
    \intertext{and since $\lnex_\theta$ is the lower envelope of $\linprevs_\theta(\nats)$,}
    &=\inf_{\pr\in\solp(\lpr)}\lnex_\theta(f_\pr)\\
    \intertext{and using Eqs.~\eqref{eq:shift-invariant-natex} and~\eqref{eq:bhaskara-gamble},}     &=\inf_{\pr\in\solp(\lpr)}\sup_{n\geq1}\inf_{k\geq0}     \left[\frac{1}{n}\sum_{\ell=k}^{k+n-1}\pr((T^\ell)^tf)\right].
  \end{align*}
  The rest of the proof is now immediate.
\end{proof}

\subsection{Generalised Banach limits}
The above results on monoids $\transfos_T$ generated by a single transformation $T$ can be generalised towards more general monoids $\transfos$ of transformations of $\values$, such that the set $\transfos^*$ of convex mixtures of the lifted linear transformations in $\transfos^t$ is directed by the successor relation $\succeeds$ on $\transfos^*$. The following discussion establishes an interesting connection between strong invariance and the notion of a generalised Banach limit.
\par
We can consider $\transfos^*$ as a monoid of transformations of itself, as follows: with any element $T^*$ we associate a transformation of $\transfos^*$, also denoted by $T^*$, such that $T^*(S^*):=S^*T^*\in\transfos^*$, for any $S^*$ in $\transfos^*$.\footnote{Usually, $T^*(S^*)$ is defined as   $T^*S^*$, see for instance \citet[Note~1 of Section~3.5.1]{walley1991}. But we have to take a different   route here because the elements of $\transfos^*$ are convex mixtures of \emph{lifted} transformations, and   as we have seen, lifting reverses the order of application of transformations.} We can, in the usual fashion, lift $T^*$ to a transformation $(T^*)^t$ on $\gambles(\transfos^*)$ by letting $(T^*)^tg=g\circ T^*$, or in other words
\begin{equation}\label{eq:abstract-lifting}
  (T^*)^tg(S^*)=g(T^*(S^*))=g(S^*T^*),
\end{equation}
for any $S^*$ in $\transfos^*$ and any gamble $g$ on $\transfos^*$, i.e., $g\in\gambles(\transfos^*)$.
\par
Now a \emph{generalised Banach limit} \citep[Sections~12.33--12.38]{schechter1997} on $\gambles(\transfos^*)$ is defined as any linear functional on $\gambles(\transfos^*)$ that dominates the limit inferior operator with respect to the directed set $\transfos^*$. Let us take a closer look at this limit inferior operator. It is defined by
\begin{equation*}
  \liminf_{\transfos^*}g
  =\liminf_{T^*\in\transfos^*}g(T^*)
  :=\sup_{S^*\in\transfos^*}\inf_{T^*\succeeds S^*}g(T^*),
\end{equation*}
for any gamble $g$ on $\transfos^*$.  Now recall that $T^*\succeeds S^*$ if and only if there is some $R^*$ in $\transfos^*$ such that $T^*=R^*S^*$, so we get, using Eq.~\eqref{eq:abstract-lifting}, that
\begin{equation*}
  \liminf_{T^*\in\transfos^*}g(T^*)
  =\sup_{S^*\in\transfos^*}\inf_{R^*\in\transfos^*}g(R^*S^*)
  =\sup_{S^*\in\transfos^*}\inf_{R^*\in\transfos^*}(S^*)^tg(R^*)
  =\lim_{S^*\in\transfos^*}\lpr_{\transfos^*}((S^*)^tg),
\end{equation*}
where $\lpr_{\transfos^*}$ is the vacuous lower prevision on $\gambles(\transfos^*)$. If we look at Corollary~\ref{cor:invnatex} for the special case $\values=\transfos^*$ and the monoid of transformations $\transfos^*$, recall that we need to lift transformations in $\transfos^*$ before we can apply them to gambles, and that the lifted transformations of $\transfos^*$ already constitute a convex set\footnote{In   general, even if $\transfos^t$ is directed by the successor relation $\succeeds$, the limit inferior   operator on $\gambles(\transfos^t)$ will not be strongly invariant.  But convexification, or going from   $\transfos^t$ to $\transfos^*$, makes the limit inferior strongly invariant. Observe in this respect that   the limit inferior operator on $\gambles(\nats)$ is not strongly shift-invariant, but its `convexified'   counterpart $\lnex_\theta$ is.}, we easily get to the following conclusion.

\begin{proposition}
  The limit inferior operator on $\gambles(\transfos^*)$ is actually the point-wise smallest strongly   $\transfos^*$-invariant coherent lower prevision on $\gambles(\transfos^*)$, and the generalised Banach   limits on $\gambles(\transfos^*)$ are the $\transfos^*$-invariant coherent previsions on   $\gambles(\transfos^*)$.
\end{proposition}

We can now apply arguments similar to the ones in the previous section, for general monoids $\transfos$ of transformations of $\values$ such that $\transfos^*$ is directed. Consider any coherent lower prevision $\lpr$ on $\gambles(\values)$ and any gamble $f$, and define the following gamble $f_\lpr$ on $\transfos^*$:
\begin{equation*}
  f_\lpr(S^*):=\lpr(S^*f)
\end{equation*}
for any $S^*$ in $\transfos^*$, which generalises Eq.~\eqref{eq:bhaskara-gamble}. Observe that, using Eq.~\eqref{eq:abstract-lifting},
\begin{equation*}
  (T^*f)_\lpr(S^*)=\lpr(S^*T^*f)=f_\lpr(S^*T^*)=(T^*)^tf_\lpr(S^*),
\end{equation*}
so
\begin{equation*}
  (T^*f)_\lpr=(T^*)^tf_\lpr,
\end{equation*}
which generalises Eq.~\eqref{eq:bhaskara-lpr}. If we consider any $\transfos^*$-invariant coherent prevision $\banlim$ on $\gambles(\transfos^*)$, or in other words a generalised Banach limit on $\gambles(\transfos^*)$, we can now define a new lower prevision $\lpr_\banlim$ on $\gambles(\values)$ by $\lpr_\banlim(f):=\banlim(f_\lpr)$, and Proposition~\ref{prop:bhaskara-banlim}, as well as Theorem~\ref{theo:all-invariants-are-banach}, can now easily be generalised from monoids of transformations with a single generator to arbitrary directed monoids. In particular, we find that
\begin{equation*}
  \lnex_{\lpr,\transfos}(f)
  =\inf_{\pr\in\solp(\lpr)}\liminf_{T^*\in\transfos^*}\pr(T^*f)
  \text{ and }
  \lnex_{\transfos}(f)
  =\inf_{\pr\in\linprevs(\values)}\liminf_{T^*\in\transfos^*}\pr(T^*f)
\end{equation*}
for any gamble $f$ on $\values$, where $\lpr$ is any weakly $\transfos$-invariant lower prevision that avoids sure loss.

\section{Permutation invariance on finite spaces}\label{sec:finite}
Assume now that $\transfos$ is a finite group $\permuts$ of permutations of $\values$.  Then we have the following characterisation result for the weakly $\permuts$-invariant coherent lower previsions.

\begin{theorem}\label{theo:finite-group}
  Let $\permuts$ be a finite group of permutations of $\values$.  All weakly $\permuts$-invariant coherent   lower previsions $\alpr$ on $\gambles(\values)$ have the form
  \begin{equation}\label{eq:invariant-mixture}
    \alpr=\frac{1}{\abs{\permuts}}\sum_{\pi\in\permuts}\pi\lpr,
  \end{equation}
  where $\abs{\permuts}$ is the number of permutations in $\permuts$, and $\lpr$ is any coherent lower   prevision on $\gambles(\values)$.
\end{theorem}

\begin{proof}
  Consider a coherent lower prevision $\lpr$ on $\gambles(\values)$, and let $\alpr$ be the corresponding   lower prevision, given by Eq.~\eqref{eq:invariant-mixture}. Then $\alpr$ is coherent, as a convex mixture of   coherent lower previsions $\pi\lpr$. Moreover, let $\varpi$ be any element of $\permuts$, then
  \begin{equation*}
    \varpi\alpr
    =\frac{1}{\abs{\permuts}}\sum_{\pi\in\permuts}(\varpi\pi)\lpr,
    =\frac{1}{\abs{\permuts}}\sum_{\pi\in\varpi\permuts}\pi\lpr,
  \end{equation*}
  where $\varpi\permuts=\set{\varpi\pi}{\pi\in\permuts}=\permuts$, because $\permuts$ is a group of   permutations. Consequently $\varpi\alpr=\alpr$, so $\alpr$ is weakly $\permuts$-invariant.
  \par
  Conversely, let $\alpr$ be any weakly $\permuts$-invariant coherent lower prevision, then we recover $\alpr$   on the left-hand side if we insert $\alpr$ in the right-hand side of Eq.~\eqref{eq:invariant-mixture}. So   any weakly $\permuts$-invariant coherent lower prevision is indeed of the form~\eqref{eq:invariant-mixture}.
\end{proof}

Next, we give an interesting representation result for the strongly $\permuts$-invariant coherent lower previsions, when in addition, $\values$ is a finite set.\footnote{We find the `permutation symmetry' between   Theorems~\ref{theo:finite-group} and~\ref{theo:group-finite-space} quite surprising: the former states that   a weakly $\permuts$-invariant coherent lower prevision is a uniform prevision (or mixture) of coherent lower   previsions, and the latter that a strongly $\permuts$-invariant coherent lower prevision is a coherent lower   prevision of uniform previsions.} As we shall see further on, this essentially simple result has many interesting consequences, amongst which a generalisation to coherent lower previsions of \citegen{finetti1937} representation result for finite sequences of exchangeable random variables (see Section~\ref{sec:exchangeability}). Recall that $\atoms_\permuts$ is the set of all $\permuts$-invariant atoms of $\values$. For each $A$ in $\atoms_\permuts$, define $\pr^u(\cdot\vert A)$ as the coherent prevision on $\gambles(\values)$ all of whose probability mass is uniformly distributed over $A$, i.e., for all gambles $f$ on $\values$:
\begin{equation*}
  \pr^u(f\vert A)=\frac{1}{\abs{A}}\sum_{x\in A}f(x).
\end{equation*}
Finally, let $\pr^u(f\vert\atoms_\permuts)$ denote the gamble on $\atoms_\permuts$ that assumes the value $\pr^u(f\vert\atoms_\permuts)(A):=\pr^u(f\vert A)$ in any element $A$ of $\atoms_\permuts$.

\begin{theorem}\label{theo:group-finite-space}
  Let $\permuts$ be a group of permutations of the finite set $\values$. A coherent lower prevision on   $\gambles(\values)$ is strongly $\permuts$-invariant if and only if   $\lpr(f)=\lpr_0(\pr^u(f\vert\atoms_\permuts))$ for all $f$ in $\gambles(\values)$, where $\lpr_0$ is an   arbitrary coherent lower prevision on $\gambles(\atoms_\permuts)$.
\end{theorem}

\begin{proof}
  We begin with the `if' part. Let $\lpr_0$ be an arbitrary coherent lower prevision on   $\gambles(\atoms_\permuts)$, and suppose that $\lpr=\lpr_0(\pr^u(\cdot\vert\atoms_\permuts))$. Then it is   easy to see that $\lpr$ is coherent. We show that $\lpr$ is strongly $\permuts$-invariant.  Consider any   gamble $f$ on $\values$ and any $\pi\in\permuts$. Then for any $A$ in $\atoms_\permuts$ and any gamble $f$   on $\values$,
  \begin{equation*}
    \pr^u(f-\pi^tf\vert A)
    =\frac{1}{\abs{A}}\sum_{x\in A}[f(x)-f(\pi x)]=0,
  \end{equation*}
  because $x\in A$ is equivalent to $\pi x\in A$. So we see that $\lpr(f-\pi^tf)=\lpr_0(0)=0$, since $\lpr_0$   is coherent. In a similar way, we can prove that $\lpr(\pi^tf-f)=0$, so $\lpr$ is indeed strongly   $\permuts$-invariant.
  \par
  To prove the `only if' part, we first concentrate on the case of a $\permuts$-invariant coherent prevision   $\pr$ on $\gambles(\values)$. Fix any gamble $f$ on $\values$. Since $\pr$ is a coherent prevision, we find   that
  \begin{equation*}
    f=\sum_{A\in\atoms_\permuts}fI_A
    \quad\text{and}\quad
    \pr(f)=\sum_{A\in\atoms_\permuts}\pr(fI_A)
    =\sum_{A\in\atoms_\permuts}\pr(f\vert A)\pr(A),
  \end{equation*}
  where we have used Bayes's rule to define $\pr(f\vert A):=\pr(fI_A)/\pr(A)$ if $\pr(A)>0$ and $\pr(f\vert   A)$ is arbitrary otherwise.
  \par
  Now assume that $\pr$ is $\permuts$-invariant. Fix any $\permuts$-invariant atom $A$ in $\atoms_\permuts$   such that $\pr(A)>0$ and let $\pi\in\permuts$.  For any gamble $f$ on $\values$, we see that   $\pi^t(fI_A)=(\pi^tf)I_A$, since $A$ is in particular $\pi$-invariant.  Consequently
  \begin{equation*}
    \pr(\pi^tf\vert A)=\pr((\pi^tf)I_A)/\pr(A)=\pr(\pi^t(fI_A))/\pr(A)
    =\pr(fI_A)/\pr(A)=\pr(f\vert A),
  \end{equation*}
  so $\pr(\cdot\vert A)$ is $\permuts$-invariant as well.\footnote{\label{fn:conditioning}This is an instance     of a more general result, namely that coherent conditioning of a coherent lower prevision on an invariant     event preserves both weak and strong invariance. A proof of this statement is not difficult, but outside     the scope of this paper.}  Now let for any $y$ in the finite set $A$, $p(y\vert A):=\pr(\{y\}\vert   A)\geq0$, then on the one hand $\sum_{x\in A}p(x\vert A)=\pr(A\vert A)=1$.  On the other hand, it follows   from the $\pi$-invariance of $\pr(\cdot\vert A)$ that $p(x\vert A)=p(\pi x\vert A)$ for any $x$ in   $A$. Since we know from Proposition~\ref{prop:invariant-atoms} that $A=\set{\pi x}{\pi\in\permuts}$, we see   that $p(\cdot\vert A)$ is constant on $A$, so $p(x\vert A)=1/\abs{A}$ for all $x$ in $A$, and consequently   $P(f\vert A)=\pr^u(f\vert A)$, whence $\pr(f)=\sum_{A\in\atoms_\permuts}\pr^u(f\vert A)\pr(A)$. So indeed   there is a coherent prevision $\pr_0$ on $\gambles(\atoms_\permuts)$, defined by $\pr_0(\{A\})=\pr(A)$ for   all $A\in\atoms_\permuts$, such that $\pr=\pr_0(\pr^u(\cdot\vert\atoms_\permuts))$.
  \par
  Finally, let $\lpr$ be any strongly $\permuts$-invariant coherent lower prevision, so any   $\pr\in\solp(\lpr)$ is $\permuts$-invariant and can therefore be written as   $\pr=\pr_0(\pr^u(\cdot\vert\atoms_\permuts))$. If we let $\lpr_0$ be the (coherent) lower envelope of the   set $\set{\pr_0}{\pr\in\solp(\lpr)}$, then since $\lpr$ is the lower envelope of $\solp(\lpr)$, we get   immediately that $\lpr=\lpr_0(\pr^u(\cdot\vert\atoms_\permuts))$.
\end{proof}

As an immediate corollary, we see that that the uniform coherent prevision $\pr^u$ on $\gambles(\values)$ is the only strongly $\permuts$-invariant coherent lower prevision on $\gambles(\values)$ if and only if $\values$ is the only $\permuts$-invariant atom, i.e., if $\atoms_\permuts=\{\values\}$.  This is for instance the case if $\permuts$ is the group of all permutations of $\values$, or more generally if $\permuts$ includes the cyclic group of permutations of $\values$. It should therefore come as no surprise that, since symmetry of beliefs is so often confused with beliefs of symmetry, the uniform distribution is so often (but wrongly so) considered to be a good model for complete ignorance.
\par
Another immediate corollary of this result is that the smallest strongly $\permuts$-invariant coherent lower prevision on $\gambles(\values)$ is given by $\lpr(f)=\inf_{A\in\atoms_\permuts}\frac{1}{\abs{A}}\sum_{x\in   A}f(x)$, which of course agrees with the uniform distribution when we let $\permuts$ be the group of all permutations.
\par
These results do not extend to the case where we have transformations of $\values$ that are not permutations; as we have said before, as soon as we have two different constant transformations in the monoid $\transfos$, there are no strongly invariant belief models.

\subsection{A few simple examples}
We now apply the theorems above in a number of interesting and simple examples.

\begin{example}\label{ex:two-elements}
  Let $\values=\values_2:=\{1,2\}$, then all coherent lower previsions on $\gambles(\values_2)$ are so-called   \emph{linear-vacuous mixtures}, i.e., convex combinations of a coherent (linear) prevision and the vacuous   lower prevision, and therefore given by
  \begin{equation*}
    \lpr(f)=\epsilon\left[\alpha f(1)+(1-\alpha)f(2)\right]
    +(1-\epsilon)\min\{f(1),f(2)\},
  \end{equation*}
  where $0\leq\alpha\leq1$ and $0\leq\epsilon\leq1$.  Let $\permuts_2$ be the set of all permutations of   $\values_2$. Then the only strongly $\permuts_2$-invariant coherent lower prevision is the uniform coherent   prevision
  \begin{equation*}
    \pr_{\frac{1}{2}}(f)=\frac{1}{2}[f(1)+f(2)],
  \end{equation*}
  corresponding to $\alpha=\frac{1}{2}$ and $\epsilon=1$. The weakly $\permuts_2$-invariant coherent lower   previsions are given by
  \begin{equation*}
    \lpr(f)=\epsilon\pr_{\frac{1}{2}}(f)+(1-\epsilon)\min\{f(1),f(2)\},
  \end{equation*}
  where $0\leq\epsilon\leq1$, so they are all the convex mixtures of the uniform coherent prevision and the   vacuous lower prevision. $\blacklozenge$
\end{example}

\begin{example}\label{ex:three-elements}
  Let $\values=\values_3:=\{1,2,3\}$, then all 2-monotone coherent lower previsions on $\gambles(\values_3)$   are given by\footnote{An explicit proof of this statement is beyond the scope of this paper, but it runs     along the following lines: (i) any coherent lower probability on the set of all events of a three-element     space is 2-monotone \citep[p.~58]{walley1981}; (ii) all 2-monotone coherent lower probabilities make up a     convex set, and are convex mixtures of the extreme points of this set \citep[Chapter~2]{maass2003} (By the     way, an argument similar to that in \citet[Chapter~2]{maass2003} shows that all strongly     $\transfos$-invariant coherent lower previsions are (infinite) convex mixtures of the extreme strongly     $\transfos$-invariant coherent lower previsions.); (iii) the 2-monotone coherent lower previsions on all     gambles are natural extensions of the 2-monotone coherent lower previsions on all events     \Citep{walley1981,cooman2005e,cooman2005d,cooman2005b}; and (iv) natural extension to gambles of     2-monotone lower probabilities preserves convex mixtures.}
  \begin{multline*}
    \lpr(f)=m_1f(1)+m_2f(2)+m_3f(3)\\
    +m_4\min\{f(1),f(2)\}+m_5\min\{f(2),f(3)\}+m_6\min\{f(3),f(1)\}\\
    +m_7\min\left\{\frac{f(1)+f(2)}{2},\frac{f(2)+f(3)}{2},
      \frac{f(3)+f(1)}{2}\right\}\\
    +m_8\min\{f(1),f(2),f(3)\}.
  \end{multline*}
  where $0\leq m_k\leq1$ and $\sum_{k=1}^8m_k=1$.  Let $\permuts_3$ be the set of all permutations of   $\values_3$. Then the only strongly $\permuts_3$-invariant coherent lower prevision is the uniform coherent   prevision
  \begin{equation*}
    \pr(f)=\frac{1}{3}[f(1)+f(2)+f(3)],
  \end{equation*}
  corresponding to $m_1=m_2=m_3=\frac{1}{3}$ and $m_4=m_5=m_6=m_7=m_8=0$ [Observe that a coherent prevision is   always 2-monotone.]. Weak $\permuts_3$- invariance, on the other hand, requires only that $m_1=m_2=m_3$ and   $m_4=m_5=m_6$, so all the weakly $\permuts_3$-invariant and 2-monotone coherent lower previsions are given   by
  \begin{multline*}
    \lpr(f)=\frac{M_1}{3}[f(1)+f(2)+f(3)]\\
    +\frac{M_2}{3}\left[\min\{f(1),f(2)\}
      +\min\{f(2),f(3)\}+\min\{f(3),f(1)\}\right]\\
    +M_3\min\left\{\frac{f(1)+f(2)}{2},\frac{f(2)+f(3)}{2},
      \frac{f(3)+f(1)}{2}\right\}\\
    +M_4\min\{f(1),f(2),f(3)\}.
  \end{multline*}
  where $0\leq M_k\leq1$ and $M_1+M_2+M_3+M_4=1$. The weakly $\permuts_3$-invariant and completely monotone   coherent lower previsions (natural extensions of belief functions) correspond to the choice $M_3=0$.   $\blacklozenge$
\end{example}

\begin{example}\label{ex:six-elements}
  Consider rolling a die for which there is evidence of symmetry between all even numbers, on the one hand,   and between all odd numbers on the other.  Let $\values=\values_6:=\{1,\dots,6\}$ and let   $\permuts_{\mathrm{eo}}$ be the set of all permutations of $\values_6$ that map even numbers to even numbers   and odd numbers to odd numbers. The $\permuts_{\mathrm{eo}}$-invariant atoms are $\{1,3,5\}$ and   $\{2,4,6\}$. By Theorem~\ref{theo:group-finite-space}, the strongly $\permuts_{\mathrm{eo}}$-invariant   coherent previsions on $\gambles(\values_6)$, which are the precise belief models that are compatible with   the subject's beliefs of symmetry, are given by
  \begin{equation*}
    \pr(f)=\frac{\alpha}{3}[f(1)+f(3)+f(5)]
    +\frac{1-\alpha}{3}[f(2)+f(4)+f(6)],
  \end{equation*}
  where $0\leq\alpha\leq1$, and more generally, the strongly $\permuts_{\mathrm{eo}}$-invariant coherent lower   previsions on $\gambles(\values_6)$ are [apply Theorem~\ref{theo:group-finite-space} and use the results in   Example~\ref{ex:two-elements}]
  \begin{multline*}
    \lpr(f)=\epsilon\left[\frac{\alpha}{3}[f(1)+f(3)+f(5)]
      +\frac{1-\alpha}{3}[f(2)+f(4)+f(6)]\right]\\
    +(1-\epsilon)\min\left\{\frac{f(1)+f(3) +f(5)}{3},\frac{f(2)+f(4)+f(6)}{3}\right\}
  \end{multline*}
  for $0\leq\epsilon\leq1$ and $0\leq\alpha\leq1$. $\blacklozenge$
\end{example}

\begin{example}
  Let us show that the point-wise smallest strongly invariant coherent lower prevision extension is not   necessarily 2-monotone. Consider $\values_4:=\{1,2,3,4\}$, and let $\pi$ be the permutation of $\values_4$   defined by $\pi(1)=2$, $\pi(2)=1$, $\pi(3)=4$ and $\pi(4)=3$. Observe that $\pi$ is its own inverse, so   $\transfos_\pi=\{\iden_{\values_4},\pi\}$ is a group. From Theorem~\ref{theo:group-finite-space} we infer   that the point-wise smallest strongly $\pi$-invariant coherent lower prevision on all gambles is given by
  \begin{equation*}
    \lnex_\pi(f)=\min\left\{\frac{f(1)+f(2)}{2},\frac{f(3)+f(4)}{2}\right\}.
  \end{equation*}
  Let us now consider the gambles $f_1$ and $f_2$ on $\values_4$, given by $f_1(1)=0$, $f_1(2)=-1$,   $f_1(3)=1$, $f_1(4)=-1$ and $f_2(1)=-1$, $f_2(2)=-0.25$, $f_2(3)=-1.5$, $f_2(4)=0$. Check that
  \begin{equation*}
    \lnex_\pi(f_1\wedge f_2)+\lnex_\pi(f_1\vee f_2)
    =-1.25-0.125=-1.375
    <-0.5-0.75=\lnex_\pi(f_1)+\lnex_\pi(f_2).
  \end{equation*}
  Hence, $\lnex_\pi$ is not $2$-monotone.  $\blacklozenge$
\end{example}

The following example shows that possibility measures are not very useful for modelling permutation invariance.

\begin{example}
  Consider a possibility measure $\Pi$ defined on all events of a finite space $\values$. Then there is a map   $\lambda\colon\values\to\reals^+$, called the \emph{possibility distribution} of $\Pi$, such that   $\lambda(x):=\Pi(\{x\})$ and moreover $\Pi(A)=\max_{x\in A}\lambda(x)$ for all non-empty events   $A\subseteq\values$. We have mentioned before that $\Pi$ is a coherent upper probability if and only if   $\Pi(\values)=\max_{x\in\values}\lambda(x)=1$.  We shall assume this is the case. Now consider any group   $\permuts$ of permutations of $\values$. Then clearly $\Pi$ is weakly $\permuts$-invariant if and only if   $\lambda$ is constant on the $\permuts$-invariant atoms of $\values$. In particular, $\Pi$ is weakly   invariant with respect to all permutations if and only is $\lambda$ is everywhere equal to one, so $\Pi$ is   the vacuous upper probability.
  \par
  For strong $\permuts$-invariance, let $\upr$ be any strongly $\permuts$-invariant coherent lower prevision   whose domain contains at least all events. Let $x$ be any element of $\values$, and let $\atom{x}{\permuts}$   be the $\permuts$-invariant atom that contains $x$.  Then it follows from   Theorem~\ref{theo:group-finite-space} that $\upr(\{x\})\leq1/\abs{\atom{x}{\permuts}}$. So for $\upr$ to   extend a possibility measure, it is necessary (but not sufficient) that there is at least one element $z$ of   $\values$ such that $\upr(\{z\})=1$, implying that $z$ should be left invariant by all the permutations in   $\permuts$, or equivalently, $\atom{z}{\permuts}=\{z\}$. $\blacklozenge$
\end{example}

\subsection{Exchangeable lower previsions}\label{sec:exchangeability}
As another example, we now discuss the case of so-called exchangeable coherent lower previsions. Consider a non-empty finite set $\values_\kappa:=\{1,\dots,\kappa\}$ of categories, and $N$ random variables $\rv_1$, \dots, $\rv_N$ taking values in the same set $\values_\kappa$, where $\kappa$ and $N$ are natural numbers with $\kappa\geq2$ and $N\geq1$. The joint random variable $\mathbf{\rv}:=(\rv_1,\dots,\rv_N)$ assumes values in the set $\values:=\values_\kappa^N$.\footnote{This means that we assume these $N$ random variables to be   \emph{logically independent}.} We want to model a subject's beliefs about the value that $\mathbf{\rv}$ assumes in $\values_\kappa^N$, and generally, we use a coherent lower prevision $\lpr$ on $\gambles(\values_\kappa^N)$ to represent such beliefs.
\par
Now assume that our subject believes that all random variables $\rv_k$ are generated by the same process at different times $k$, and that the properties of this process do not depend on the time $k$. So, the subject assesses that there is permutation symmetry between the different times $k$. How can such \emph{beliefs of   symmetry} be modelled?
\par
With a permutation $\pi$ of $\{1,\dots,N\}$, we can associate (by the usual procedure of lifting) a permutation of $\values=\values_\kappa^N$, also denoted by $\pi$, that maps any $\sample{x}=(x_1,\dots,x_N)$ in $\values_\kappa^N$ to $\pi\sample{x}:=(x_{\pi(1)},\dots,x_{\pi(N)})$. The belief models that are compatible with the subject's beliefs of symmetry, are therefore the coherent lower previsions on (subsets of) $\gambles(\values_\kappa^N)$ that are \emph{strongly} $\permuts_\kappa^N$-invariant, where $\permuts_\kappa^N$ is the group of liftings to $\values_\kappa^N$ of all permutations of $\{1,\dots,N\}$. \citet[Chapter~9]{walley1991} calls such lower previsions \emph{exchangeable}, as they generalise \citegen{finetti1937} notion of exchangeable coherent previsions. We intend to characterise the exchangeable lower previsions using Theorem~\ref{theo:group-finite-space}. This will lead us to a generalisation (Eq.~\eqref{eq:exchangeability}) of \citegen{finetti1937} representation result for finite numbers of exchangeable random variables.
\par
It should be mentioned here that we should, as always, clearly distinguish between `beliefs of symmetry' and `symmetry of beliefs'. The latter imposes much weaker requirements on coherent lower previsions, namely those of weak $\permuts_\kappa^N$-invariance, which is called \emph{permutability} by \citet[Chapter~9]{walley1991}.\footnote{See \citet[Chapter~9]{walley1991} for a much more detailed discussion   of the difference between permutability and exchangeability.} In particular, the permutation symmetry that goes along with ignorance can only be invoked to justify permutability, but not, of course, exchangeability. Observe in this respect that the vacuous lower prevision on $\gambles(\values_\kappa^N)$ is permutable, but not exchangeable.  It is well-known (see for instance \cite{zabell1989,zabell1992}), that Laplace's Rule of Succession can be obtained by updating a particular exchangeable coherent prevision, but it should be clear from the discussion in this paper that ignorance alone (the Principle of Insufficient Reason) cannot be invoked to justify using such an exchangeable prevision, as (with considerable hindsight) Laplace implicitly seems to have done (see for instance \citet{howie2002,zabell1989,zabell1992}).
\par
For any $\sample{x}=(x_1,\dots,x_N)$ in $\values_\kappa^N$, the $\permuts_\kappa^N$-invariant atom $\atom{\sample{x}}{\permuts_\kappa^N}$ is the set of all permutations of (the components of) $\sample{x}$. If we define the set of possible \emph{count vectors}
\begin{equation*}
  \counts^N_\kappa
  =\set{(m_1,\dots,m_\kappa)}
  {m_k\in\nats^+\text{ and }\sum_{k=1}^\kappa m_k=N}
\end{equation*}
and the \emph{counting map} $\cntf\colon\values_\kappa^N\to\counts^N_\kappa$ such that $\cntf(x_1,\dots,x_N)$ is the $\kappa$-tuple, whose $k$-th component is given by
\begin{equation*}
  T_k(x_1,\dots,x_N)=\abs{\set{\ell\in\values_\kappa}{x_\ell=k}},
\end{equation*}
i.e., the number of components of $x$ whose value is $k$, then the number of elements of the invariant atom $\atom{\sample{x}}{\permuts_\kappa^N}$ is precisely
\begin{equation*}
  \nu(\cntf(\sample{x}))
  :=\dbinom{N}{T_1(\sample{x})\dots T_\kappa(\sample{x})}
  =\dfrac{N!}{T_1(\sample{x})!\dots T_\kappa(\sample{x})!}
\end{equation*}
and $\cntf$ is a bijection (one-to-one and onto) between $\atoms_{\permuts_\kappa^N}$ and $\counts^N_\kappa$. An invariant atom is therefore completely identified by the count vector $\cntf(\sample{x})$ of any of its elements $\sample{x}$, and we shall henceforth denote the invariant atoms of $\values_\kappa^N$ by $[\cnt{m}]$, where $\cnt{m}=(m_1,\dots,m_\kappa)\in\counts^N_\kappa$, and $\sample{x}\in[\cnt{m}]$ if and only if $\cntf(\sample{x})=\cnt{m}$.
\par
The coherent prevision $\pr^u(\cdot\vert\cnt{m})$ on $\gambles(\values_\kappa^N)$ whose probability mass is uniformly distributed over the invariant atom $[\cnt{m}]$ is given by
\begin{equation*}
  \pr^u(f\vert\cnt{m})
  =\frac{1}{\nu(\cnt{m})}\sum_{\sample{x}\in[\cnt{m}]}f(\sample{x}).
\end{equation*}
Interestingly, this is the precise prevision that is associated with taking $N$ a-select drawings without replacement from an urn with $N$ balls, $m_1$ of which are of type $1$, \dots, and $m_\kappa$ of which are of type $\kappa$.  Theorem~\ref{theo:group-finite-space} now tells us that any exchangeable coherent lower prevision $\lpr$ on $\gambles(\values_\kappa^N)$ can be written as
\begin{equation}\label{eq:exchangeability}
  \lpr(f)=\lpr^N_\kappa(\pr^u(f\vert\counts^N_\kappa)),
\end{equation}
where $\lpr^N_\kappa$ is some coherent lower prevision on $\gambles(\counts^N_\kappa)$. This means that such \emph{an exchangeable lower prevision can be associated with $N$ a-select drawings from an urn with $N$ balls   of types $1$, \dots, $\kappa$, whose composition $\cnt{m}$ is unknown, but for which the available   information about the unknown composition is modelled by a coherent lower prevision $\lpr^N_\kappa$.}
\par
That exchangeable coherent previsions can be interpreted in terms of sampling without replacement from an urn with unknown composition, is actually well-known, and essentially goes back to \citet{finetti1937}. \citet{heath1976} give a simple proof for random variables that may assume two values. But we believe our proof\footnote{\citet[Chapter~9]{walley1991} also mentions this result for exchangeable coherent lower   previsions. The essence of his argument is similar to what we do in the last paragraph of the proof of   Theorem~\ref{theo:group-finite-space}.} for the more general case of exchangeable coherent \emph{lower} previsions and random variables that may assume \emph{more than two values}, is conceptually even simpler than Heath and Sudderth's proof, even though it is a special case of a much more general representation result (Theorem~\ref{theo:group-finite-space}). The essence of the present proof in the special case of coherent previsions $\pr$ is captured wonderfully well by \citegen[Section~3.1]{zabell1992} succinct statement: ``Thus $\pr$ is exchangeable if and only if two sequences having the same frequency vector have the same probability.''
\par
Our subject's beliefs could, in addition, be symmetrical in the categories in $\values_\kappa=\{1,\dots,\kappa\}$, for instance as a result of her ignorance about the process that generates the outcomes $\rv_k$ at each time $k$. As we have seen, this will be typically represented by using a type of \emph{weakly} invariant belief models, in this case with respect to permutations of the categories, rather than the times.  Any permutation $\varpi$ of $\values_\kappa$ induces a permutation of $\values_\kappa^N$, also denoted by $\varpi$, through
\begin{equation*}
  \varpi\sample{x}=\varpi(x_1,\dots,x_N)
  :=(\varpi(x_1),\dots,\varpi(x_N)).
\end{equation*}
What happens if we require that $\lpr$, in addition to being exchangeable, should also be \emph{weakly} invariant under all such permutations? It is not difficult to prove that
\begin{equation*}
  \pr^u(\varpi^{-1}f\vert\cnt{m})=\pr^u(f\vert\varpi\cnt{m}),
\end{equation*}
where we let $\varpi\cnt{m}=\varpi(m_1,\dots,m_\kappa) :=(m_{\varpi(1)},\dots,m_{\varpi(\kappa)})$ in the usual fashion.  This implies that there is such weak invariance if and only if the coherent lower prevision $\lpr_\kappa^N$ on $\gambles(\counts_\kappa^N)$ is \emph{weakly} invariant with respect to all category permutations! In particular, this weak invariance is satisfied for the vacuous lower prevision on $\gambles(\counts_\kappa^N)$. Another type of lower coherent prevision that exhibits such a combination of strong invariance for time permutations and weak invariance for category permutations, and which also has other very special and interesting properties, is constructed by taking lower envelopes of specific sets of Dirichlet-Multinomial distributions, leading to the so-called Imprecise Dirichlet-Multinomial Model (IDMM, see \citet{walley1999}).
\par
In the literature, however, it is sometimes required that a coherent precise prevision should be invariant with respect to the combined action of the permutations of times and categories. These are the so-called \emph{partition exchangeable} previsions (see \citet{zabell1992} for an interesting discussion and historical overview). Of course, the generalisation of this notion to coherent lower previsions should be strongly invariant with respect to such combined permutations, and therefore be a lower envelope of partition exchangeable previsions. For such \emph{partition exchangeable} lower previsions, Theorem~\ref{theo:group-finite-space} can be invoked to prove a representation result that is similar to that for coherent lower previsions that are only exchangeable. It should be clear that they correspond to exchangeable lower previsions for which the corresponding coherent lower prevision $\lpr_\kappa^N$ on $\gambles(\counts_\kappa^N)$ is \emph{strongly} rather than just weakly invariant with respect to all category permutations.  Of course, any justification for such models should be based on beliefs that there is permutation symmetry in the categories behind the process that generates the outcomes $\rv_k$ at different times $k$, and \emph{cannot be justified by mere ignorance about this process.}

\subsection{Updating exchangeable lower previsions: predictive inference}
Finally, let us discuss possible applications of the discussion in this paper to predictive inference. Assume that we have $n^*$ random variables $\rv_1$, \dots $\rv_{n^*}$, that may assume values in the set $\values_\kappa=\{1,\dots,\kappa\}$. We assume that these random variables are assessed to be exchangeable, in the sense that any coherent lower prevision that describes the available information about the values that the joint random variable $\sample{\rv}^*=(\rv_1,\dots,\rv_{n^*})$ assumes in $\values_\kappa^{n^*}$ should be exchangeable, i.e., strongly $\permuts_\kappa^{n^*}$-invariant. This requirement could be called \emph{pre-data exchangeability}. So we know from the previous section that such a coherent lower prevision must be of the form $\lpr=\lpr_\kappa^{n^*}(\pr^u(\cdot\vert\counts_\kappa^{n^*}))$, where $\lpr_\kappa^{n^*}$ is some coherent lower prevision on $\gambles(\counts_\kappa^{n^*})$. We shall assume that $\lpr_\kappa^{n^*}$ is a lower envelope of a set of coherent previsions $\solp_\kappa^{n^*}$ on $\gambles(\counts_\kappa^{n^*})$.
\par
Suppose we now observe the values $\sample{x}=(x_1,\dots,x_n)$ of the first $n$ random variables $\sample{X}=(\rv_1,\dots,\rv_n)$, where $1\leq n<n^*$. We ask ourselves how we should coherently update the belief model $\lpr$ to a new model $\lpr(\cdot\vert\sample{x})$ which describes our beliefs about the values of the remaining random variables $\sample{\rv}'=(\rv_{n+1},\dots,\rv_{n^*})$.  This is, generally speaking, the problem of \emph{predictive inference}. In order to make things as easy as possible, we shall assume that $\lpr(\{\sample{x}\})>0$, so our subject has some reason, prior to observing $\sample{x}$, to believe that this observation will actually occur, because she is willing to bet on its occurrence at non-trivial odds.
\par
Let us denote by $n'=n^*-n$ the number of remaining random variables, then we know that $\sample{\rv}'$ assumes values in $\values_\kappa^{n'}$, and $\lpr(\cdot\vert\sample{x})$ will be a lower prevision on $\gambles(\values_\kappa^{n'})$.
\par
We shall first look at the problem of updating the coherent prevision $\pr=\apr(\pr^u(\cdot\vert\counts_\kappa^{n^*}))$ for any coherent prevision $\apr$ in $\solp_\kappa^{n^*}$. So consider any gamble $g$ on $\values_\kappa^{n'}$.  It follows from coherence requirements (Bayes's rule) that the updated coherent prevision $\pr(\cdot\vert\sample{x})$ is given by
\begin{equation}\label{eq:bayes-1}
  \pr(g\vert\sample{x})
  =\frac{\pr(gI_{\sample{x}})}{\pr(I_{\sample{x}})}
  =\frac{\apr(\pr^u(gI_{\sample{x}}\vert\counts_\kappa^{n^*}))}
  {\apr(\pr^u(I_{\sample{x}}\vert\counts_\kappa^{n^*}))},
\end{equation}
where $I_{\sample{x}}(\sample{x}^*)=1$ if the first $n$ components of the vector $\sample{x}^*\in\values_\kappa^{n^*}$ are given by the vector $\sample{x}$, and zero otherwise. Observe, by the way, that by assumption, $\pr(I_{\sample{x}})\geq\lpr(I_{\sample{x}})=\lpr(\{\sample{x}\})>0$.
\par
Now for any $\cnt{m}^*$ in $\counts_\kappa^{n^*}$ we find that, with obvious notations,
\begin{equation}\label{eq:bayes-2}
  \pr^u(gI_{\sample{x}}\vert\cnt{m}^*)
  =\frac{1}{\nu(\cnt{m}^*)}
  \sum_{\cntf'(\sample{x}')+\cnt{m}=\cnt{m}^*}g(\sample{x}')
  =\frac{\nu(\cnt{m}^*-\cnt{m})}{\nu(\cnt{m}^*)}\pr^u(g\vert\cnt{m}^*-\cnt{m})
\end{equation}
where we let $\cnt{m}=\cntf(\sample{x})$,and where $\cntf'$ maps samples $\sample{x}'$ in $\values_\kappa^{n'}$ to their corresponding count vectors $\cntf'(\sample{x}')$ in $\counts_\kappa^{n'}$. Of course $\nu(\cnt{m}^*-\cnt{m})$ is non-zero only if $\cnt{m}^*\geq\cnt{m}$, or equivalently if $\cnt{m}^*-\cnt{m}\in\counts_\kappa^{n'}$, or in other words if it is possible to select $n$ balls of composition $\cnt{m}$ without replacement from an urn with composition $\cnt{m^*}$. In this expression, $\pr^u(\cdot\vert\cnt{m}')$ stands for the coherent prevision on $\gambles(\values_\kappa^{n'})$ whose probability mass is uniformly distributed over the $\permuts_\kappa^{n'}$-invariant atom $[\cnt{m}']$, for any $\cnt{m}'$ in $\counts_\kappa^{n'}$. Now for $g=1$ we find that
\begin{equation}\label{eq:bayes-3}
  \pr^u(I_{\sample{x}}\vert\cnt{m}^*)
  =\frac{\nu(\cnt{m}^*-\cnt{m})}{\nu(\cnt{m}^*)}
  =p(\cnt{m}\vert\cnt{m}^*)=:L_{\cnt{m}}(\cnt{m}^*)
\end{equation}
is the probability of observing a sample of size $n$ with composition $\cnt{m}$ by sampling without replacement from an urn with composition $\cnt{m}^*$. $L_{\cnt{m}}$ is the corresponding likelihood function on $\counts_\kappa^{n^*}$. We may as well consider $L_{\cnt{m}}$ as a likelihood function on $\counts_\kappa^{n'}$, and for any $\cnt{m}'$ in $\counts_\kappa^{n'}$ we let
\begin{equation*}
  L_{\cnt{m}}(\cnt{m}'):=L_{\cnt{m}}(\cnt{m}+\cnt{m}')
  =\frac{\nu(\cnt{m}')}{\nu(\cnt{m}+\cnt{m}')}
\end{equation*}
be the probability that there remain $n'$ balls of composition $\cnt{m}'$ after drawing (without replacement) $n$ balls of composition $\cnt{m}$ from an urn with $n^*$ balls. We may then rewrite Eq.~\eqref{eq:bayes-1}, using Eqs.~\eqref{eq:bayes-2} and~\eqref{eq:bayes-3}, as
\begin{equation}\label{eq:bayes-5}
  \pr(g\vert\sample{x})
  =\frac{\apr(L_{\cnt{m}}\pr^u(g\vert\counts_\kappa^{n'}))}{\apr(L_{\cnt{m}})}
  =\apr(\pr^u(g\vert\counts_\kappa^{n'})\vert\cnt{m}),
\end{equation}
where $\apr(L_{\cnt{m}})=\pr(I_{\sample{x}})>0$ by assumption, and $\apr(\cdot\vert\cnt{m})$ is the coherent prevision on $\gambles(\counts_\kappa^{n'})$ defined by
\begin{equation}\label{eq:bayes-6}
  \apr(h\vert\cnt{m})
  :=\frac{\apr(L_{\cnt{m}}h)}{\apr(L_{\cnt{m}})},
\end{equation}
for any gamble $h$ on $\counts_\kappa^{n'}$, i.e., $\apr(\cdot\vert\cnt{m})$ is the coherent prevision obtained after using Bayes's rule to update $\apr$ with the likelihood function $L_{\cnt{m}}$. This means that \emph{if $\apr$ is a belief model for the unknown composition of an urn with $n^*$ balls, then   $\apr(\cdot\vert\cnt{m})$ is the corresponding model for the unknown composition of the remaining $n'$ balls   in the urn, after $n$ balls with composition $\cnt{m}$ have been taken from it.}
\par
Now if we have a coherent lower prevision $\lpr_\kappa^{n^*}$ on $\gambles(\counts_\kappa^{n^*})$ that is a lower envelope of a set $\solp_\kappa^{n^*}$ of coherent previsions $\apr$, then coherence\footnote{This   follows from \citegen[Section~6.5]{walley1991} Generalised Bayes Rule.} tells us that the updated lower prevision $\lpr(\cdot\vert\sample{x})$ is precisely the lower envelope of the corresponding updated coherent previsions $\pr(\cdot\vert\sample{x})$, and consequently, using Eqs.~\eqref{eq:bayes-5} and~\eqref{eq:bayes-6}, we find that
\begin{equation}\label{eq:bayes-7}
  \lpr(g\vert\sample{x})
  =\lpr_\kappa^{n^*}(\pr^u(g\vert\counts_\kappa^{n'})\vert\cnt{m}),
\end{equation}
where $\lpr_\kappa^{n^*}(\cdot\vert\cnt{m})$ is the coherent lower prevision on $\gambles(\counts_\kappa^{n'})$ given by
\begin{equation}\label{eq:bayes-8}
  \lpr_\kappa^{n^*}(h\vert\cnt{m})
  :=\inf\set{\frac{\apr(L_{\cnt{m}}h)}{\apr(L_{\cnt{m}})}}
  {\apr\in\solp_\kappa^{n^*}}
  =\inf\set{\apr(h\vert\cnt{m})}{\apr\in\solp_\kappa^{n^*}},
\end{equation}
for any gamble $h$ on $\counts_\kappa^{n'}$. In other words, $\lpr_\kappa^{n^*}(\cdot\vert\cnt{m})$ is the coherent lower prevision obtained after using coherence (the so-called Generalised Bayes Rule) to update $\lpr_\kappa^{n^*}$ with the likelihood function $L_{\cnt{m}}$. This means again that \emph{if   $\lpr_\kappa^{n^*}$ is a belief model for the unknown composition of an urn with $n^*$ balls, then   $\lpr_\kappa^{n^*}(\cdot\vert\cnt{m})$ is the corresponding belief model for the unknown composition of the   remaining $n'$ balls in the urn, after $n$ balls with composition $\cnt{m}$ have been taken from it.}
\par
If we compare Eq.~\eqref{eq:bayes-7} with Eq.~\eqref{eq:exchangeability}, we see that the updated belief model $\lpr(\cdot\vert\sample{x})$ is still strongly $\permuts_\kappa^{n'}$-invariant,\footnote{See also   footnote~\ref{fn:conditioning}.} so there still is \emph{post-data exchangeability} for the remaining random variables $\sample{\rv}'=(\rv_{n+1},\dots,\rv_{n^*})$. Moreover, by looking at Eq.~\eqref{eq:bayes-1} and Eqs.~\eqref{eq:bayes-7} and~\eqref{eq:bayes-8}, we see that the updated (lower) previsions $\pr(\cdot\vert\sample{x})$ and $\lpr(\cdot\vert\sample{x})$ only depend on the observed sample $\sample{x}$ through the \emph{likelihood function} $L_{\cntf(\sample{x})}$. This tells us that this type of predictive inference satisfies the so-called \emph{likelihood principle}, and moreover that the count vector $\cnt{m}=\cntf(\sample{x})$, or more generally the map $\cntf$ is a \emph{sufficient statistic}.

\section{Conclusions}
\label{sec:conclusions}
We have tried to argue that there is a clear distinction between the symmetry of belief models, and models of beliefs of symmetry, and that both notions can be distinguished between when indecision is taken seriously, as is the case in \citegen{walley1991} behavioural theory of imprecise probabilities. Our present attempt to distinguish between these notions, and capture the distinction in a formal way, is inspired by \citegen[Chapter~9]{walley1991} discussion of the difference between permutable and exchangeable lower previsions, and \citegen{pericchi1991} discussion of `classes of reasonable priors' versus `reasonable classes of priors'.
\par
Indeed, there seems to be a difference of type between the two notions. The former (symmetry of models) is a property that belief models may have, and we may require, as a principle of rationality, or as a principle of `faithful modelling', that if the available evidence is symmetrical, then our corresponding belief models should be symmetrical too. A case in point is that of complete ignorance, where the `evidence' is completely symmetrical, and we may therefore require that corresponding belief model should be completely symmetrical too. This leads to the various principles discussed in Section~\ref{sec:complete-ignorance}, all of which seem to single out the vacuous belief model for representing complete ignorance, and which extend \citegen[Section~5.5]{walley1991} treatment of this matter.
\par
The latter notion (models of symmetry) is more properly related to a type of structural assessment: if a subject believes there is symmetry, how should she model that, and how should assessments of symmetry be combined with other assessments? We have tried to answer such questions in Sections~\ref{sec:strongly-invariant-lpr}, where we discuss the strongly invariant natural extension.
\par
It is well-known that if we only use Bayesian, or precise, probability models, requiring invariance of the probability measures with respect to all types of symmetry in the evidence may be impossible; examples were given by Boole, Bertrand and Fisher (see \citet{zabell1989} for discussion and references).  This has led certain researchers to abandon requiring the above-mentioned `faithfulness' of belief models, or to single out certain types of symmetry which are deemed to be better than others. We have tried to argue that this is unnecessary: the vacuous belief model has no such problems, and is symmetrical with respect to any transformation you care to name. And of course, our criticism of the Principle of Insufficient Reason is not new. Our ideas were heavily influenced by \citegen{walley1991} book on imprecise probabilities, whose Chapter~5 contains a wonderful overview of arguments against restricting ourselves to precise probability models. \citet{zabell1989a} also gives an excellent discussion of much older criticism, dating back to the middle of the 19th century. In particular, \citegen{ellis1844} \textit{ex nihilo nihil} --- you cannot make decisions or inferences based on complete ignorance --- finds a nice confirmation in the fact that the vacuous belief model captures complete indecision, and that updating a vacuous belief model leads to a vacuous belief model \citep[Section~6.6.1]{walley1991}. But what we have tried to do here is provide a framework and mathematical apparatus that allows us to better understand and discuss the problems underlying the Principle of Insufficient Reason, and more general problems of dealing with any type of symmetry in belief models.
\par
This study of symmetry in relation to belief models is far from being complete however, and our notions of weak and strong invariance may have to be refined, and perhaps even modified, as well as complemented by other notions of symmetry. It might for instance be of interest to study the notion of symmetry that captures the \emph{insufficient reason to strictly prefer} that is briefly touched upon near the end of Section~\ref{sec:symmetry-of-models}.  Also, we may seem more certain than we actually are about the appropriateness (in terms of having a sound behavioural justification and interpretation) of our notions of weak and (especially) strong invariance for random variables that may assume an infinite number of values. This is the point where our intuition deserts us, and where a number of interesting questions and problems leave us speechless. To name but one such problem, brought to the fore by the discussion in Section~\ref{sec:strongly-invariant-lpr}: for certain types of monoids, it is completely irrational to impose strong invariance (because doing so makes us subject to a sure loss). We can understand why this is the case for the monoid of all transformations, even on a finite set (Theorem~\ref{theo:no-complete-strong-invariance}). But why, for instance, are there no (strongly) permutation invariant coherent (lower) previsions on the set of natural (and \textit{a fortiori} real) numbers? Why are we (consequently) reduced to using (strong) shift or translation invariance of coherent (lower) previsions when we want to try and capture the idea of a uniform distribution on the set of natural (or real) numbers? And even then, why, as is hinted at in footnote~\ref{fn:congomerability}, are there situations where updating a (strongly) shift-invariant coherent (lower) prevision produces a sure loss?  Are there appropriately weakened versions of our strong invariance condition that avoid these problems?

\section*{Acknowledgements}
This paper has been partially supported by research grant G.0139.01 of the
Flemish Fund for Scientific Research (FWO), and projects MTM2004-01269,
TSI2004-06801-C04-01.


\bibliographystyle{plainnat}

\end{document}